\documentclass[11pt]{article}
\usepackage[english]{babel}
\usepackage[a4paper,scale=0.75,hcentering]{geometry}
\usepackage{amsmath,amsthm,amssymb,mathtools}
\usepackage{aliascnt,needspace,xcolor}
\usepackage[hyperfootnotes=false]{hyperref}
\ifdefined\pdfinfoomitdate\pdfinfoomitdate=1\fi
\hypersetup{colorlinks=true,linkcolor=blue,urlcolor=red,citecolor=blue,
 pdftitle={Gradient stability for the Sobolev inequality with optimal dimensional dependence},
 pdfauthor={Yuxuan Zhou and Wenming Zou}}
\usepackage[nameinlink,capitalise,sort]{cleveref}
\crefname{equation}{}{}
\crefname{enumi}{}{}
\allowdisplaybreaks
\newcounter{results}[section]

\theoremstyle{plain}
\newaliascnt{theorem}{results}
\newtheorem{theorem}[theorem]{Theorem}
\aliascntresetthe{theorem}
\newaliascnt{proposition}{results}
\newtheorem{proposition}[proposition]{Proposition}
\aliascntresetthe{proposition}
\newaliascnt{lemma}{results}
\newtheorem{lemma}[lemma]{Lemma}
\aliascntresetthe{lemma}
\crefname{theorem}{Theorem}{Theorems}
\crefname{proposition}{Proposition}{Propositions}
\crefname{lemma}{Lemma}{Lemmas}
\numberwithin{equation}{section}
\newcommand{\R}{\mathbb R}

\makeatletter
\DeclarePairedDelimiter{\@tmpabs}{\lvert}{\rvert}
\newcommand{\@absstar}[1]{{\@tmpabs*{#1}}}
\newcommand{\@absnostar}[2][]{{\@tmpabs[#1]{#2}}}
\newcommand{\abs}{\@ifstar\@absstar\@absnostar}
\DeclarePairedDelimiter{\@tmpnorm}{\lVert}{\rVert}
\newcommand{\@normstar}[1]{{\@tmpnorm*{#1}}}
\newcommand{\@normnostar}[2][]{{\@tmpnorm[#1]{#2}}}
\newcommand{\norm}{\@ifstar\@normstar\@normnostar}
\makeatother
\title{\bfseries Gradient stability for the Sobolev inequality with optimal dimensional dependence%
\thanks{This work is supported by National Key R\&D Program of China (Grant 2023YFA1010001) and NSFC (Grant 12571123).
E-mail addresses: \texttt{zhou-yx22@mails.tsinghua.edu.cn} (Zhou),
\texttt{zou-wm@mail.tsinghua.edu.cn} (Zou).}}

\author{{\bfseries Yuxuan Zhou, Wenming Zou}\\[1ex]
{\footnotesize\itshape Department of Mathematical Sciences, Tsinghua University,}\\
{\footnotesize\itshape Beijing 100084, China.}}
\date{}
\begin{document}
\maketitle
\begin{abstract}
Figalli and Zhang~\cite{FZ} established gradient stability
for the sharp Sobolev inequality with the optimal exponent
$r_p=\max\{2,p\}$. We prove that, for every fixed $p>1$,
the corresponding optimal stability constant $\mathsf C_{n,p}$
admits effective upper and lower bounds of order $n^{-1}$
for all integers $n>p$. In particular, the dimensional
dependence $n^{-1}$ is optimal. Our result extends the
optimal dimensional dependence established for $p=2$
by Dolbeault, Esteban, Figalli, Frank, and Loss~\cite{DEFLL}
to every $p>1$.

\medskip
\noindent\textbf{Key words:} Sobolev inequality, quantitative stability,
gradient stability, optimal dimensional dependence, $p$-Laplacian.

\medskip
\noindent\textbf{2020 Mathematics Subject Classification:}
Primary 46E35, 26D10; Secondary 35J92, 35B35.
\end{abstract}

\tableofcontents

\section{Introduction}\label{sec:introduction}

Let $n\geq2$ and $1<p<n$. Denote by
$\dot W^{1,p}(\mathbb R^n)$ the completion of
$C_c^\infty(\mathbb R^n)$ with respect to the norm
$\|\nabla u\|_{L^p(\mathbb R^n)}$. Set
\begin{equation}\label{eq:exponents}
 p'=\frac p{p-1},\qquad q=\frac{np}{n-p},\qquad
 r_p=\max\{2,p\},\qquad q-p=\frac{p^2}{n-p}.
\end{equation}
The sharp Sobolev inequality is
\begin{equation}\label{eq:sobolev}
 S_{n,p}\norm*{u}_{L^q(\mathbb R^n)}\leq\norm*{\nabla u}_{L^p(\mathbb R^n)},
\end{equation}
where $S_{n,p}$ is the optimal Sobolev constant.
Let $\mathcal M_{n,p}$ be the family of functions
\begin{equation}\label{eq:manifold}
 \mathcal M_{n,p}=\left\{
 A\lambda^{(n-p)/p}
 (1+\abs*{\lambda(x-y)}^{p'})^{-(n-p)/p}:
 A\in\mathbb R\setminus\{0\},
 \lambda>0,\ y\in\mathbb R^n
 \right\}.
\end{equation}
It is known~\cite{Aubin,Talenti,CNV} that
$\mathcal M_{n,p}$ is precisely the set of all
nonzero extremals of \eqref{eq:sobolev}.

In this paper, we study quantitative stability for
the Sobolev inequality \eqref{eq:sobolev}.
For $u\in\dot W^{1,p}(\mathbb R^n)\setminus\{0\}$,
define
\begin{align}
 E(u)&=\norm*{\nabla u}_{L^p(\mathbb R^n)}^p,&
 F(u)&=E(u)-S_{n,p}^p\norm*{u}_{L^q(\mathbb R^n)}^p,&
 \mathcal D(u)&=\frac{F(u)}{E(u)},
 \label{eq:deficit}\\
 \operatorname{dist}_{\nabla,p}(u)
 &=\inf_{V\in\mathcal M_{n,p}}
       \norm*{\nabla(u-V)}_{L^p(\mathbb R^n)},&
 d(u)&=\frac{\operatorname{dist}_{\nabla,p}(u)}
             {E(u)^{1/p}}.
 &&\label{eq:distance}
\end{align}
Here $\mathcal D(u)$ is the normalized Sobolev
deficit, and $d(u)$ is the relative gradient
distance to $\mathcal M_{n,p}$.
Both quantities are invariant under multiplication
by nonzero constants, translations, and critical
dilations.

The question of quantitative stability for \eqref{eq:sobolev} was raised by Brezis and
Lieb~\cite{BL}. They asked whether, for $p=2$,
the Sobolev deficit controls a suitable distance
to the manifold of extremals.
Bianchi and Egnell~\cite{BE} answered this question
by proving that, for $p=2$,
\[
 \mathcal D(u)\geq c(n)d(u)^2.
\]
Their result is optimal both in the choice of
the gradient distance and in the exponent $2$.
The proof uses spectral estimates near the
extremal manifold and a compactness argument.
The first stability result covering the full range
$1<p<n$ was obtained by Cianchi, Fusco, Maggi,
and Pratelli~\cite{CFMP}. Their estimate,
which combines symmetrization and mass transport,
can be written in the form
\[
 \mathcal D(u)\geq c(n,p)
 \left(
 \inf_{V\in\mathcal M_{n,p}}
 \frac{\norm*{u-V}_{L^q(\mathbb R^n)}}{\norm*{u}_{L^q(\mathbb R^n)}}
 \right)^{\alpha_{n,p}},
\]
for an explicit exponent $\alpha_{n,p}>0$.

Figalli and Neumayer~\cite{FN} first obtained
gradient stability for $p>2$, proving
\[
 \mathcal D(u)\geq c(n,p)d(u)^{\beta_{n,p}}
\]
for a suitable exponent $\beta_{n,p}>0$.
Their argument combines weighted spectral estimates
with the $L^q$ stability result of~\cite{CFMP}.
Neumayer~\cite{Neumayer} subsequently extended
gradient stability to the full range $1<p<n$,
using a uniform-convexity argument.
The exponents in these results were not optimal.
Finally, Figalli and Zhang~\cite{FZ} proved that,
for every $1<p<n$,
\begin{equation}\label{eq:fz-stability}
 \mathcal D(u)\geq c(n,p)d(u)^{r_p},
\end{equation}
where the exponent $r_p=\max\{2,p\}$ is optimal.
Their proof combines nonlinear vector inequalities,
weighted spectral analysis, and compactness arguments.

Quantitative stability has also been studied for
isoperimetric and Wulff inequalities
\cite{FMPiso,FMPmass,FigalliZhangWulff},
the Brunn--Minkowski, Pr\'ekopa--Leindler, and
Borell--Brascamp--Lieb inequalities
\cite{FigalliJerisonBM,FvHTBM,FvHTBMtransport,
BoroFigRamosPL,FigalliRamosPL,FvHTBBL},
and various functional inequalities, including
Gagliardo--Nirenberg, Caffarelli--Kohn--Nirenberg,
Hardy--Littlewood--Sobolev, and logarithmic
Sobolev inequalities
\cite{CarlenFigalliGNS,NguyenGNS,
CazacuFlynnLamLuCKN,CarlenDuality,
FathiIndreiLedouxLSI,EldanLehecShenfeldLSI}.
We refer to these works, as well as
\cite{FMPanisotropic,BrigatiDolbeaultSimonovLSI},
for further results and references.

A further question is to determine the optimal
constant in \eqref{eq:fz-stability}. We define
\begin{equation}\label{eq:bestconstant}
 \mathsf C_{n,p}=
 \inf_{\substack{
 u\in\dot W^{1,p}(\mathbb R^n)\setminus\{0\}\\
 d(u)>0}}
 \frac{\mathcal D(u)}{d(u)^{r_p}}.
\end{equation}
For $p=2$, K\"onig~\cite{KoenigSharp} proved
 the strict upper bound
\[
 \mathsf C_{n,2}<\frac{4}{n+4}
 \qquad(n\geq3).
\]
In subsequent work~\cite{KoenigMinimizer},
he also showed that
\[
 \mathsf C_{n,2}<2(1-2^{-2/n})
\]
and that the infimum defining
$\mathsf C_{n,2}$ is attained.

Obtaining effective lower bounds is more difficult.
Although Bianchi and Egnell~\cite{BE} and
Figalli and Zhang~\cite{FZ} proved the existence
of positive stability constants, their arguments
rely on compactness and do not yield explicit
lower bounds for these constants. Moreover, as
observed in~\cite{DEFLL}, the classical local
estimates for $p=2$ do not give the optimal
dimensional dependence.

Recently, Dolbeault, Esteban, Figalli, Frank,
and Loss~\cite{DEFLL} proved that
\[
 \mathsf C_{n,2}\geq\frac{\beta}{n}
 \qquad(n\geq3),
\]
where $\beta>0$ is explicit and independent
of $n$. Their proof uses refined local estimates
based on spectral analysis near the manifold of
extremals, together with competing symmetries
and continuous Steiner symmetrization to pass
from local to global stability.

Related results on stability constants for
fractional and higher-order Sobolev inequalities,
fractional Sobolev trace inequalities,
Hardy--Littlewood--Sobolev inequalities,
Caffarelli--Kohn--Nirenberg inequalities,
and Sobolev inequalities in other settings
can be found in
\cite{CLTExplicitHLS,CLTFractionalDim,CLT,
CLTHeisenberg,CPPCones,ZZZTrace,LiWeiCKN}
and the references therein.

In this paper, we extend the optimal
$n^{-1}$ dimensional dependence to every fixed
$p>1$. Our main result is the following.

\Needspace{12\baselineskip}
\begin{theorem}\label{thm:main}\label{thm:full}
For every $p>1$, there exists a constant $c_p>0$
depending only on $p$ such that, for every integer
$n>p$,
\begin{equation}\label{eq:mainconstants}
 \frac{c_p}{n}
 \leq\mathsf C_{n,p}
 \leq2^{r_p/p}(1-2^{-p/n})
 \leq\frac{2^{r_p/p}p\log 2}{n},
 \qquad r_p=\max\{2,p\}.
\end{equation}
The constant $c_p$ in the lower bound can be
chosen explicitly by tracing the estimates in
the proof.
\end{theorem}

We first recall the main ideas of the proofs by
Dolbeault et al.~\cite{DEFLL} and Figalli and
Zhang~\cite{FZ}. For $p=2$, Dolbeault et al.
establish dimension-sharp local stability using
the conformal formulation of the Sobolev inequality
on the sphere, a spectral gap, and a decomposition
of perturbations according to their size.
They then use competing symmetries and continuous
Steiner symmetrization to pass from local to
global stability for nonnegative functions,
before extending the result to sign-changing
functions.

For general $1<p<n$, Figalli and Zhang
\cite{FZ} proved gradient stability with the
optimal exponent $r_p=\max\{2,p\}$.
Their proof uses nonlinear vector inequalities
and weighted spectral estimates to control
perturbations near the manifold of Sobolev
extremals. They also use a variational argument
to obtain the required orthogonality conditions.
Compactness arguments are then used to establish
the nonlinear local estimates and reduce the
global problem to the local one.

Our proof follows the local-to-global approach
of \cite{DEFLL} and uses nonlinear estimates
and weighted spectral analysis from \cite{FZ}.
However, the local argument in \cite{DEFLL}
relies on the Hilbert space structure available
when $p=2$ and cannot be applied directly to
general $p$. Figalli and Zhang overcome the
corresponding nonlinear difficulties for fixed
$n,p$, but their compactness arguments do not
provide effective stability constants with
optimal dimensional dependence.
We therefore develop quantitative local
estimates that yield a coefficient of order
$n^{-1}$.

We first consider nonnegative functions in sufficiently large dimensions. Following the decomposition of perturbations in \cite{DEFLL}, we treat bounded, intermediate, and large values separately. A coupled estimate for the radial and angular terms (\cref{lem:core}), together with estimates for small-support and large-value perturbations, allows us to control the nonlinear errors while retaining a coefficient of order $n^{-1}$. Applying \cref{prop:iv-invariant} to the resulting
orthogonal estimate gives \cref{prop:local},
with a neighborhood size depending only on $p$.

For the remaining finitely many dimensions, we use Zhang's spectral description of the linearized operator \cite{ZhangSpectrum}. We combine the spectral gap in
\cref{prop:stage4-spectral-toolkit}
with estimates using finitely many
eigenfunctions to control the nonlinear
remainder. For $p\ge2$, this gives \cref{thm:super-local}. When $1<p<2$, we also use truncation and finite spectral approximation (\cref{prop:short-scalar-capture}), since a general perturbation need not belong to the quadratic form domain. This yields \cref{thm:short-subquadratic-local}. Both estimates allow sign-changing perturbations. We then apply \cref{prop:iv-invariant} to obtain
local stability estimates in terms of the
gradient distance to the full manifold of
Sobolev extremals.

We next establish global stability for nonnegative functions. Following the radial transport method of Cianchi, Fusco, Maggi, and Pratelli \cite{CFMP}, we first prove quantitative approximation by a centered Sobolev extremal for radially nonincreasing functions (\cref{prop:both-radial}). We then use finitely many polarizations \cite{BF}, connected by paths continuous in the gradient norm, to treat general nonnegative functions. The vector Clarkson inequalities \cite{Neumayer} show that the endpoint is close to an extremal in the gradient norm. Since the paths preserve the gradient energy, continuity allows us to reach the boundary of the local stability neighborhood and apply the local estimate. This gives \cref{prop:exact-crossing}.

Finally, we extend the global estimate to arbitrary functions by separating their positive and negative parts. Combining the dimension-uniform estimate with the finitely many remaining dimensions gives the lower bound in \eqref{eq:mainconstants}, with $c_p$ defined in \eqref{eq:finalcp}. The upper bound follows from two widely separated Sobolev extremals (\cref{prop:upper}).

\section{Preliminaries}\label{sec:preliminaries}

In this section, we collect the estimates needed for
 the proof. In Subsection~\ref{sec:bubble-identities},
we introduce coordinates associated with Sobolev
extremals and derive the weighted representation
\eqref{eq:local-representation}. In
Subsection~\ref{sec:basic-tools}, we establish
a vector Bregman bound and record estimates for
the coordinate operators. Finally, in
Subsection~\ref{sec:modulation}, we choose a nearby
Sobolev extremal satisfying the required
orthogonality conditions and show how to pass
from an orthogonal local estimate to a stability
estimate involving the gradient distance to
the manifold of extremals.

Unless otherwise specified, all norms and integrals are taken over
$\mathbb R^n$ with respect to Lebesgue measure. For
$1\leq r\leq\infty$, we write
$\|f\|_r=\|f\|_{L^r(\mathbb R^n)}$.
Norms and integrals with respect to other measures are indicated
explicitly. For simplicity, $C_p$ denotes a positive constant
depending only on $p$ and may change from line to line; additional
dependencies are indicated by subscripts.

\subsection{Sobolev extremals and weighted identities}\label{sec:bubble-identities}
Set
\begin{equation}\label{eq:standard-bubble}
v(x)=(1+\abs*{x}^{p'})^{-(n-p)/p},\qquad
E_v=E(v),\qquad V_q=\int_{\R^n}v^q\,dx.
\end{equation}
The sharp Sobolev constant is
\[
 S_{n,p}=\sqrt\pi\,n^{1/p}
 \left(\frac{n-p}{p-1}\right)^{(p-1)/p}
 \left(\frac{\Gamma(n/p)\Gamma(1+n-n/p)}
 {\Gamma(1+n/2)\Gamma(n)}\right)^{1/n}.
\]
Equivalently, $S_{n,p}=E_v^{1/p}/V_q^{1/q}$.

Define the probability measures
\begin{equation}\label{eq:local-measures}
d\mu=V_q^{-1}v^q\,dx,\qquad
d\nu=E_v^{-1}\abs*{\nabla v}^p\,dx.
\end{equation}
Set
\[
t=\frac{\abs*{x}^{p'}}{1+\abs*{x}^{p'}},\qquad
\tau=\frac{p-1}{p},\qquad a(t)=t(1-t),\qquad
d_*=\tau(1-\tau),\qquad k(t)=t^{1/p}(1-t)^\tau.
\]
Write $\omega=x/\abs*{x}$ off the origin, and let $\nabla_S$ denote the
gradient on the unit sphere. Let
$d\sigma=dS/\abs*{\mathbb S^{n-1}}$ be normalized area measure.
Direct radial integration gives
\begin{equation}\label{eq:beta-measures}
\begin{gathered}
 d\mu=d\mu_t\,d\sigma,\qquad d\nu=d\nu_t\,d\sigma,\\
 d\mu_t(t)=\frac{t^{n\tau-1}(1-t)^{n/p-1}}
                    {\mathrm B(n\tau,n/p)}\,dt,\\
 d\nu_t(t)=\frac{t^{n\tau}(1-t)^{n/p-2}}
                    {\mathrm B(n\tau+1,n/p-1)}\,dt.
\end{gathered}
\end{equation}
Both radial measures have mass $1$, as does $\sigma$.
Set $m(t)=d\mu_t/dt$.

Write $e=\nabla v/\abs*{\nabla v}=-\omega$ off the origin.
Define the operator $\mathsf Z$ by
\begin{equation}\label{eq:ground-gradient}
\mathsf Zb=-\frac{p(1-t)}{n-p}b_te
 +\frac{p-1}{(n-p)t}\nabla_Sb.
\end{equation}
The identity
\begin{equation}\label{eq:normalized-gradient-product}
\frac{\nabla(vb)}{\abs*{\nabla v}}
=be+\mathsf Zb
\end{equation}
follows from the definition of $\mathsf Z$.

For $s\geq0$ and nonnegative functions $f$, define
\begin{equation}\label{eq:weighted-remainder}
G_s(z)=\abs*{se+z}^p-s^p-ps^{p-1}e\cdot z,
\qquad G_0(z)=\abs*{z}^p,
\qquad \mathcal R(f)=\int G_f(\mathsf Zf)\,d\nu.
\end{equation}
For $s\geq0$, define $D_{p,q}$ by
\begin{equation}\label{eq:local-potential}
D_{p,q}(s)=\frac pq s^q-s^p+1-\frac pq,
\qquad \epsilon=q-p=\frac{p^2}{n-p}.
\end{equation}
Since $D_{p,q}(1)=0$ and
$D_{p,q}'(s)=ps^{p-1}(s^\epsilon-1)$,
we have $D_{p,q}(s)\geq0$ for all $s\geq0$.

For $1<p<n$, Hardy's inequality, first on $C_c^\infty(\R^n)$
and then by density, gives for every $h\in\dot W^{1,p}$
\begin{equation}\label{eq:hardy-basic}
 \norm*{h/\abs*{x}}_p\leq\frac p{n-p}\norm*{\nabla h}_p.
\end{equation}

We call a function $b$ admissible if
$b\in W^{1,p}_{\rm loc}(\R^n)$,
$vb\in\dot W^{1,p}(\R^n)$ and
$b,\mathsf Zb\in L^p(\nu)$. By Hardy's inequality
\eqref{eq:hardy-basic}, the conditions
$b,\mathsf Zb\in L^p(\nu)$ follow from
$vb\in\dot W^{1,p}(\R^n)$.
In particular, every admissible function
belongs to $L^q(\mu)$.

The bubble $v$ satisfies
\begin{equation}\label{eq:bubble-euler-lagrange}
-\operatorname{div}(\abs*{\nabla v}^{p-2}\nabla v)
 =n\left(\frac{n-p}{p-1}\right)^{p-1}v^{q-1}.
\end{equation}
For every admissible function $b$, we have
\begin{equation}\label{eq:local-picone-signed}
\int\abs*{b}^p\,d\mu
 =\int\abs*{b}^p\,d\nu
  +p\int\abs*{b}^{p-2}b\,e\cdot\mathsf Zb\,d\nu.
\end{equation}

For nonnegative admissible $f$,
\begin{align}
\mathcal R(f)&=\frac{E(vf)}{E_v}-\int f^p\,d\mu,
 \label{eq:local-picone}\\
\frac{F(vf)}{E_v}&=\mathcal R(f)-P_q(f),
\qquad
P_q(f)=\left(\int f^q\,d\mu\right)^{p/q}-\int f^p\,d\mu
 \leq\int D_{p,q}(f)\,d\mu.
 \label{eq:local-representation}
\end{align}
The last inequality is the tangent-line inequality for $x^{p/q}$ at
$x=1$; it does not require a normalization of $f$.

To justify \eqref{eq:local-picone-signed}, let $b$
be admissible and set
$b_j=\min\{j,\max\{-j,b\}\}$. Then the test $v\abs*{b_j}^p$
belongs to $\dot W^{1,p}$, since
\[
\abs*{\nabla(v\abs*{b_j}^p)}
 \leq j^{p-1}(\abs*{b_j\nabla v}+p\abs*{v\nabla b_j}),
\qquad v\abs*{b_j}^p\in L^q.
\]
Testing \eqref{eq:bubble-euler-lagrange} with $v$ fixes its coefficient as
$E_v/V_q$; testing with $v\abs*{b_j}^p$ gives
\cref{eq:local-picone-signed} for $b_j$.
The cross term is dominated by
$\abs*{b}^{p-1}\abs*{\mathsf Zb}\in L^1(\nu)$, so dominated convergence
proves the identity for $b$.
For nonnegative admissible $f$, we use
\eqref{eq:normalized-gradient-product} and expand
$|fe+\mathsf Zf|^p$. Together with
\eqref{eq:local-picone-signed}, this gives
\eqref{eq:local-picone}. Since $p|b|^{p-2}b$ is continuous at
zero for $p>1$, the test does not require a positive lower bound
on $f$ or additional weighted quadratic integrability.

For $h=v\varphi$ put $\delta=\norm*{\nabla h}_p/E_v^{1/p}$.
By \eqref{eq:sobolev} and \eqref{eq:hardy-basic}, we have
\begin{equation}\label{eq:multiplier-bounds}
\norm*{\varphi}_{L^q(\mu)}\leq\delta,\qquad
\norm*{\varphi}_{L^p(\nu)}\leq p'\delta,\qquad
\norm*{\mathsf Z\varphi}_{L^p(\nu)}\leq(1+p')\delta.
\end{equation}
The factor $(n-p)/(p-1)$ in $|\nabla v|/v$ cancels the Hardy
constant $p/(n-p)$ in \eqref{eq:hardy-basic}.

For any $1$-Lipschitz function $T$ with $T(0)=0$,
the Sobolev chain rule gives
$\mathsf ZT(\varphi)=T'(\varphi)\mathsf Z\varphi$
almost everywhere. The value of $T'$ at points
where $T$ is not differentiable can be chosen
arbitrarily.
Thus $T(\varphi)$ is admissible and
\begin{equation}\label{eq:truncated-gradient-bound}
\frac{\norm*{\nabla(vT(\varphi))}_p}{E_v^{1/p}}
 \leq (1+2p')\delta.
\end{equation}

For a bounded admissible function $g$, we may also
 test \eqref{eq:bubble-euler-lagrange} with $vg^2$.
This gives
\begin{equation}\label{eq:bubble-quadratic-identity}
\int g^2\,d\mu
=\int g^2\,d\nu
+2\int g\,e\cdot\mathsf Zg\,d\nu.
\end{equation}

\subsection{Basic inequalities and coordinate operators}\label{sec:basic-tools}

For $\ell>1$ and vectors $a,b$, write
\begin{equation}\label{eq:vector-bregman}
 B_\ell(a,b)=|a+b|^\ell-|a|^\ell-\ell|a|^{\ell-2}a\cdot b,
\end{equation}
where the linear term is zero at $a=0$. Thus $G_s(z)=B_p(se,z)$.
A lower bound of this type follows from
\cite[Lemma~2.1]{FZ}. Here we give explicit
constants that will be used later.

\begin{lemma}\label{lem:eff-bregman}
For every $p>1$, there exists $\gamma_p>0$ such that,
for every $s\geq0$ and $z\in\R^n$,
\begin{equation}\label{eq:full-bregman-bound}
G_s(z)\geq\gamma_p(s+\abs*{z})^{p-2}\abs*{z}^2.
\end{equation}
The right-hand side is understood as zero when
$s=\abs*{z}=0$. One may take
$\gamma_p=p(p-1)/2$ for $1<p<2$, and
$\gamma_p=(p/8)24^{2-p}$ for $p\geq2$.
\end{lemma}

\begin{proof}
For $1<p<2$, the smallest eigenvalue of the Hessian
of $x\mapsto |x|^p$ is $p(p-1)|x|^{p-2}$ for $x\ne0$.
Taylor's formula along the segment $se+tz$,
$0\leq t\leq1$, together with
$|se+tz|\leq s+|z|$, gives
\eqref{eq:full-bregman-bound} with
$\gamma_p=p(p-1)/2$.
The case where the segment passes through the origin
follows by approximation.

For $p\ge2$, Taylor's formula gives
\[
 G_s(z)\ge p|z|^2\int_0^1(1-t)|se+tz|^{p-2}\,dt.
\]
If $|z|\le s/2$, the integral is at least
$\tfrac12((s+|z|)/3)^{p-2}$. Otherwise consider $0\le t\le1/2$.
The set where $|se+tz|<|z|/8$ has length at most $1/4$, as follows by
projection onto $z/|z|$. On its complement, of length at least $1/4$,
$1-t\ge1/2$ and $|z|>(s+|z|)/3$. This gives coefficient
$(p/8)24^{2-p}$. The cases $s=0$ or $z=0$ are immediate.
\end{proof}

\medskip
\noindent\emph{Coordinate operators.}
Write $\alpha=n\tau$, $\beta=n/p$, and
$\kappa=\min\{\tau,1-\tau\}$, with $a,m,d_*$ as in
\cref{sec:bubble-identities}. In the remainder of this subsection, assume $n\geq\max\{4,\kappa^{-1}\}$.
Let $A$ and $B_S$ be the nonnegative
self-adjoint operators associated with the closed forms
\begin{equation}\label{eq:coordinate-forms}
 \mathcal E_A(f,g)=\int_0^1 a f'g'\,d\mu_t,
 \qquad
 \mathcal E_{B_S}(f,g)=\int_{\mathbb S^{n-1}}
                           \nabla_S f\cdot\nabla_S g\,d\sigma.
\end{equation}
The first form is the closure of its restriction to $C^\infty([0,1])$,
while the second is the
usual $H^1$ form on the sphere.
For smooth functions, the operators are given by
\begin{equation}\label{eq:coordinate-operator-formulas}
 A=-a\partial_{tt}-(\alpha-nt)\partial_t,
 \qquad B_S=-\Delta_S.
\end{equation}

Their spectra are
\begin{equation}\label{eq:coordinate-spectra}
 \lambda_j(A)=j(j+n-1),\qquad
 \lambda_\ell(B_S)=\ell(\ell+n-2),\qquad j,\ell=0,1,\ldots.
\end{equation}
Consequently the first two positive spectral levels are $n,2(n+1)$ for
$A$ and $n-1,2n$ for $B_S$.

We use $A^{-1}$ and $B_S^{-1}$ on the
orthogonal complements of the constants
in their respective $L^2$ spaces.
The operator
\begin{equation}\label{eq:coordinate-resolvent}
 T=n(n+A)^{-1}
\end{equation}
is positive, preserves the mean, and is a contraction on $L^\infty$.
On the mean-zero subspace of $L^2(\mu_t)$,
\begin{equation}\label{eq:coordinate-resolvent-smoothing}
 \norm*{T}_{2\to2}\leq\frac12,
 \qquad \norm*{A^{1/2}T}_{2\to2}\leq\frac{\sqrt n}{2}.
\end{equation}
On the product space the two coordinate operators, their spectral
projections, and their resolvents commute.

The density satisfies
\begin{equation}\label{eq:coordinate-density-derivative}
 (a m)'=(\alpha-nt)m=n(\tau-t)m.
\end{equation}
By \eqref{eq:coordinate-density-derivative},
integration by parts in the first form
in \eqref{eq:coordinate-forms} gives
the expression for $A$ in
\eqref{eq:coordinate-operator-formulas}. The boundary term
vanishes since
$am=t^\alpha(1-t)^\beta/\mathrm B(\alpha,\beta)$.
The form is closable. Indeed, if $f_j\to0$
in $L^2(\mu_t)$ and $f_j'$ converges in
$L^2(a\,d\mu_t)$, then its limit vanishes
in the distributional sense on every
compact subinterval of $(0,1)$.
Polynomials are also dense in the form domain.
Indeed, we can approximate the derivative of a
smooth function uniformly by polynomials and
then integrate, fixing one value.

Under $x=1-2t$, the weight is a constant multiple of
$(1-x)^{\alpha-1}(1+x)^{\beta-1}$, so this closure is the Jacobi
realization in \cite[Section~2.7.4, pp.~113--116]{BGL}.
Together with the sphere spectrum in
\cite[Section~2.2.3, pp.~85--86]{BGL}, this gives
\cref{eq:coordinate-spectra}.
The eigenfunctions corresponding to the first
positive eigenvalues are $t-\tau$ for $A$
and $\omega_1,\ldots,\omega_n$ for $B_S$.

Both forms decrease under $1$-Lipschitz functions that vanish at zero.
Their semigroups and
normalized resolvents are therefore positive, preserve constants and
contract $L^\infty$ by
\cite[Sections~3.1.4--3.1.6, pp.~125--130]{BGL}.
By \eqref{eq:coordinate-resolvent},
$T=n\int_0^\infty e^{-ns}e^{-sA}\,ds$.
The spectral gap $\lambda\geq n$ on the mean-zero subspace
gives \cref{eq:coordinate-resolvent-smoothing}; tensor products
of the coordinate eigenbases give commutation.

\medskip
\noindent\emph{Curvature and entropy estimates.}
For $D=A$ or $B_S$, set
$\Gamma_A(f)=a\abs*{f'}^2$ and
$\Gamma_{B_S}(f)=\abs*{\nabla_Sf}^2$.
Set $\theta=2\arcsin\sqrt t$, so that $dt/d\theta=\sqrt a$.
The Jacobi form becomes $\int\abs*{f_\theta}^2\,d\mu_t$, and the density
in this coordinate is a constant multiple of
\[
 \sin^{2\alpha-1}(\theta/2)\cos^{2\beta-1}(\theta/2),
 \qquad 0<\theta<\pi.
\]
Writing the density as $e^{-U}$ gives
\begin{equation}\label{eq:coordinate-curvature}
 U''(\theta)=
 \frac{2\alpha-1}{4\sin^2(\theta/2)}
 +\frac{2\beta-1}{4\cos^2(\theta/2)}
 \geq\frac{n-1}{2}\geq\frac n4.
\end{equation}
Also $U''\geq n\kappa$, since $2\alpha-1,2\beta-1\geq n\kappa$
and $\sin^{-2}(\theta/2)+\cos^{-2}(\theta/2)\geq4$.

For $L=-D$ and
$\Gamma_2(f)=\frac12L\Gamma(f)-\Gamma(f,Lf)$, direct differentiation
and the sphere Bochner formula give
\[
 \Gamma_2^A(f)=\abs*{f_{\theta\theta}}^2+U''\abs*{f_\theta}^2,\qquad
 \Gamma_2^{B_S}(f)=\abs*{\nabla_S^2f}^2+(n-2)\abs*{\nabla_Sf}^2.
\]
Thus the curvature bounds are $K=n\kappa$ for $A$ and $K=n-2$
for $B_S$, and both allow $K=n/4$.

For a smooth function of $t$, $f_\theta=O(\theta)$ at zero and
$O(\pi-\theta)$ at $\pi$. Since the density vanishes with positive powers at both endpoints,
the boundary terms in the curvature and entropy calculations vanish.
We first work with polynomials and finite sums of spherical harmonics,
then extend the estimates by form closure and bounded approximation.

For $(D,\pi)=(A,\mu_t)$ or $(B_S,\sigma)$, we have
\begin{equation}\label{eq:coordinate-lsi}
 \operatorname{Ent}_\pi(f^2)
       \leq\frac8n\int\Gamma_D(f)\,d\pi,
       \qquad D=A\ \text{or}\ B_S.
\end{equation}
Here $\operatorname{Ent}_\pi(H)=
\int H\log(H/\int H\,d\pi)\,d\pi$ for $H\geq0$.
The corresponding semigroup $P_s=e^{-sD}$ obeys
\begin{equation}\label{eq:coordinate-hypercontractivity}
 \norm*{P_{s/2}z}_{L^2(\pi)}\leq\norm*{z}_{L^{1+e^{-ns/4}}(\pi)},\qquad s>0.
\end{equation}
On $\mu=\mu_t\otimes\sigma$ one also has
\begin{equation}\label{eq:coordinate-product-lsi}
 \operatorname{Ent}_\mu(F^2)
 \leq\frac8n\int\big(a\abs*{F_t}^2+\abs*{\nabla_S F}^2\big)\,d\mu.
\end{equation}

The logarithmic Sobolev inequalities
\eqref{eq:coordinate-lsi} and
\eqref{eq:coordinate-product-lsi}
hold on the corresponding closed form domains.

The curvature criterion \cite[Proposition~5.7.1, p.~268]{BGL} gives
$\operatorname{Ent}_\pi(f^2)\leq(2/K)\int\Gamma_D(f)\,d\pi$.
Thus $K=n/4$ gives \cref{eq:coordinate-lsi}. Gross's equivalence
\cite{Gross}, in \cite[Theorem~5.2.3, p.~246]{BGL}, gives
$\norm*{P_sz}_{L^{r_s}(\pi)}\leq\norm*{z}_{L^r(\pi)}$ with $r_s-1=(r-1)e^{2Ks}$.
Taking final exponent $2$ at time $s/2$ and $K=n/4$ gives initial
exponent $1+e^{-ns/4}$ and \cref{eq:coordinate-hypercontractivity}.
By entropy tensorization \cite[Proposition~5.2.7, p.~249]{BGL},
the constant $8/n$ is unchanged on the product, giving
\cref{eq:coordinate-product-lsi}.

\subsection{Orthogonality and local stability}\label{sec:modulation}

Following the variational approach of
Figalli and Zhang~\cite[proof of Lemma~4.1]{FZ},
we choose a nearby Sobolev extremal satisfying
the required orthogonality conditions.
Here we maximize a weighted pairing over the
family $w_\theta$ defined below.
The derivative estimates in \cref{lem:bubble}
allow us to control the choice of parameters,
uniformly in sufficiently large dimensions.

Set
\[
 \sigma_t^2=\frac{d_*}{n+1},\qquad K^2=\int k^2\,d\mu.
\]
Writing $b=(n-p)/(p-1)$, define
\[
 \gamma=(b\sigma_t)^{-1},\qquad \eta=\frac{\sqrt n}{bK},\qquad
 \lambda=e^{-\gamma\theta_0},\qquad
 w_\theta(x)=\lambda^{(n-p)/p}v(\lambda(x-\eta z)),
 \quad\theta=(\theta_0,z).
\]

\begin{lemma}\label{lem:iv-modulation}
Let $L,M,N$ be the constants in \cref{lem:bubble}, and set
\[
\begin{aligned}
 g&=(q-1)e^{-2},
 & R&=\min\{1,\gamma^{-1},g/(4NL)\},\\
 \delta_*&=\min\{1/4,g/(4N),gR/(4M)\},
 & C&=2+2LM/g.
\end{aligned}
\]
If $u\in\dot W^{1,p}$ satisfies
$\delta=\|\nabla(u-v)\|_p/E_v^{1/p}\leq\delta_*$,
then there exists a positive Sobolev extremal
$V=A w_{\theta_*}$ such that
$\int_{\R^n}V^{q-2}(u-V)\xi\,dx=0$ for every
$\xi\in T_V\mathcal M_{n,p}$. Moreover,
\[
\frac{\|\nabla(u-V)\|_p}{E(V)^{1/p}}
       \leq\frac{4C}{3}\delta.
\]
For each fixed $p>1$, $\delta_*$ is bounded below by
 a positive constant and $C$ is bounded above by
 a finite constant, independently of $n$ for all
 sufficiently large $n$.
\end{lemma}

\begin{proof}
Set $J(\theta;u)=V_q^{-1}\int u w_\theta^{q-1}\,dx$
and maximize it over $|\theta|\leq R$.
Since $R\leq\gamma^{-1}$, we have $|\log\lambda|\leq1$
on this ball, so \cref{lem:bubble} applies.
The identities $\int w_\theta^q=V_q$ and
$E(w_\theta)=E_v$ hold for every $\theta$.
When differentiating $J(\theta;u)$ with
respect to $\theta$, we keep $u$ fixed.
Differentiating the first identity twice
and using \eqref{eq:iv-chart-gram}, we obtain
\[
 D_\theta^2J(\theta;w_\theta)
 =-(q-1)\operatorname{diag}(1,\lambda^2I_n).
\]
By \eqref{eq:sobolev} and \eqref{eq:iv-chart-gradient}, we have
$\|u-w_\theta\|_q/V_q^{1/q}\leq\delta+L|\theta|$.
Consequently \eqref{eq:iv-chart-second} implies
\[
 \|D^2J(\theta;u)-D^2J(\theta;w_\theta)\|_{\rm op}
 \leq N(\delta+L|\theta|)\leq g/2.
\]

Since $\lambda^2\geq e^{-2}$ on the parameter ball,
$D^2J(\theta;u)\leq-gI/2$.
By \eqref{eq:iv-chart-first},
$|DJ(0;u)|\leq M\delta\leq gR/4$.
Together with the Hessian bound, this gives
$DJ(\theta;u)\cdot\theta
\leq-gR^2/4<0$ on $|\theta|=R$.
Thus $J$ has a unique maximum at an interior
point $\theta_*$, and the Hessian bound gives
$|\theta_*|\leq2M\delta/g$.

Set $A=J(\theta_*;u)$.
The maximum property and H\"older's inequality imply
$1-\delta\leq J(0;u)\leq A
\leq\|u\|_q/V_q^{1/q}\leq1+\delta$.
The choice of $A$ gives orthogonality in the
amplitude direction. Since $\lambda>0$,
the stationarity of $J$ and
\eqref{eq:iv-chart-gram} give orthogonality
in the translation and dilation directions.  Finally
\[
 \|\nabla(u-Aw_{\theta_*})\|_p/E_v^{1/p}
 \leq\delta+L|\theta_*|+|A-1|\leq C\delta.
\]
Dividing the residual bound by $A\geq3/4$ gives the stated relative
error.

The bounds $g\geq(p-1)e^{-2}$,
$\gamma\leq1$, and the uniform estimates
for $L,M,N$ in \cref{lem:bubble} give
the stated bounds for $R,\delta_*$ and $C$,
independently of $n$ for all sufficiently
large $n$.
\end{proof}

\begin{proposition}\label{prop:iv-invariant}
Let $\mathcal U$ be either $\dot W^{1,p}(\R^n)$ or
 the set of nonnegative functions in this space.
Suppose there exist $a,\rho>0$ and $r\geq1$ such that,
for every tangent-orthogonal $h$ satisfying
$v+h\in\mathcal U$ and
$\|\nabla h\|_p/E_v^{1/p}\leq\rho$, we have
\[
 F(v+h)/E_v\geq a\bigl(\|\nabla h\|_p/E_v^{1/p}\bigr)^r.
\]
Take the constants in Lemma~\ref{lem:iv-modulation}, put
$K_{\rm dist}=16C/3$, and choose
\[
 \varepsilon=\min\{1/8,\delta_*/4,\rho/K_{\rm dist},1/(2K_{\rm dist})\}.
\]
Then every nonzero $u\in\mathcal U$ satisfies
\[
 d(u)\leq\varepsilon\quad\Longrightarrow\quad
 \mathcal D(u)\geq a(2/3)^p2^{-r}d(u)^r.
\]
The proposition also applies if, under the same orthogonality
and admissibility assumptions, the local estimate is given by
$F(v+h)\geq c\|\nabla h\|_p^r$ for
$\|\nabla h\|_p\leq t_*$. In this case take
$a=cE_v^{r/p-1}$ and $\rho=t_*/E_v^{1/p}$.
\end{proposition}

\begin{proof}
Assume $d(u)\leq\varepsilon$.
For $d(u)>0$, choose an approximate extremal
$V_0$ with error at most
$2d(u)\|\nabla u\|_p$.

If $u\geq0$, we first show that $V_0$ is
positive. Indeed, for every positive
Sobolev extremal $W$, the Euler--Lagrange
equation gives

\[
 \int|\nabla W|^{p-2}\nabla W\cdot\nabla(u+W)
 =\frac{E(W)}{\int W^q}\int W^{q-1}(u+W)\geq E(W).
\]
H\"older's inequality and the triangle inequality therefore show that a negative
bubble $-W$ has error at least both $E(W)^{1/p}$ and
$E(u)^{1/p}-E(W)^{1/p}$, hence at least $E(u)^{1/p}/2$.
Since $d(u)\leq\varepsilon\leq1/8$, the chosen error is at most
$\|\nabla u\|_p/4$, so $V_0$ is positive.

Using the invariances of $\mathcal D$ and $d$,
we may now normalize $V_0$ to $v$.
The normalized initial error $\delta_0=\|\nabla(u-v)\|_p/E_v^{1/p}$ satisfies
\[
 \delta_0\leq\frac{2d(u)}{1-2d(u)}\leq4d(u).
\]

Lemma~\ref{lem:iv-modulation} gives
$\delta_V=\|\nabla(u-V)\|_p/E(V)^{1/p}\leq K_{\rm dist}d(u)$, because
its recovered amplitude is at least $3/4$.
The choices ensure $\delta_V\leq\min\{\rho,1/2\}$.
Using the translation, dilation, and amplitude
invariances, we may now assume that $V=v$.
We then apply the orthogonal local estimate.
Since
\[
 d(u)\leq\frac{\delta_V}{1-\delta_V}\leq2\delta_V,
 \qquad
 \frac{E(u)}{E(V)}\leq(1+\delta_V)^p\leq(3/2)^p,
\]
the stated lower bound follows.

If $d(u)=0$, there is a sequence of
Sobolev extremals converging to $u$ in
$\dot W^{1,p}(\R^n)$. Since $F$ is continuous
under this convergence and vanishes on
Sobolev extremals, we have $F(u)=0$.
\end{proof}

\section{Local stability in high dimensions}\label{sec:local}

In this section, we prove \cref{prop:local}, a local
stability estimate for nonnegative functions with a
radius independent of $n$ and a coefficient of order
$n^{-1}$. In Subsection~\ref{subsec:coupled-coercivity},
we estimate bounded perturbations by controlling the
radial and angular terms together.
Subsection~\ref{subsec:localized-estimates} gives
estimates for functions supported on small sets,
while Subsection~\ref{subsec:large-value-estimates}
treats large perturbations. Finally, in
Subsection~\ref{subsec:local-proof}, we combine
these estimates and apply \cref{prop:iv-invariant}.

\begin{proposition}\label{prop:local}\label{thm:effective-local}
For every $p>1$, there exist $N_p^L>p$, $\varepsilon_p>0$,
and $b_p>0$, depending only on $p$, such that for every
integer $n\geq N_p^L$ and every nonzero nonnegative
$u\in\dot W^{1,p}(\R^n)$,
\begin{equation}\label{eq:localinvariant}
 d(u)\leq\varepsilon_p
 \quad\Longrightarrow\quad
 \mathcal D(u)\geq\frac{b_p}{n}d(u)^{r_p}.
\end{equation}
\end{proposition}

To prove \cref{prop:local}, we take
$n\geq2p$ and increase the dimension
thresholds as needed.
By \eqref{eq:exponents}, this gives
$p<q\leq2p$ and
$p^2/n\leq q-p\leq2p^2/n$.

\subsection{Estimates for bounded perturbations}\label{subsec:coupled-coercivity}

Let $v,\mu,\nu$ and $\mathsf Z$ be as in
Subsection~\ref{sec:bubble-identities},
and let $\mathcal R$ and $D_{p,q}$ be
defined by \eqref{eq:weighted-remainder}
and \eqref{eq:local-potential}.
Let $A,B_S,T$ be the operators introduced in
Subsection~\ref{sec:basic-tools}.
The coordinate inverses act on the orthogonal complements of the
constants, and $T$ preserves means and contracts $L^\infty$.

As in Figalli and Zhang~\cite[Section~3]{FZ},
we use orthogonality to the tangent directions
of the manifold of Sobolev extremals.
We control the radial and angular terms
together using the auxiliary functions
$\psi_r$ and $\psi_S$ introduced below.

For a bounded admissible function $b$, set
$\bar b=\int b\,d\mu$. We choose $\psi_r$ and
$\psi_S$ so that
$A\psi_r+B_S\psi_S=b-\bar b$,
where the right-hand sides in the two equations
are bounded and have zero mean in the
corresponding variable.
Set
\[
d_n=\frac{p(1-t)}{n-p},\qquad C_n=\frac{p^2}{n(n-p)},\qquad
M_p=\operatorname{Id}+(p-2)e\otimes e,
\]
and, for $L>0$, define
\begin{equation}\label{eq:core-dual-field}
Y=-L(p-1)d_n(\psi_r)_te
  +LC_nn(1-t)\nabla_S\psi_S.
\end{equation}
The Jacobi-density identities
\[
(p-1)d_n^2d\nu_t=C_na\,d\mu_t,\qquad
\frac{p-1}{(n-p)t}\,d\nu_t=\frac{d\mu_t}{n(1-t)}
\]
give
\begin{equation}\label{eq:core-pairing-cost}
\begin{aligned}
\int Y\cdot\mathsf Zb\,d\nu
 &=LC_n\int\abs*{b-\bar b}^2\,d\mu,\\
\int Y^TM_p^{-1}Y\,d\nu&=L^2C_nI,\qquad
I=\int\left[
a\abs*{(\psi_r)_t}^2+\frac a{d_*}\abs*{\nabla_S\psi_S}^2\right]\,d\mu.
\end{aligned}
\end{equation}
The inverse-gradient bounds make $Y$ bounded, and its angular
component contains the factor $1-t$. The pairings are therefore
absolutely integrable when $\mathsf Zb\in L^p(\nu)$.
The slicing argument in \cref{lem:form-domain}
and the weak integration-by-parts identity
\eqref{eq:coordinate-weak-pairings} justify these
calculations without assuming that $b$ has
finite quadratic energy.

Write $\kappa=\min\{\tau,1-\tau\}$ and $k_0=k(\tau)$.
The inverse-gradient bounds in
\cref{eq:coordinate-radial-gradient,eq:coordinate-sphere-gradient}
and the construction in the proof of
\cref{lem:capacity} allow us to choose
\[
 R_p=\frac2\kappa+\frac\pi{\sqrt{\kappa d_*}},\qquad
 K_p^{\rm supp}=\max\{80+d_*^{-1},2+20d_*^{-1}\},\qquad
 Y_p=28p(p-1)R_p+8\sqrt2\pi p^2.
\]
If the radial and spherical right-hand sides are
bounded by $14\rho$ and $8\rho$, respectively,
the inverse-gradient estimates give
$\|Y\|_\infty\leq Y_pL\rho/(n\sqrt n)$.
For the two right-hand sides in the
\hyperref[sec:conditional-support]{construction}
in the proof of \cref{lem:capacity}, each bounded
by $2a_0$, the same estimates give
$\|Y\|_\infty\leq Y_pLa_0/(n\sqrt n)$.
These bounds hold for $n\geq\max\{2p,4,\kappa^{-1}\}$.

Set
\[
\sigma_t^2=\frac{d_*}{n+1},\qquad K^2=\int_0^1k(t)^2\,d\mu_t(t),
\qquad e_0=1,\qquad e_d=\frac{t-\tau}{\sigma_t},\qquad
e_z=\frac{\sqrt n}{K}k(t)(z\cdot\omega).
\]
Let $\mathbf e_1,\ldots,\mathbf e_n$ denote the Euclidean basis.
The functions $e_0,e_d,e_{\mathbf e_1},\ldots,e_{\mathbf e_n}$ are orthonormal in
$L^2(\mu)$. In the coercivity argument, unqualified $L^2$ norms and
inner products are taken with respect to $\mu$. Let $\mathcal T$ be the
span of these functions and write
\[
M_{\mathcal T}(g)=\norm*{\Pi_{\mathcal T}g}_2^2,
\qquad
\mathcal B_p(g)=\int(1+\abs*{\mathsf Zg})^{p-2}\abs*{\mathsf Zg}^2\,d\nu.
\]
For a function in $L^p(\mu)$, the same notation denotes the sum of
the squared tangent moments; their integrability follows from
\cref{app:bubble}.

For a convex function $H:\R^n\to\R$, write
\[
H^*(Y)=\sup_{z\in\R^n}\{Y\cdot z-H(z)\}.
\]
The Fenchel--Young inequality gives
$H(z)\geq Y\cdot z-H^*(Y)$.

\begin{lemma}\label{lem:iv-hessian}
Let $p>1$ and let $e\in\R^n$ be a unit vector. Set
$m=\min\{1,p-1\}$ and
\[
 H_p=|p-2|(4/3)^{|p-2|+1}\bigl(\max\{1,p-1\}+6\bigr),\qquad
 a_{\rm c}=\min\{1/4,m/[17(1+H_p)]\}.
\]
For $s\geq0$, let
$G_s(z)=|se+z|^p-s^p-ps^{p-1}e\cdot z$, and put
$M_p=I+(p-2)e\otimes e$.
If $|s-1|\leq a_{\rm c}/2$ and
$|Y|\leq4pm a_{\rm c}/17$ for $Y\in\R^n$, then
\[
 G_s^*(Y)\leq\frac{17}{32p}Y^\mathsf{T}M_p^{-1}Y.
\]
\end{lemma}

\begin{proof}
On $|x-e|\leq1/4$, the mean-value theorem gives
$||x|^{p-2}-1|\leq|p-2|(4/3)^{|p-2|+1}|x-e|$.
Moreover, for $|x-e|\leq1/4$,
\[
\left\|
\frac{x}{|x|}\otimes\frac{x}{|x|}
-e\otimes e
\right\|_{\rm op}
\leq6|x-e|.
\]
The Hessian formula for $|x|^p$ therefore gives
\[
 \|p^{-1}D^2(|x|^p)-M_p\|_{\rm op}\leq H_p|x-e|.
\]
For $|z|\leq a_{\rm c}/2$, $|se+z-e|\leq a_{\rm c}$, so
$D^2G_s(z)\geq(16p/17)M_p$.
For $|z|=a_{\rm c}/2$, integrating the Hessian
bound along the segment from $0$ to $z$ gives
\[
z\cdot D(G_s(z)-Y\cdot z)
\geq(16pm/17)|z|^2-|Y||z|>0.
\]
Since $G_s(z)-Y\cdot z$ is convex and grows
superlinearly, its global minimizer lies inside
this ball. At the minimizer, Taylor's formula gives
$G_s(z)\geq(8p/17)z^\mathsf{T}M_pz$.
Completing the square proves the stated bound
for $G_s^*(Y)$.
\end{proof}

\Needspace{8\baselineskip}
\begin{lemma}\label{lem:core}
For every $p>1$, there exist $\rho_p>0$, $\eta_p>0$,
and $N_p^{\rm bnd}>p$, depending only on $p$, such that
for every integer $n\geq N_p^{\rm bnd}$ and every
admissible function $g$ satisfying $|g|\leq\rho_p$
almost everywhere, with $\bar g=\int g\,d\mu$, we have
\begin{equation}\label{eq:defective-core}
\begin{gathered}
\mathcal R(1+g)-\int D_{p,q}(1+g)\,d\mu
 \geq a_p\epsilon\int\abs*{g-\bar g}^2\,d\mu
       +\eta_p\mathcal B_p(g)-C_p\epsilon M_{\mathcal T}(g),\\
a_p=p\left(\frac{40}{51}-\frac{17}{32}\right)=\frac{413p}{1632}>0.
\end{gathered}
\end{equation}
\end{lemma}

\begin{proof}
Take $\rho_p\leq1/2$ and
$N_p^{\rm bnd}\geq\max\{2p,4,2/\min\{\tau,1-\tau\}\}$,
with further increases specified below. Let $\xi=g-\bar g$ and set
\[
r(t)=\int_{\mathbb S^{n-1}}\xi(t,\omega)\,d\sigma(\omega),\qquad
s(\omega)=\int_0^1\xi(t,\omega)\,d\mu_t(t),\qquad w=\xi-r-s.
\]
The functions $r,s,w$ are pairwise orthogonal
in $L^2(\mu)$, and their $L^\infty$ norms are
bounded by $C\rho_p$. Moreover, $w$ has zero
mean in both variables. Set
\[
V_r=\|r\|_2^2,\qquad V_s=\|s\|_2^2,\qquad
V_w=\|w\|_2^2,
\]
and
\[
V=V_r+V_s+V_w
 =\int|g-\bar g|^2\,d\mu.
\]
Define
\begin{equation}\label{eq:core-potentials}
\psi_r=A^{-1}r+(n+A)^{-1}w,
\qquad \psi_S=B_S^{-1}(s+Tw).
\end{equation}
By \eqref{eq:core-potentials},
\[
A\psi_r=r+(\operatorname{Id}-T)w,\qquad B_S\psi_S=s+Tw,
\qquad A\psi_r+B_S\psi_S=\xi.
\]
Since $\|r\|_{L^\infty(\mu)},\|s\|_{L^\infty(\mu)}\leq2\rho_p$
and $\|w\|_{L^\infty(\mu)}\leq6\rho_p$, the $L^\infty$
contraction of $T$ bounds the radial and angular right-hand sides
by $14\rho_p$ and $8\rho_p$, respectively. Both have zero mean
in the variable being inverted. Thus \eqref{eq:core-pairing-cost}
applies with $b=g$ and variance $V$.
Let $P_1$ project onto the first spherical harmonics, and set
$U=\abs*{\langle g,e_d\rangle}^2$ and
$D_{\rm tr}=\sum_i\abs*{\langle g,e_{\mathbf e_i}\rangle}^2$.
The normalization of the translation directions gives
\[
\norm*{P_1\int_0^1 k(t)\xi(t,\cdot)\,d\mu_t(t)}_{L^2(\sigma)}^2
=n\sum_{i=1}^n\abs*{\int k\xi\omega_i\,d\mu}^2=K^2D_{\rm tr}.
\]
With $\bar k=\int_0^1k(t)\,d\mu_t(t)$,
\[
\bar k P_1s=P_1\int_0^1k(t)\xi(t,\cdot)\,d\mu_t(t)
 -\int_0^1(k(t)-\bar k)P_1w(t,\cdot)\,d\mu_t(t).
\]
Indeed,
\[
a\abs*{k'}^2=\frac{k^2(1/p-t)^2}{a}\leq\frac1a,\qquad
\int_0^1\frac1{a(t)}\,d\mu_t(t)
 =\frac{n-1}{n\tau-1}+\frac{n-1}{n/p-1}\leq\frac2{d_*}.
\]
Hence $\int_0^1\abs*{k-\bar k}^2\,d\mu_t\leq2/(nd_*)$ by the Jacobi gap.
The interval $[\tau/2,(1+\tau)/2]$ has mass at least $1/2$ for large
$n$, by the Beta variance $d_*/(n+1)$ and Chebyshev's inequality, and there
$k\geq k(\tau)/2$. Thus $\bar k\geq k(\tau)/4$ and, after increasing
the fixed threshold, $K^2/\bar k^2=1+\bar k^{-2}\int_0^1\abs*{k-\bar k}^2\,d\mu_t\leq2$.
Consequently,
\begin{equation}\label{eq:weighted-first-harmonic}
\norm*{P_1s}_2^2\leq4D_{\rm tr}+\frac{C-1}{n}V_w,
\qquad C=1+\frac{64}{d_*k_0^2}.
\end{equation}
The estimates leading to
\eqref{eq:weighted-first-harmonic} hold for
\[
n\geq\max\{2/\kappa^2,32/(d_*k_0^2),2/\kappa\}.
\]
A single Cauchy--Schwarz inequality in $L^2(\sigma)$ controls all
translation moments together, without a factor $n$.

The first two spectral levels in \cref{eq:coordinate-spectra} give
\begin{align*}
I_0&=\int[a\abs*{(\psi_r)_t}^2+\abs*{\nabla_S\psi_S}^2]\,d\mu
      =I_r+I_s+I_w,\\
I_r&=\langle r,A^{-1}r\rangle
 \leq\frac Un+\frac{V_r-U}{2(n+1)}
 \leq\frac{V_r}{2n}+\frac U{2n},\\
I_s&=\langle s,B_S^{-1}s\rangle
 \leq\frac{V_s}{2n}
 +\left(\frac1{n-1}-\frac1{2n}\right)\norm*{P_1s}_2^2,\\
I_w&=\norm*{A^{1/2}(n+A)^{-1}w}_2^2
          +\norm*{B_S^{-1/2}Tw}_2^2.
\end{align*}
For each product eigenfunction with eigenvalues
$\lambda_A\geq n$ and $\lambda_S\geq n-1$
for $A$ and $B_S$, respectively, we have
\begin{equation}\label{eq:mixed-cost-multiplier}
\begin{aligned}
\frac{\lambda_A}{(n+\lambda_A)^2}
 +\frac{n^2}{\lambda_S(n+\lambda_A)^2}
 &\leq\frac1{n+\lambda_A}
       +\frac{n}{(n-1)(n+\lambda_A)^2}\\
 &\leq\frac1{2n}+\frac1{4n(n-1)}.
\end{aligned}
\end{equation}
Combining \eqref{eq:mixed-cost-multiplier}
with \eqref{eq:weighted-first-harmonic},
we obtain
\[
I_0\leq\left(\frac1{2n}+\frac{C}{n^2}\right)V
       +\frac4n(U+D_{\rm tr}).
\]

To estimate the term with the angular weight,
set $G=\nabla_S\psi_S$. Using the spectral
decompositions of $A$ and $B_S$, we obtain
\[
\int a\abs*{\partial_tG}^2\,d\mu
 =\norm*{A^{1/2}B_S^{-1/2}Tw}_2^2
 \leq\frac{n}{4(n-1)}V_w\leq\frac12V.
\]
The inverse-gradient estimate makes $G$ bounded, so
\cref{eq:coordinate-stein} applies with this finite radial form energy.
Using $a=d_*+(1-2\tau)(t-\tau)-(t-\tau)^2$ and
$\sqrt a\leq1/2$ to cancel the factor $2$ in that identity gives
\[
\int\frac a{d_*}\abs*{G}^2\,d\mu
 \leq\norm*{G}_2^2+
 \frac{\abs*{1-2\tau}}{d_*n}\norm*{G}_2
 \left(\int a\abs*{\partial_tG}^2\,d\mu\right)^{1/2}.
\]
The mixed-mode bound $\|T\|_{2\to2}\leq1/2$ also gives
$\|G\|_2^2\leq(V+4U+4D_{\rm tr})/n$.
With $h_p=|1-2\tau|/d_*$, the preceding estimates imply
\[
 I\leq\frac1n\left(\frac12+\frac Cn+\frac{h_p}{\sqrt n}\right)V
       +\frac1n\left(4+\frac{2h_p}{\sqrt n}\right)(U+D_{\rm tr}).
\]
Thus $n\geq\max\{32C,(32h_p)^2\}$ gives
\begin{equation}\label{eq:strict-core-cost}
I\leq\frac9{16n}V+\frac9n(U+D_{\rm tr}).
\end{equation}
We now apply \cref{lem:iv-hessian}.
Choose
\[
 \rho_p\leq
 \min\{a_{\rm c}/2,1/4,\log(17/16)/(2p+1)\}.
\]
With $L=256pn/153$, the bound for $Y$
satisfies the hypothesis of that lemma provided
\[
 \sqrt n\geq
 \frac{256pY_p\rho_p}
 {153\,(4p\min\{1,p-1\}a_{\rm c}/17)}.
\]
Consequently, the dual field $Y$ satisfies
\begin{equation}\label{eq:core-small-conjugate}
 G_s^*(Y)\leq\frac{17}{32p}Y^TM_p^{-1}Y
 \quad(|s-1|\leq\rho_p).
\end{equation}
Combining the Fenchel--Young inequality
with \eqref{eq:core-small-conjugate},
\eqref{eq:core-pairing-cost}, and
\eqref{eq:strict-core-cost}, we obtain
\[
\mathcal R(1+g)
 \geq\frac{128p}{153}\epsilon V-C_p\epsilon(U+D_{\rm tr}).
\]
On the other hand, the full vector Bregman bound in
\cref{lem:eff-bregman} and $1/2\leq1+g\leq3/2$ give
$\mathcal R(1+g)\geq\widehat\gamma_p\mathcal B_p(g)$, where
$\widehat\gamma_p=\gamma_p\min\{2^{2-p},(3/2)^{p-2}\}>0$.
Combining these two lower bounds with weights $15/16$ and $1/16$,
\[
\mathcal R(1+g)
 \geq\frac{40p}{51}\epsilon V
   +\frac{\widehat\gamma_p}{16}\mathcal B_p(g)-C_p\epsilon(U+D_{\rm tr}).
\]
Finally, for $\abs*{g}\leq\rho_p$ and $p<q\leq2p$,
\[
D_{p,q}(1+g)
 \leq\frac p2\epsilon\frac{(1+\rho_p)^{2p-1}}{1-\rho_p}g^2
 \leq\frac{17p}{32}\epsilon g^2.
\]
Since $e_0=1$, the term $|\bar g|^2$ is included in $M_{\mathcal T}(g)$. Taking $\eta_p=\widehat\gamma_p/16$
proves \eqref{eq:defective-core}.
\end{proof}

\subsection{Small-support estimates}\label{subsec:localized-estimates}

\begin{lemma}\label{lem:capacity}
Fix $p>1$ and $L_0\geq2$. There exist
$\ell_p>0$, depending only on $p$, and
$N_p^{\rm supp}(L_0)>p$, depending only on $p,L_0$,
such that the following estimates hold for every integer
$n\geq N_p^{\rm supp}(L_0)$.

If $1<p<2$, let $a_0>0$ and let $\chi$ be an admissible
function satisfying $0\leq\chi\leq a_0$ and
$\mu(\chi\ne0)\leq\mathrm e^{-L_0}$.
Let $f\in W^{1,p}_{\rm loc}(\R^n)$ be nonnegative.
On the set $D_\chi=\{\mathsf Z\chi\ne0\}$, suppose that
$0\leq f\leq2a_0$ and $\mathsf Z\chi=\pm\mathsf Zf$
almost everywhere, with one fixed choice of sign. Then
\begin{equation}\label{eq:capacity-low}
\int_{D_\chi}G_f(\mathsf Zf)\,d\nu
 \geq\frac{\ell_pL_0}{n}a_0^{p-2}\int\chi^2\,d\mu.
\end{equation}

If $p\geq2$, let $X$ belong to the domain of the closed
quadratic form $\int|\mathsf ZX|^2\,d\nu$ on $L^2(\mu)$.
If $\mu(X\ne0)\leq\mathrm e^{-L_0}$, then
\begin{equation}\label{eq:capacity-high}
\int\abs*{\mathsf ZX}^2\,d\nu
 \geq\frac{\ell_pL_0}{n}\int X^2\,d\mu.
\end{equation}
\end{lemma}

\begin{proof}
\phantomsection\label{sec:conditional-support}
Take $N_p^{\rm supp}(L_0)
\geq\max\{4,\kappa^{-1}\}$.
Fix $a_0>0$ and consider a bounded
measurable function $\chi$ satisfying
$0\leq\chi\leq a_0$ and
$\mu(\chi\ne0)\leq e^{-L_0}$.
We construct $\psi_r$ and $\psi_S$
without using derivatives of $\chi$.

Put $m_0=e^{-L_0}$ and
\[
\Omega_-=\{\omega:\mu_t(\chi(\cdot,\omega)\ne0)>\sqrt{m_0}\}.
\]
Fubini's theorem gives $\sigma(\Omega_-)\leq\sqrt{m_0}$.
Set
\[
\chi_1=\chi1_{\Omega_-^c},\quad \chi_2=\chi1_{\Omega_-},\qquad
r_1(\omega)=\int_0^1\chi_1(t,\omega)\,d\mu_t(t),\quad
r_2(t)=\int_{\mathbb S^{n-1}}\chi_2(t,\omega)\,d\sigma(\omega).
\]
Thus
$\norm*{r_1}_{L^2(\sigma)}^2\leq\sqrt{m_0}\norm*{\chi_1}_{L^2(\mu)}^2$ and
$\norm*{r_2}_{L^2(\mu_t)}^2\leq\sqrt{m_0}\norm*{\chi_2}_{L^2(\mu)}^2$.

Define $\psi_r$ and $\psi_S$ by
\begin{equation}\label{eq:capacity-potentials}
A\psi_r=\chi_1-r_1+r_2-\int_0^1r_2(t)\,d\mu_t(t),
\qquad
B_S\psi_S=\chi_2-r_2+r_1-\int_{\mathbb S^{n-1}}r_1(\omega)\,d\sigma(\omega).
\end{equation}
Adding the two equations in \eqref{eq:capacity-potentials}
gives $A\psi_r+B_S\psi_S=\chi-\int\chi\,d\mu$.
Both right-hand sides in \eqref{eq:capacity-potentials}
are bounded by $2a_0$ and have zero mean in their
respective variables.
Write
\[
\begin{aligned}
\psi_r^{(1)}&=A^{-1}(\chi_1-r_1),&
\psi_r^{(2)}&=A^{-1}\left(r_2-\int_0^1r_2(t)\,d\mu_t(t)\right),\\
\psi_S^{(1)}&=B_S^{-1}(\chi_2-r_2),&
\psi_S^{(2)}&=B_S^{-1}\left(r_1-\int_{\mathbb S^{n-1}}r_1(\omega)\,d\sigma(\omega)\right).
\end{aligned}
\]
Both right-hand sides defining $\psi_r^{(1)}$
and $\psi_r^{(2)}$ have zero mean with respect
to $\mu_t$. Similarly, those defining
$\psi_S^{(1)}$ and $\psi_S^{(2)}$ have zero
mean with respect to $\sigma$. Define
\[
I_r^{(j)}=\int a\abs*{(\psi_r^{(j)})_t}^2\,d\mu,\qquad
I_S^{(j)}=\int\frac a{d_*}\abs*{\nabla_S\psi_S^{(j)}}^2\,d\mu.
\]

By construction,
$\mu_t(\chi_1(\cdot,\omega)\ne0)\leq\sqrt{m_0}$
for almost every $\omega$, and
$\sigma(\chi_2(t,\cdot)\ne0)\leq\sqrt{m_0}$
for almost every $t$.
We apply \eqref{eq:coordinate-support-inverse}
with $L=L_0/2$ to $\chi_1$ in the radial
variable and to $\chi_2$ in the angular variable.
For the remaining two terms, we use the
spectral gaps in \eqref{eq:coordinate-spectra}.
Since $a\leq1/4$, we obtain
\begin{align*}
I_r^{(1)}&\leq\frac{40}{nL_0}\norm*{\chi_1}_{L^2(\mu)}^2,&
I_r^{(2)}&\leq\frac{\sqrt{m_0}}n\norm*{\chi_2}_{L^2(\mu)}^2,\\
I_S^{(1)}&\leq\frac{10}{d_*nL_0}\norm*{\chi_2}_{L^2(\mu)}^2,&
I_S^{(2)}&\leq\frac{\sqrt{m_0}}{4d_*(n-1)}\norm*{\chi_1}_{L^2(\mu)}^2.
\end{align*}
Since $I\leq2\sum_{j=1}^2(I_r^{(j)}+I_S^{(j)})$ and
$\sqrt{m_0}\leq1/L_0$, these estimates give
\begin{equation}\label{eq:capacity-cost}
\begin{aligned}
I&=\int\left[a\abs*{(\psi_r)_t}^2
  +\frac a{d_*}\abs*{\nabla_S\psi_S}^2\right]\,d\mu\\
&\leq\frac{K_p^{\rm supp}}{nL_0}\int\chi^2\,d\mu,\\
K_p^{\rm supp}&=\max\{80+d_*^{-1},2+20d_*^{-1}\}.
\end{aligned}
\end{equation}
The indicators in \eqref{eq:capacity-potentials} appear only on
the right-hand sides, so they are never differentiated.

For $p<2$, apply this construction to the admissible function $\chi$
in the statement. We use the dual field $Y$ defined in
\eqref{eq:core-dual-field} with the functions $\psi_r,\psi_S$
constructed in \eqref{eq:capacity-potentials}.
Then \eqref{eq:core-pairing-cost} applies with $b=\chi$, and
$\norm*{Y}_\infty\leq Y_pLa_0/(n\sqrt n)$. Since $p<2$,
$M_p^{-1}\geq\operatorname{Id}$, so $\abs*{Y}^2\leq Y^TM_p^{-1}Y$.
Let $\kappa_p=p(p-1)3^{p-2}/2$. By \eqref{eq:full-bregman-bound}, for $0\leq f\leq2a_0$
we have
\[
G_f(z)\geq
\begin{cases}
\kappa_pa_0^{p-2}\abs*{z}^2,&\abs*{z}\leq a_0,\\
\kappa_p\abs*{z}^p,&\abs*{z}\geq a_0.
\end{cases}
\]
Therefore, if $\abs*{Y}\leq(\kappa_p/2)a_0^{p-1}$,
\[
G_f^*(Y)\leq\frac{a_0^{2-p}}{4\kappa_p}\abs*{Y}^2.
\]
For $\abs*{z}\geq a_0$, the quantity $Y\cdot z-G_f(z)$ is nonpositive.
Take
$L=\kappa_pnL_0a_0^{p-2}/K_p^{\rm supp}$.
The bound for $Y$ is then
$\norm*{Y}_\infty\leq Y_p\kappa_pL_0a_0^{p-1}/(K_p^{\rm supp}\sqrt n)$.
The bound for $G_f^*$ applies when
$n\geq(2Y_pL_0/K_p^{\rm supp})^2$; the height $a_0$ cancels.
If $\mathsf Z\chi=-\mathsf Zf$ on $D_\chi$,
we use $-Y$ instead of $Y$ in the
Fenchel--Young inequality.
Set $\bar\chi=\int\chi\,d\mu$. Then
\[
\int\abs*{\chi-\bar\chi}^2\,d\mu
 \geq(1-m_0)\int\chi^2\,d\mu
 \geq\frac12\int\chi^2\,d\mu.
\]
To estimate the conjugate term, we bound
the integral of $|Y|^2$ over $D_\chi$ by its
integral over the whole space. Using the energy
bound \eqref{eq:capacity-cost}, the Fenchel--Young
inequality on $D_\chi$ then gives
\[
\int_{D_\chi}G_f(\mathsf Zf)\,d\nu
 \geq LC_n\int|\chi-\bar\chi|^2\,d\mu
 -\frac{a_0^{2-p}L^2C_n}{4\kappa_p}I
 \geq\frac{\kappa_pp^2L_0}{4K_p^{\rm supp}n}
      a_0^{p-2}\int\chi^2\,d\mu.
\]
\Needspace{8\baselineskip}
For $p\geq2$, set $m_p=\min\{(p-1)^{-1},4d_*\}$.
The product logarithmic Sobolev inequality
\eqref{eq:coordinate-product-lsi}, the bound
$C_n\geq p^2/n^2$, and Jensen's inequality on the support give
\begin{align}
\int\abs*{\mathsf ZX}^2\,d\nu
 &=C_n\int\left[
 \frac a{p-1}\abs*{X_t}^2+\frac{d_*}{a}\abs*{\nabla_SX}^2\right]\,d\mu
 \label{eq:exact-Z-form}\\
 &\geq\frac{p^2m_p}{8n}\operatorname{Ent}_\mu(X^2)
 \geq\frac{p^2m_pL_0}{8n}\int X^2\,d\mu.\notag
\end{align}
Applying the estimate first to bounded
truncations of $X$ and then passing to the
limit in the closed quadratic form gives
the same inequality for every $X$ in the
stated form domain.
We may take
\[
 \ell_p=\frac{\kappa_pp^2}{4K_p^{\rm supp}}\quad(p<2),\qquad
 \ell_p=\frac{p^2}{8}\min\{(p-1)^{-1},4d_*\}\quad(p\geq2).\qedhere
\]
\end{proof}
\subsection{Estimates for large perturbations}\label{subsec:large-value-estimates}

Following the treatment of large perturbations in
Dolbeault--Esteban--Figalli--Frank--Loss
\cite[Sections~2.5 and 2.6.2]{DEFLL},
we split the part of the perturbation above a
large threshold into two terms. This allows us
to obtain a positive gradient contribution
while keeping the full nonlinear term
$D_{p,q}(1+\varphi)$.

\begin{lemma}\label{lem:high}
Fix $p>1$ and $0<\rho\leq\rho_p$.
There exist $M=M(p)\geq\mathrm e$, $\delta_T(p)>0$,
an integer $N_p^{\rm large}>p$, and positive constants
$c_{p,M},C_{p,M}$, all depending only on $p$,
with the following property.
For every integer $n\geq N_p^{\rm large}$,
let $u=v(1+\varphi)\in\dot W^{1,p}(\R^n)$ be nonnegative
and satisfy
$\delta=\|\nabla(u-v)\|_p/E_v^{1/p}\leq\delta_T(p)$.
Set
\[
H=\{\varphi>M\},\quad h_H=v(\varphi-M)_+,
\quad \chi=\min\{(\varphi-\rho)_+,M-\rho\},
\quad R_H=\int_HG_{1+\varphi}(\mathsf Z\varphi)\,d\nu.
\]
Then
\begin{equation}\label{eq:high-balance}
R_H-\int_HD_{p,q}(1+\varphi)\,d\mu
 \geq\frac{c_{p,M}}n
 \left(\frac{\norm*{\nabla h_H}_p}{E_v^{1/p}}\right)^{r_p}
 -\frac{C_{p,M}}n\int\chi^2\,d\mu.
\end{equation}
No tangent-orthogonality condition or upper bound on
$\varphi$ is assumed.
\end{lemma}

\begin{proof}
Take $N_p^{\rm large}\geq\max\{2p,N_p^{\rm supp}(2)\}$ and choose
\[
 C_{\rm pot}=p2^{2p}+p^{-1},\qquad
 M\geq\max\{e,(32pC_{\rm pot})^2\}.
\]
Let $c=M+1$, $\zeta=\frac14\min\{1,p-1\}$, and
$\vartheta=\min\{1,p/2\}$. Define
\[
\Phi_c(s)=(s^\zeta-c^\zeta)_+^{1/\zeta},\qquad
\Psi_c(s)=(s-c)_+-\Phi_c(s),\qquad
\xi=\Phi_c(1+\varphi),\qquad\eta=\Psi_c(1+\varphi).
\]
By \cref{lem:transfer}(ii), for $s>c$ and
$z\in\R^n$, we have
\begin{equation}\label{eq:high-coefficient-one}
G_s(z)\geq G_{\Phi_c(s)}(\Phi_c'(s)z)
       +\vartheta G_{c+\Psi_c(s)}(\Psi_c'(s)z).
\end{equation}
The functions $\Phi_c$ and $\Psi_c$ are $1$-Lipschitz,
vanish for $s\leq c$, and are bounded above by $(s-c)_+$.
By the chain rule, $\xi$ and $\eta$ are
admissible, and $h_H=v\xi+v\eta$.
Integrating \eqref{eq:high-coefficient-one},
we obtain
\begin{equation}\label{eq:high-residual-split}
R_H\geq\mathcal R(\xi)+\vartheta\mathcal R(c+\eta).
\end{equation}
By \eqref{eq:multiplier-bounds} and the truncation chain rule,
$\|\xi\|_{L^q(\mu)},\|\eta\|_{L^q(\mu)}\leq\delta$,
$\|\xi\|_{L^p(\nu)},\|\eta\|_{L^p(\nu)}\leq p'\delta$,
and $\|\mathsf Z\xi\|_{L^p(\nu)},\|\mathsf Z\eta\|_{L^p(\nu)}
\leq(1+p')\delta$.

Put $a^\sharp=p/q+\epsilon/(8p)$ and
$b^\sharp=p/q+\epsilon/(4p)$.
By \cref{lem:transfer}(iii), we have
\begin{equation}\label{eq:high-scalar-transfer}
a^\sharp r^q-r^p
 \leq b^\sharp\Phi_c(1+r)^q-\Phi_c(1+r)^p
       +C_{p,M}\epsilon1_{\{M<r<R_0\}},\qquad r>M,
\end{equation}
where $R_0$ is fixed by $p,M$ before $n$ varies.
Set $E_\xi=E(v\xi)/E_v$.
Integrating \eqref{eq:high-scalar-transfer}
and using \eqref{eq:local-picone}
and \eqref{eq:high-residual-split}, we obtain
\begin{align}\label{eq:high-paired-budget}
R_H+\int_H\varphi^p\,d\mu-a^\sharp\int_H\varphi^q\,d\mu
 &\geq E_\xi-b^\sharp\int\xi^q\,d\mu
    +\vartheta\mathcal R(c+\eta)-C_{p,M}\epsilon\mu(H).
\end{align}
By \eqref{eq:sobolev},
\[
\left(\int\xi^q\,d\mu\right)^{p/q}\leq E_\xi.
\]
Since $q>p$ and
$\|\xi\|_{L^q(\mu)}\leq\delta\leq1$,
we have $\delta^{q-p}\leq1$, and hence
\begin{equation}\label{eq:high-sobolev-reserve}
E_\xi-b^\sharp\int\xi^q\,d\mu
 \geq(1-b^\sharp)E_\xi,
\qquad
1-b^\sharp=\epsilon\left(\frac1q-\frac1{4p}\right)
 \geq\frac p{4n}\quad(n\geq2p).
\end{equation}

Choose $0<\delta_T\leq[2(1+2p')]^{-1}$. Then $\eta$ is supported in $H$,
$\mu(H)\leq(\delta/M)^q<e^{-2}$, and its $L^p(\mu)$ norm and the
$L^p(\nu)$ norms of $\eta,\mathsf Z\eta$ are at most $c=M+1$.
By \cref{lem:shifted} and the bound \eqref{eq:shifted-coefficient},
\[
\mathcal R(c+\eta)\geq\frac{k_{\rm sh}}n
 \left(\frac{\norm*{\nabla(v\eta)}_p}{E_v^{1/p}}\right)^{r_p}.
\]
By \eqref{eq:truncated-gradient-bound}
and the choice of $\delta_T$, we have
$E_\xi^{1/p}\leq1/2$.
Since $r_p\geq p$, it follows that
$E_\xi\geq E_\xi^{r_p/p}$.
With the value $k_{\rm sh}>0$ in \cref{eq:shifted-coefficient}, the
resulting gradient coefficient may be taken as
\[
 c_{p,M}=2^{1-r_p}\min\{p/4,\min(1,p/2)k_{\rm sh}\}>0.
\]
Combining \eqref{eq:high-paired-budget}
with \eqref{eq:high-sobolev-reserve}
and the preceding gradient estimates,
and using $h_H=v\xi+v\eta$, we obtain
\begin{equation}\label{eq:high-before-potential}
R_H+\int_H\varphi^p\,d\mu
 \geq a^\sharp\int_H\varphi^q\,d\mu
 +\frac{c_{p,M}}n\left(\frac{\norm*{\nabla h_H}_p}{E_v^{1/p}}\right)^{r_p}
 -C_{p,M}\epsilon\mu(H).
\end{equation}
Since $\chi=M-\rho$ on $H$, $\rho\leq1/2$, and
$\epsilon\leq2p^2/n$, we have
\begin{equation}\label{eq:high-measure-error}
\epsilon\mu(H)
\leq\frac{2p^2}{n(M-1/2)^2}
\int\chi^2\,d\mu.
\end{equation}

It remains to compare the scalar term with
$D_{p,q}(1+r)$. Integrating with respect
to the exponent, we obtain
\begin{align*}
D_{p,q}(1+r)-\left(\frac pq r^q-r^p\right)
 &=p\int_p^q\left[\int_r^{r+1}x^{s-1}\log x\,dx+\frac1{s^2}\right]ds\\
 &\leq C_{\rm pot}\epsilon\frac{\log r}{r}r^q
 \leq\frac{\epsilon}{16p}r^q\qquad(r\geq M).
\end{align*}
The last inequality uses $\log r/r\leq\log M/M\leq2/\sqrt M$
and the chosen bound on $M$. Thus
\[
a^\sharp r^q-r^p-D_{p,q}(1+r)
\geq\frac{\epsilon}{16p}r^q\geq0\qquad(r\geq M).
\]
Combining this inequality with
\eqref{eq:high-before-potential}
and \eqref{eq:high-measure-error}
proves \eqref{eq:high-balance}.
\end{proof}
\subsection{Proof of the local estimate}\label{subsec:local-proof}

\begin{proof}
Following the decomposition used in
\cite[Section~2.5]{DEFLL}, we first prove the
estimate for $h=v\varphi$ satisfying
$v+h\geq0$ and the orthogonality conditions
\begin{equation}\label{eq:weighted-orthogonality}
\int\varphi\,d\mu=\int t\varphi\,d\mu=0,
\qquad \int k\omega_i\varphi\,d\mu=0\quad(1\leq i\leq n).
\end{equation}
Choose $\rho\in(0,\rho_p]$ small enough for
\cref{lem:core}, then choose $M$ as in
\cref{lem:high}. Let
\[
g=\min\{\rho,\max\{-\rho,\varphi\}\},\quad
B=\{\abs*{\varphi}>\rho\},\quad H=\{\varphi>M\},\quad I=B\setminus H,
\]
and set
\[
\chi=\min\{(\varphi-\rho)_+,M-\rho\},\quad
\eta_-=(1-\rho-(1+\varphi))_+,\quad w_I=\chi-\eta_-,
\quad h_H=v(\varphi-M)_+.
\]
This gives the exact decomposition
\begin{equation}\label{eq:local-three-parts}
h=vg+vw_I+h_H,
\qquad
\mathcal R(1+\varphi)=\mathcal R(1+g)+R_I+R_H,
\end{equation}
where
$R_I=\int_I G_{1+\varphi}(\mathsf Z\varphi)\,d\nu$.
Put
\[
S_I=\int(\chi^2+\eta_-^2)\,d\mu
=\int w_I^2\,d\mu.
\]
The integral includes the contribution from $H$,
where $\chi=M-\rho$.

Put $b=\varphi-g$ and suppose $\delta\leq\rho$.
By \eqref{eq:multiplier-bounds},
$\mu(B)\leq(\delta/\rho)^p$ and
$\|b\|_{L^p(\mu)}\leq\delta$. For every unit tangent
combination $e_c$, H\"older's inequality and \cref{lem:bubble} give
\[
\abs*{\int be_c\,d\mu}
 \leq M_T\delta\,\mu(B)^{3/(4p')}
 \leq M_T\rho^{-3(p-1)/4}\delta^{(3p+1)/4}.
\]
Taking the supremum over unit coefficient
vectors and using the orthogonality conditions
\eqref{eq:weighted-orthogonality}, we obtain
\begin{equation}\label{eq:clipped-tangent-defect}
M_{\mathcal T}(g)\leq T_p\delta^{t_p},\qquad
 T_p=M_T^2\rho^{-3(p-1)/2},
\qquad t_p=\frac{3p+1}{2}>r_p.
\end{equation}
This estimate also includes the coefficient corresponding to $e_0=1$.

We first estimate $D_{p,q}(1+\varphi)$ on $I$. On $I$, $g=\pm\rho$
and $\varphi-g=w_I$. For $0\leq f\leq1$,
\[
f^{p-1}\abs*{\log f}
 \leq\left(1+\frac1{e(p-1)}\right)(1-f).
\]
The left-hand side is defined by continuity at $f=0$.
For $1\leq f\leq M+1$ use $\log f\leq f-1$.
Integrating $D'_{p,q}(f)=pf^{p-1}(f^\epsilon-1)$ along the
interval between $1+g$ and $1+\varphi$ therefore gives
\[
0\leq D_{p,q}(1+\varphi)-D_{p,q}(1+g)
 \leq C_{p,M}\epsilon\left(\rho\abs*{w_I}+\frac12w_I^2\right)
 \quad\hbox{on }I.
\]
Writing $m_g=\int g\,d\mu$ and $V=\int\abs*{g-m_g}^2\,d\mu$,
we have $\rho^2\mu(I)\leq V+m_g^2$. Young's inequality gives
\begin{equation}\label{eq:middle-potential}
\int_I[D_{p,q}(1+\varphi)-D_{p,q}(1+g)]\,d\mu
 \leq\frac{a_p}{4}\epsilon V+C_{p,M}\epsilon S_I
                  +C_p\epsilon\delta^{t_p}.
\end{equation}

We next estimate $R_I$ using
\cref{lem:capacity}. Fix $L_0\geq2$.
We will choose $L_0$ and then restrict
$\delta$ so that $\mu(B)\leq e^{-L_0}$. If $p<2$, apply
\cref{lem:capacity} to $\chi$ with $a_0=M+1$ and to
$\eta_-$ with $a_0=1-\rho$, using the negative sign for the
latter. Their gradient supports are disjoint and contained in $I$.
If $p\geq2$, use instead the capped powers
\[
X_+=[\min(1+\varphi,M+1)^{p/2}-(1+\rho)^{p/2}]_+,
\qquad X_-=[(1-\rho)^{p/2}-(1+\varphi)^{p/2}]_+.
\]
Their gradients are supported only in $I$. If $R_I<\infty$,
\eqref{eq:full-bregman-bound} and the chain rule show that
$X_+$ and $X_-$ have finite quadratic energy. Their membership
in the closed form domain follows from \cref{lem:form-domain};
for $R_I=\infty$, the desired bound is immediate. Moreover,
\[
R_I\geq\frac{4\gamma_p}{p^2}\int
(\abs*{\mathsf ZX_+}^2+\abs*{\mathsf ZX_-}^2)\,d\nu.
\]
Since
$X_+^2\geq\chi^2$ and
$X_-^2\geq(1-\rho)^{p-2}\eta_-^2$,
both cases give
\begin{equation}\label{eq:middle-capacity}
R_I\geq\frac{\Lambda_{p,M}L_0}{n}S_I,
\qquad \Lambda_{p,M}>0.
\end{equation}
Using $\ell_p$ from \cref{lem:capacity} and $\gamma_p$ from
\cref{lem:eff-bregman}, a sufficient choice is
\[
\Lambda_{p,M}=\begin{cases}
\ell_p(M+1)^{p-2},&1<p<2,\\
(4\gamma_p\ell_p/p^2)2^{2-p},&p\geq2.
\end{cases}
\]
Both gradients vanish on $H$.

Let $C_1,C_2$, depending on $p,M$, be the coefficients of
$\epsilon S_I$ in \cref{eq:middle-potential} and
$n^{-1}\int\chi^2\,d\mu$ in \cref{eq:high-balance}, respectively. Choose
\[
L_0\geq\max\left\{2,\frac{4(2p^2C_1+C_2)}{\Lambda_{p,M}}\right\}.
\]
Since $\epsilon\leq2p^2/n$ and $\int\chi^2\,d\mu\leq S_I$,
\cref{eq:middle-capacity} gives
\[
C_1\epsilon S_I+\frac{C_2}{n}\int\chi^2\,d\mu
 \leq\frac14R_I.
\]
Thus $L_0$ is independent of $n$.
We now combine the estimates for the three parts
in \eqref{eq:local-three-parts}.
By \eqref{eq:local-representation}, we subtract
the $D_{p,q}(1+\varphi)$ term from the Bregman remainder.
On $H$, we may replace
$D_{p,q}(1+\varphi)-D_{p,q}(1+g)$ by
$D_{p,q}(1+\varphi)$, since
$D_{p,q}(1+g)\geq0$.
Combining \cref{lem:core,lem:high} with
\eqref{eq:middle-potential}, we obtain
\begin{equation}\label{eq:local-ledger}
\frac{F(v+h)}{E_v}
 \geq\frac{3a_p}{4}\epsilon V+\eta_p\mathcal B_p(g)
       +\frac34R_I+\frac{c_{p,M}}n\delta_H^{r_p}
       -\frac{C_p}{n}\delta^{t_p},
\qquad \delta_H=\frac{\norm*{\nabla h_H}_p}{E_v^{1/p}}.
\end{equation}

We now bound the gradients of the three terms in
\eqref{eq:local-three-parts} by the positive terms in
\eqref{eq:local-ledger}.
Let $\delta_g=\norm*{\nabla(vg)}_p/E_v^{1/p}$.
By \eqref{eq:bubble-quadratic-identity},
\[
V+m_g^2=\int g^2\,d\nu+2\int g\,e\cdot\mathsf Zg\,d\nu.
\]
For $p\geq2$, Young's inequality with $\abs*{g}\leq1$ and
$(1+\abs*{z})^{p-2}\abs*{z}^2\geq\max\{\abs*{z}^2,\abs*{z}^p\}$
gives \cref{eq:good-gradient-conversion} directly.
If $p<2$, apply Young's inequality on $\abs*{\mathsf Zg}\leq1$ and use
$(1+\abs*{z})^{p-2}\abs*{z}^2\geq2^{p-2}\abs*{z}$ on its complement, together
with $\abs*{g}\leq1$. This gives
$\int g^2\,d\nu\leq2(V+m_g^2)+C_p\mathcal B_p(g)$.
Also
\[
\norm*{\mathsf Zg}_{L^p(\nu)}^2
 \leq\mathcal B_p(g)(1+\norm*{\mathsf Zg}_{L^p(\nu)})^{2-p}.
\]
The chain rule and \eqref{eq:multiplier-bounds} give
$\norm*{\mathsf Zg}_{L^p(\nu)}\leq\norm*{\mathsf Z\varphi}_{L^p(\nu)}
\leq(1+p')\delta$.
Since $\nu(\R^n)=1$ and $p<2$, one has
$\norm*{g}_{L^p(\nu)}\leq\norm*{g}_{L^2(\nu)}$.
Thus, in both cases,
\begin{equation}\label{eq:good-gradient-conversion}
\delta_g^{r_p}\leq C_p[V+m_g^2+\mathcal B_p(g)].
\end{equation}

For the middle function let
$X=\norm*{w_I}_{L^p(\nu)}$, $U=\norm*{w_I}_{L^p(\mu)}$ and
$Y=\norm*{\mathsf Zw_I}_{L^p(\nu)}$.
Set $J_p=p/(1-2^{-p})$.
Identity \cref{eq:local-picone-signed} gives
$X^p\leq U^p+pX^{p-1}Y$.
For $p\geq2$, $U^p\leq C_{p,M}S_I$ and
$Y^p\leq C_pR_I$. Young's inequality and
\cref{eq:middle-capacity} give $\delta_I^p\leq C_{p,M}nR_I$.
For $p<2$, $U^2\leq S_I$ and \eqref{eq:full-bregman-bound} gives
\[
Y^2\leq C_pR_I(\norm*{1+\varphi}_{L^p(\nu)}+Y)^{2-p}
     \leq C_pR_I.
\]
Indeed $\|1+\varphi\|_{L^p(\nu)}+Y\leq2+2p'$
by \eqref{eq:multiplier-bounds} and $\delta\leq1$.
If $X<2U$, we use $U^2\leq S_I$.
If $X\geq2U$, the inequality
$X^p\leq U^p+pX^{p-1}Y$ gives
$X\leq J_pY$.
Combining these bounds with
\eqref{eq:middle-capacity} and the estimate
for $Y$ gives the same conclusion with
exponent $2$.
Thus
\begin{equation}\label{eq:middle-gradient-conversion}
\delta_I^{r_p}\leq C_{p,M}nR_I,
\qquad \delta_I=\frac{\norm*{\nabla(vw_I)}_p}{E_v^{1/p}}.
\end{equation}

We can choose explicit constants in
\eqref{eq:good-gradient-conversion} and
\eqref{eq:middle-gradient-conversion}.
With $J_p$ as above, one may take
\[
G_p=\begin{cases}
7\cdot2^{p-1},&p\geq2,\\
4+2^{6-p}+2(2+p')^{2-p},&1<p<2,
\end{cases}
\]
and
\[
H_{p,M}=\begin{cases}
2^{p-1}\bigl(2^{p-1}M^{p-2}/\Lambda_{p,M}
 +(J_p+1)^p/\gamma_p\bigr),&p\geq2,\\
4/\Lambda_{p,M}+2(J_p+1)^2\gamma_p^{-1}
 (2+2p')^{2-p},&1<p<2.
\end{cases}
\]
Indeed, \eqref{eq:bubble-quadratic-identity}
and Young's inequality control
$\int g^2\,d\nu$ by the variance and
$\mathcal B_p(g)$. For $p\geq2$, use $|g|\leq1$ and
$\int|\mathsf Zg|^2\,d\nu\leq\mathcal B_p(g)$.
For $p<2$, split at $|\mathsf Zg|=1$ and use
$\|\mathsf Zg\|_{L^p(\nu)}\leq(1+p')\delta$.
For the middle term, \eqref{eq:local-picone-signed} gives
$X\leq2U+J_pY$, while \eqref{eq:middle-capacity}
controls $U$ and the Bregman bound controls $Y$.
These choices use $L_0\geq2$ and $\delta\leq1$.

We now choose the constants in the final estimate.
Set
\[
 k_p=\min\{3a_pp^2/4,\eta_p\},\qquad
 A_p=3^{1-r_p}\min\{k_p/G_p,3/(4H_{p,M}),c_{p,M}\},\qquad
 B_p=1+C_p+k_pT_p.
\]
Here $C_p$ is the error coefficient in \cref{eq:local-ledger}.
Using $p^2/n\leq\epsilon\leq2p^2/n$,
\cref{eq:clipped-tangent-defect}, and
$\delta\leq\delta_g+\delta_I+\delta_H$ gives
\[
 \frac{F(v+h)}{E_v}\geq\frac{A_p}{n}\delta^{r_p}
                             -\frac{B_p}{n}\delta^{t_p},
 \qquad A_p>0,\quad B_p\geq1.
\]
Having fixed $\rho$, $M$, and $L_0$ above, take
$N_p^{\rm orth}\geq
\max\{N_p^{\rm bnd},N_p^{\rm supp}(L_0),N_p^{\rm large},N_p\}$,
where $N_p$ is the threshold in \cref{lem:bubble}, and choose
\[
 0<\delta_p^{\rm orth}\leq
 \min\left\{\delta_T,\frac\rho2,\rho e^{-L_0/p},
  \left(\frac{A_p}{2B_p}\right)^{1/(t_p-r_p)}\right\}.
\]
These choices ensure that the support
conditions in \cref{lem:capacity,lem:high}
and the form-domain conditions used above
are satisfied. Since
$t_p-r_p>0$, we obtain
\begin{equation}\label{eq:orthogonal-local}
\frac{F(v+h)}{E_v}\geq\frac{A_p}{2n}
 \left(\frac{\norm*{\nabla h}_p}{E_v^{1/p}}\right)^{r_p}
\end{equation}
under \cref{eq:weighted-orthogonality},
$\delta\leq\delta_p^{\rm orth}$ and $n\geq N_p^{\rm orth}$.
All parameters, including the radius, are fixed before $n$ varies.

We apply \cref{prop:iv-invariant} to
\eqref{eq:orthogonal-local}, with
$a=A_p/(2n)$, $r=r_p$, and
$\rho=\delta_p^{\rm orth}$.
By \cref{lem:iv-modulation,lem:bubble},
the constants needed in this application
are uniform in all sufficiently large
dimensions, with a threshold depending
only on $p$.
Thus $\varepsilon_p$ can be chosen
independently of $n$. Taking
\[
 b_p=\frac{A_p}{2}(2/3)^p2^{-r_p}>0
\]
gives \cref{eq:localinvariant}. Finally, choose $N_p^L$ large enough to
satisfy all the dimension restrictions
used above.
\end{proof}

\section{Fixed-dimensional local stability}\label{sec:fixed-local}

In this section, we prove local stability estimates
for perturbations orthogonal to the tangent space
in each fixed dimension $n>p$.
In Subsection~\ref{sec:exact-spectrum}, we establish
the spectral gap \eqref{eq:stage4-weaker-gap},
the weak identities \eqref{eq:stage4-weak-pairings},
and the finite spectral estimate
\eqref{eq:stage4-spectral-tail}.
In Subsection~\ref{sec:superquadratic-effective},
we prove \cref{thm:super-local} for $p\geq2$,
while Subsection~\ref{sec:subquadratic-effective}
establishes \cref{thm:short-subquadratic-local}
for $1<p<2$.

\subsection{Spectral estimates}\label{sec:exact-spectrum}

Fix an integer $n>p>1$, and let $v$ be the
Sobolev extremal defined in
\eqref{eq:standard-bubble}.
Write $q'=q/(q-1)$, $\mathcal V=V_q$, $\mathcal G=E_v^{1/p}$, and
\[
 \Lambda=E_v/V_q=n\left(\frac{n-p}{p-1}\right)^{p-1},
 \qquad W=v^{q-2}.
\]
We write $L^2(W)=L^2(\mathbb R^n,W(x)\,dx)$ and set
\[
 \mathbb A=|\nabla v|^{p-2}(I+(p-2)\omega\otimes\omega),\qquad
 Q(h)=\int\mathbb A\nabla h\cdot\nabla h.
\]
Let $D(Q)$ be the closure of
\[
 \mathcal D_0=\{h\in C_c^\infty(\mathbb R^n):
                    h\text{ is constant near }0\}
\]
with respect to the norm
$(Q(h)+\|h\|_{L^2(W)}^2)^{1/2}$.
We apply Zhang's spectral theorem
\cite[Theorem~1.1]{ZhangSpectrum} to the
resulting closed form. In particular, the
form domain is determined by the same test-function class
$\mathcal D_0$, including the endpoint conditions selected by
the closure when $n(p-1)/p\leq1$.

Set
\[
 \alpha=\frac{n(p-1)}p,\quad \beta=\frac np,\quad
 d=\frac{p-1}{p^2},\quad C_n=\frac{p^2}{n(n-p)},\quad
 \mu_{\mathrm{tan}}=(q-1)\Lambda,
\]
and let $\mathcal T$ be the span of $v$ and its infinitesimal dilation
and translation derivatives.

\begin{proposition}\label{prop:stage4-spectral-toolkit}
The self-adjoint operator associated with $Q$ in $L^2(W)$ has
an orthonormal basis $(e_i)$ of eigenfunctions, indexed by
angular and radial degrees $(\ell,j)$ and the multiplicity of
spherical harmonics. The function $v$ has eigenvalue
$(p-1)\Lambda$, while its infinitesimal dilation and
translation derivatives have eigenvalue $\mu_{\mathrm{tan}}$.
For every $h\in D(Q)$ orthogonal to $\mathcal T$ in $L^2(W)$,
\begin{equation}\label{eq:stage4-weaker-gap}
 Q(h)\ge\mu_\dagger\|h\|_{L^2(W)}^2,\qquad
 \mu_\dagger:=\mu_{\mathrm{tan}}+
 \frac{4\Lambda(p-1)(n+1)}{n(n-p)}>\mu_{\mathrm{tan}}.
\end{equation}
For each $i$, write $e_i=vb_i$ and set $T_i=\mathbb A\nabla e_i$.
Then $b_i$ is bounded, and the norms $\|b_i\|_\infty$,
$\|T_i\|_{p'}$, and $\|We_i\|_{q'}$ admit upper bounds
obtained by finite algebraic operations, powers, and Beta integrals.
Moreover,
\begin{equation}\label{eq:stage4-weak-pairings}
 \int T_i^T\mathbb A^{-1}T_k=\mu_i\delta_{ik},\qquad
 \int T_i\cdot\nabla h=\mu_i\int We_i h
 \quad(h\in\dot W^{1,p}).
\end{equation}
For $L>0$, choose an integer $M$ satisfying
\begin{equation}\label{eq:stage4-short-cutoff}
 M>\sqrt{\frac{L}{4\Lambda C_nd}}.
\end{equation}
Let $P_M$ be the orthogonal projection onto the span of the
eigenfunctions indexed by $0\le\ell,j\le M$, including all
angular multiplicities.
Every omitted eigenvalue exceeds $L$, and hence
\begin{equation}\label{eq:stage4-spectral-tail}
 \|(I-P_M)h\|_{L^2(W)}^2\le L^{-1}Q(h)
\end{equation}
for every $h\in D(Q)$.
\end{proposition}

\begin{proof}
Write $r=|x|$ and use $t=r^{p'}/(1+r^{p'})$ from
Subsection~\ref{sec:bubble-identities}. Set $s_0=r_0=0$; for $\ell\ge1$, let
$s_\ell,r_\ell$ be the positive roots of
\[
 s_\ell(s_\ell+\alpha-1)
 =r_\ell(r_\ell+\beta-1)
 =d\ell(\ell+n-2),\qquad z_\ell=s_\ell+r_\ell.
\]
In our notation, Zhang's eigenfunctions and
eigenvalues take the form
\begin{align}
 e_{\ell,j,k}
 &\propto v\,t^{s_\ell}(1-t)^{r_\ell}
 P_j^{(\alpha+2s_\ell-1,\,\beta+2r_\ell-1)}(1-2t)Y_{\ell,k},
 \label{eq:stage4-finite-eigenfunctions}\\
 \mu_{\ell,j}
 &=\Lambda\{p-1+C_n(j+z_\ell)(j+z_\ell+n-1)\}.
 \label{eq:stage4-imported-eigenvalues}
\end{align}
Indeed, Zhang's parameters satisfy
$s_{\rm Z}=p's_\ell$, $\tau_{\rm Z}=(n-p)/(p-1)+p'r_\ell$,
$m_{\rm Z}=(n-p)/p+z_\ell$, and $C_{\rm Z}=\Lambda C_n$.
We choose the eigenfunctions in
\eqref{eq:stage4-finite-eigenfunctions}
to be normalized in $L^2(W)$.

To prove \eqref{eq:stage4-weaker-gap}, write $h=vb$.
Recall from \eqref{eq:local-measures} and
\eqref{eq:beta-measures} that
$d\mu=V_q^{-1}v^q\,dx$ is the product of
the $\operatorname{Beta}(\alpha,\beta)$ measure
and the normalized measure on the sphere.
Using \eqref{eq:bubble-euler-lagrange},
we obtain
\[
 \frac{Q(vb)}{\Lambda\int v^q}
 =(p-1)\int b^2\,d\mu+C_n\int\left\{
 t(1-t)b_t^2+\frac d{t(1-t)}|\nabla_Sb|^2\right\}d\mu.
\]
This identity holds first on $\mathcal D_0$ and hence on its closure.
Let $H_\ell$ be the operator corresponding to angular degree
$\ell$ in the last form. The spectrum for $\ell=0$ is $j(j+n-1)$;
since $s_1=1/p$ and $r_1=(p-1)/p$, the spectrum for $\ell=1$ is
$(j+1)(j+n)$. After excluding the eigenfunctions in $\mathcal T$, both are bounded
below by $2(n+1)$. 

For $\ell\ge2$, the lowest eigenvalue of
$H_1$ is $n$. Since $t(1-t)\le1/4$,
comparison on the common dense class of test functions for $H_\ell$ and $H_1$ gives
\[
 H_\ell=H_1+
 \frac{d\{\ell(\ell+n-2)-(n-1)\}}{t(1-t)}
 \ge n+4d(n+1).
\]
The comparison extends to the closed forms.  Since $4d\le1$,
$2(n+1)>n+4d(n+1)$.  Summing over the spherical harmonic degrees proves
\eqref{eq:stage4-weaker-gap}.

The coefficients and normalization factors in
\eqref{eq:stage4-finite-eigenfunctions}
can be computed from Beta and spherical moments.
Spherical harmonics can be constructed by finite
linear algebra on homogeneous harmonic polynomials,
whose coefficients give bounds for the harmonics
and their first derivatives on the sphere.
Since $s_\ell,r_\ell\ge0$, this also gives
a bound for $b_i$.
Differentiating \eqref{eq:stage4-finite-eigenfunctions}, we obtain
bounds with these coefficient sums as constants:
\[
 |T_i(x)|\le K_i
 \begin{cases}
 |x|,& |x|\le1,\ \ell=0,\\
 |x|^{p's_\ell-1/(p-1)},& |x|\le1,\ \ell\ge1,\\
 |x|^{-(n-1)-p'r_\ell},& |x|\ge1.
 \end{cases}
\]
The exponent at zero is nonnegative because $s_\ell\ge s_1=1/p$;
at infinity, $(n-1)p'-n=(n-p)/(p-1)>0$.
These power bounds give finite upper bounds for $\|T_i\|_{p'}$.

They also give computable upper bounds for
$\int_{|x|<\rho}T_i^{\mathsf T}\mathbb A^{-1}T_i\,dx$ with decay
powers $\rho^{n+p'}$ for $\ell=0$ and
$\rho^{p'(\alpha+2s_\ell-1)}$ for $\ell\ge1$; both exponents are
positive. The corresponding bound for
$\int_{|x|>R}T_i^{\mathsf T}\mathbb A^{-1}T_i\,dx$ has decay power
$R^{-((n-p)/(p-1)+2p'r_\ell)}$, again with a positive decay exponent.
We use these bounds in
Subsection~\ref{sec:superquadratic-effective}
to control the errors from restricting the
fields $T_i$ to an annulus.
Also $|We_i|\le\|b_i\|_\infty v^{q-1}$, giving the required
$L^{q'}$ bound even when $q\le2$.
The eigenvalue equation and
$L^2(W)$-orthonormality give the first identity
in \eqref{eq:stage4-weak-pairings}.
The second holds for $h\in\mathcal D_0$
by the cited operator equation.  This class is dense in $\dot W^{1,p}$:
flattening a smooth compactly supported function to its value at zero
on a ball of radius $\varepsilon$ changes its gradient norm by
$O(\varepsilon^{n/p})$.  Since $T_i\in L^{p'}$ and $We_i\in L^{q'}$,
H\"older's inequality and \eqref{eq:sobolev}
allow us to extend the identity to every
$h\in\dot W^{1,p}$.

Finally, \eqref{eq:stage4-imported-eigenvalues} gives
$\mu_{\ell,j}\ge\Lambda C_nj^2$.  The nonnegative radial derivative
term and $t(1-t)\le1/4$ give
$\mu_{\ell,j}\ge4\Lambda C_nd\ell(\ell+n-2)$.
Because $4d\le1$ and $n\ge2$, either $j>M$ or $\ell>M$ implies
$\mu_{\ell,j}>L$.  Parseval's identity and the spectral expansion
of $Q$ give \eqref{eq:stage4-spectral-tail}.
\end{proof}

\subsection{The case \texorpdfstring{$p\ge2$}{p >= 2}}
\label{sec:superquadratic-effective}

Fix $p\ge2$ and an integer $n>p$; then $q>2$. Set
$\mu_*=\mu_\dagger$ as in \eqref{eq:stage4-weaker-gap}.

\begin{theorem}\label{thm:super-local}
Let $p\ge2$ and let $n>p$ be an integer.
There exist $t_*(n,p)>0$ and $c_{\mathrm{orth}}(n,p)>0$
such that, for every $h\in\dot W^{1,p}(\R^n)$ satisfying
\[
 \|\nabla h\|_p\le t_*,\qquad
 \int Wh e_i=0\quad\text{for every tangent eigenfunction }e_i,
\]
we have
\[
 F(v+h)\ge c_{\mathrm{orth}}\|\nabla h\|_p^p.
\]
\end{theorem}

With $B_p$ as defined in \eqref{eq:vector-bregman},
the Bregman bound \eqref{eq:full-bregman-bound} gives
\begin{equation}\label{eq:super-elementary-convexity}
 B_p(a,b)\ge\gamma_p(|a|+|b|)^{p-2}|b|^2.
\end{equation}

Set
\[
 \mathfrak B(h)=\int B_p(\nabla v,\nabla h),\qquad
 K_p^{\mathrm{B}}=\frac{p-1}{\gamma_p}.
\]
By \eqref{eq:super-elementary-convexity}, we have
\begin{equation}\label{eq:super-form-domination}
 \mathfrak B(h)\ge \gamma_p\|\nabla h\|_p^p,
 \qquad Q(h)\le K_p^{\mathrm{B}}\mathfrak B(h).
\end{equation}
For $p\ge2$, every $h\in\dot W^{1,p}(\mathbb R^n)$ belongs
to $D(Q)$. Indeed, H\"older's inequality gives
\[
 Q(h)\le(p-1)\mathcal G^{p-2}\|\nabla h\|_p^2,
 \qquad
 \int Wh^2\le\mathcal V^{(q-2)/q}S_{n,p}^{-2}\|\nabla h\|_p^2.
\]
Approximate $h$ strongly in $\dot W^{1,p}$ by compactly supported
smooth functions, and flatten each approximant to a constant on a
shrinking ball at zero. The latter modification converges in
$\dot W^{1,p}$. The two displayed bounds give convergence in the
weighted form norm.

For fixed $a$, the convex conjugate of $b\mapsto B_p(a,b)$ is
\begin{align}\label{eq:super-conjugate}
 B_p(a,\cdot)^*(Y)
 ={}&(p-1)p^{-p'}
 \left(|p|a|^{p-2}a+Y|^{p'}-|p|a|^{p-2}a|^{p'}\right)-a\cdot Y.
\end{align}
For $a\ne0$, its Hessian at zero is $\mathbb A(a)^{-1}/p$, where
\[
 \mathbb A(a)=|a|^{p-2}\left(I+(p-2)\frac{a\otimes a}{|a|^2}\right).
\]
For $1<s\le2$ and $z\ne0$, we have
$B_s(z,Y)\le20|z|^{s-2}|Y|^2$.
When $|Y|\le|z|/2$, this follows from the
Hessian of $|\cdot|^s$.
When $|Y|>|z|/2$, we use
\[
 B_s(z,Y)\le
 (|z|+|Y|)^s+|z|^s+s|z|^{s-1}|Y|
\]
and $s\le2$.
Substituting $s=p'$ in~\eqref{eq:super-conjugate} yields, for $a\ne0$,
\begin{equation}\label{eq:super-global-dual}
 B_p(a,\cdot)^*(Y)\le L_pY^T\mathbb A(a)^{-1}Y,
 \qquad L_p=20\frac{(p-1)^2}{p^2}.
\end{equation}

For the sharper bound needed below, fix
$0<\eta<1$ and set
\[
 \vartheta=\min\{1/4,\eta/[100(p-1)]\}.
\]
If $a\ne0$ and $|Y|\le\vartheta p|a|^{p-1}$, comparison of the radial factor and
the rank-one projector in the Hessian gives the sharper bound
\begin{equation}\label{eq:super-local-dual}
 B_p(a,\cdot)^*(Y)\le\frac{1+\eta/2}{2p}Y^T\mathbb A(a)^{-1}Y.
\end{equation}
To prove this, set $z_0=p|a|^{p-2}a$ and
$z=z_0+tY$ for $0\le t\le1$.
Since $|Y|\le\vartheta|z_0|$ and
$\vartheta\le1/4$, we have
$1/2\le|z|/|z_0|\le3/2$.
On this range, the relative change in
$|z|^{p'-2}$ and the operator-norm change
in $z\otimes z/|z|^2$ are each at most
$4|z-z_0|/|z_0|$.
The smallest eigenvalue of the normalized
Hessian at $z_0$ is $p'-1=1/(p-1)$.
Thus our choice of $\vartheta$ bounds the
relative change in the Hessian by $\eta/2$.
Integrating the Hessian along the segment
from $z_0$ to $z_0+Y$ proves
\eqref{eq:super-local-dual}.

\begin{lemma}\label{lem:super-finite-dual}
Let $I$ be a finite set of indices corresponding to
non-tangent eigenfunctions in the $L^2(W)$-orthonormal basis
of \cref{prop:stage4-spectral-toolkit}, and let $0<\eta<1$.
There exists an effectively determinable $t_0>0$ such that,
for every $h\in\dot W^{1,p}(\R^n)$ with
$\|\nabla h\|_p\le t_0$,
\[
 \mathfrak B(h)\ge\frac p{2(1+\eta)}
       \sum_{i\in I}\mu_i\left|\int We_i h\right|^2.
\]
\end{lemma}

\begin{proof}
We first choose $0<\rho<1<R$ so that, with
$\Omega=\{\rho<|x|<R\}$ and $T_i=\mathbb A\nabla e_i$,
\begin{equation}\label{eq:dual-tail-stopping-test}
 \sum_{i\in I}\int_{\Omega^c}T_i^T\mathbb A^{-1}T_i
       \le\frac{\eta\mu_*}{8pL_p}.
\end{equation}
By the bounds for $T_i$ proved in
\cref{prop:stage4-spectral-toolkit}, the integrals
over $|x|<\rho$ and $|x|>R$ are bounded by
finite sums of positive powers of $\rho$ and
$R^{-1}$, respectively.

For angular degree $\ell\ge1$, the inner
convergence exponent is
$p'(\alpha+2s_\ell-1)>0$; for $\ell=0$,
it is $n+p'>0$.
The outer exponent is at least
$(n-p)/(p-1)+2p'r_\ell>0$.
We can therefore choose $\rho$ and $R$
in finitely many dyadic steps so that
\eqref{eq:dual-tail-stopping-test} holds.

Let $B_I\ge1$ be an explicit upper bound for
$(\sum_i|T_i|^2)^{1/2}$ on $\Omega$. Such a bound is obtained by
bounding each finite monomial at the endpoints of $[\rho,R]$ and each
spherical polynomial by the absolute sum of its coefficients.
An explicit positive lower
bound for $|\nabla v|$ there is
\[
 a_{\min}=\frac{n-p}{p-1}\,
       \rho^{p'-1}(1+R^{p'})^{-n/p}.
\]
Put $c_0=\vartheta p a_{\min}^{p-1}/B_I$.
For $Y=\sum_i c_iT_i$ and $|c|\le c_0$,
\eqref{eq:super-local-dual} applies on $\Omega$. On its complement,
\eqref{eq:dual-tail-stopping-test} and $\mu_i\ge\mu_*$ imply
\[
 \int_{\Omega^c}Y^T\mathbb A^{-1}Y
 \le |c|^2\sum_i\int_{\Omega^c}T_i^T\mathbb A^{-1}T_i
 \le\frac{\eta}{8pL_p}\sum_i\mu_i c_i^2.
\]
Combining \eqref{eq:super-local-dual} on $\Omega$
with \eqref{eq:super-global-dual} on $\Omega^c$,
and using \eqref{eq:dual-tail-stopping-test},
we obtain
\[
 \int B_p(\nabla v,\cdot)^*(Y)
 \le\frac{1+\eta}{2p}\sum_i\mu_i c_i^2.
\]
For $h\in\dot W^{1,p}$, set
$a_i=\int We_i h$.
Bessel's inequality, H\"older's inequality,
and \eqref{eq:sobolev} give
\[
 \sum_{i\in I}|a_i|^2
 \le\int Wh^2
 \le\mathcal V^{(q-2)/q}S_{n,p}^{-2}
       \|\nabla h\|_p^2.
\]
Therefore, if $\|\nabla h\|_p\le t_0$, where
\[
 t_0=\frac{c_0S_{n,p}}
 {p\mathcal V^{(q-2)/(2q)}},
\]
then the coefficients
$c_i=pa_i/(1+\eta)$ satisfy $|c|\le c_0$.
Applying the Fenchel--Young inequality with
$Y=\sum_{i\in I}c_iT_i$ and using
\eqref{eq:stage4-weak-pairings}
gives the required estimate.
\end{proof}

\begin{proof}[Proof of \cref{thm:super-local}]
Let $g=\mu_*-\mu_{\mathrm{tan}}>0$, and choose
\[
 \eta=\frac g{16\mu_*},\qquad
 L=\max\{2\mu_*,8pK_p^{\mathrm{B}}\mu_*^2/g\}.
\]
Choose an integer $M$ satisfying
\eqref{eq:stage4-short-cutoff}, and let $I$
consist of the indices of the non-tangent eigenfunctions
with $0\le\ell,j\le M$, including their
angular multiplicities.
Apply \cref{lem:super-finite-dual} with
this $I$ and $\eta$, and denote the
resulting radius by $t_0$.

By tangent orthogonality and Parseval's identity,
\eqref{eq:stage4-spectral-tail}
and \eqref{eq:super-form-domination} imply
\[
 \sum_{i\notin I}\left|\int We_i h\right|^2
 \le\frac{Q(h)}L\le\frac{K_p^{\mathrm{B}}}L\mathfrak B(h).
\]
Combining this estimate with
\cref{lem:super-finite-dual} and using
$\mu_i\ge\mu_*$ for $i\in I$, we obtain,
whenever $\|\nabla h\|_p\le t_0$,
\[
 \mathfrak B(h)\ge a\left(\int Wh^2-
                      \frac{K_p^{\mathrm{B}}}L\mathfrak B(h)\right),
 \qquad a=\frac{p\mu_*}{2(1+\eta)}.
\]
Since $aK_p^{\mathrm{B}}/L\le\eta$ and
\[
 \frac{\mu_*}{(1+\eta)^2}
 \ge\mu_*-\frac g8
 \ge\mu_{\mathrm{tan}}+\frac g2,
\]
we can absorb the last term into the
left-hand side and obtain
\begin{equation}\label{eq:super-nonlinear-gap}
 \mathfrak B(h)\ge\frac p2\left(\mu_{\mathrm{tan}}+\frac g2\right)\int Wh^2.
\end{equation}

Set $T=\mu_{\mathrm{tan}}+g/2$, $\theta=g/(8T)$, and
$\sigma=\min\{1/2,g/(8\mu_{\mathrm{tan}})\}$. Define
\begin{align*}
 \tau&=\min\{1/2,((1+\sigma)^{1/(q-2)}-1)/2\},\\
 C_{q,\sigma}&=2^{q-1}(1+\tau^{-q})+q\tau^{1-q}.
\end{align*}
Taylor's formula on $|s|\le\tau$ and the triangle inequality on its
complement give, for every real $s$,
\[
 |1+s|^q\le1+qs+\frac{q(q-1)}2(1+\sigma)s^2
                         +C_{q,\sigma}|s|^q.
\]
Using \eqref{eq:bubble-euler-lagrange}
and the concavity of $t\mapsto t^{p/q}$,
we obtain
\[
 F(v+h)\ge\mathfrak B(h)-\frac p2\mu_{\mathrm{tan}}(1+\sigma)\int Wh^2
                    -\frac pq\Lambda C_{q,\sigma}\|h\|_q^q.
\]
We apply \eqref{eq:super-nonlinear-gap}
to $(1-\theta)\mathfrak B(h)$ and the
first bound in \eqref{eq:super-form-domination}
to $\theta\mathfrak B(h)$.
Since
\[
 (1-\theta)T-(1+\sigma)\mu_{\mathrm{tan}}
 \ge g/4>0,
\]
we may discard the remaining nonnegative
term involving $\int Wh^2$.
By \eqref{eq:sobolev}, we obtain
\[
 F(v+h)\ge\theta \gamma_p\|\nabla h\|_p^p
        -\frac pq\Lambda C_{q,\sigma}S_{n,p}^{-q}\|\nabla h\|_p^q.
\]
Since $q>p$, the term of order
$\|\nabla h\|_p^q$ can be absorbed
by choosing $\|\nabla h\|_p$ sufficiently
small. More precisely, set
\begin{equation}\label{eq:super-final-radius}
 t_* =\min\left\{t_0,
 \left(\frac{\theta \gamma_p qS_{n,p}^q}
 {2p\Lambda C_{q,\sigma}}\right)^{1/(q-p)}\right\},
 \qquad c_{\mathrm{orth}}=\frac{\theta \gamma_p}{2}.
\end{equation}
With $t_*$ and $c_{\mathrm{orth}}$ chosen in
\eqref{eq:super-final-radius}, the preceding
estimate gives
$F(v+h)\ge c_{\mathrm{orth}}\|\nabla h\|_p^p$
whenever $\|\nabla h\|_p\le t_*$,
as claimed.
\end{proof}

\subsection{The case \texorpdfstring{$1<p<2$}{1 < p < 2}}
\label{sec:subquadratic-effective}

Fix $1<p<2$ and an integer $n>p$.
We use $q,v,W,\mathbb A$ and $\Lambda$
as in Subsection~\ref{sec:exact-spectrum}.
Set
\[
 k=p',\qquad m=\frac{n-p}{p},\qquad K=\frac{n-p}{p-1},\qquad
 a_v(x)=|\nabla v(x)|=K|x|^{k-1}(1+|x|^k)^{-n/p}.
\]
Thus $\mathbb A=a_v^{p-2}(I+(p-2)e\otimes e)$, where
$e=\nabla v/a_v$ off the origin. Write $S=S_{n,p}$ and
$\mathcal G=\|\nabla v\|_p$.
Set $\mu_*=\mu_\dagger>\mu_{\mathrm{tan}}$
as in \eqref{eq:stage4-weaker-gap}, and use
the spectral estimates in
\cref{prop:stage4-spectral-toolkit}.
For $q\le2$, we adapt the truncation argument
in \cite[proof of Lemma~3.4]{FZ}
to signed perturbations, keeping track
of the constants.

\begin{theorem}\label{thm:short-subquadratic-local}
Let $1<p<2$ and let $n>p$ be an integer.
There exist $\delta_0,c_0>0$, depending only on $n,p$, such that,
for every $h\in\dot W^{1,p}(\R^n)$ satisfying
\[
 \|\nabla h\|_p\le\delta_0,\qquad
 \int Wh e_j=0\quad\text{for every tangent eigenfunction }e_j,
\]
we have
\[
 F(v+h)\ge c_0\|\nabla h\|_p^2.
\]
\end{theorem}

\subsubsection{Gradient estimates and truncation}

For $a\ge0$, let
$F_a(z)=(a+|z|)^{p-2}|z|^2$,
with $F_0(0)=0$.
This function is convex and satisfies
$F_a(2z)\le4F_a(z)$. Consequently,
\begin{equation}\label{eq:sq-modular-subadditive}
F_a(z-z')\le2F_a(z)+2F_a(z').
\end{equation}
For $\psi\in\dot W^{1,p}$ and $0<t\le1$, put
\begin{equation}\label{eq:modular-energy}
\mathfrak E_t(\psi)
=\int(a_v+t|\nabla\psi|)^{p-2}|\nabla\psi|^2
=t^{-2}\int F_{a_v(x)}(t\nabla\psi)\,dx.
\end{equation}
H\"older's and Minkowski's inequalities show that
\begin{equation}\label{eq:sq-gradient-bound}
\mathfrak E_t(\psi)\le L
\quad\Longrightarrow\quad
\|\nabla\psi\|_p^2\le L(\mathcal G+t\|\nabla\psi\|_p)^{2-p},
\qquad
\|\nabla\psi\|_p\le B_L:=\mathcal G+(2^{2-p}L)^{1/p}.
\end{equation}
Indeed, if $\|\nabla\psi\|_p>\mathcal G$,
the first inequality gives
$\|\nabla\psi\|_p^p\le2^{2-p}L$;
otherwise the bound by $B_L$ is immediate.
In particular, \eqref{eq:sq-gradient-bound}
and the Sobolev inequality \eqref{eq:sobolev}
give, whenever $\mathfrak E_t(\psi)\le1$,
\begin{equation}\label{eq:sq-qbound}
\|\psi\|_q\le Q_0:=S^{-1}B_1.
\end{equation}

We truncate $\psi$ at the levels $\pm v/t$
by setting
\[
w=\operatorname{sgn}(\psi)\min\{|\psi|,v/t\},
\qquad
z=\psi-w,\qquad
\mathcal A=\{t|\psi|>v\}.
\]
The Sobolev chain rule for minima and maxima
gives $w,z\in\dot W^{1,p}$, since
$\psi,v/t\in\dot W^{1,p}$.
We claim that
\begin{equation}\label{eq:sq-clipping}
\mathfrak E_t(w)\le\mathfrak E_t(\psi),\qquad
\mathfrak E_t(z)\le4\mathfrak E_t(\psi).
\end{equation}
To prove the first inequality, set $\phi=t\psi$. On
$\mathcal A_+=\{\phi>v\}$, convexity gives
\[
\int_{\mathcal A_+}F_{a_v(x)}(\nabla\phi)
\ge\int_{\mathcal A_+}F_{a_v(x)}(\nabla v)
+\int_{\mathcal A_+}DF_{a_v(x)}(\nabla v)\cdot\nabla(\phi-v).
\]
Let $\eta_\pm=(\pm\phi-v)_+$.
The Sobolev chain rule gives
$\eta_\pm\in\dot W^{1,p}$.
Since $a_v^{p-2}\nabla v\in L^{p'}$,
truncation and density justify testing
\eqref{eq:bubble-euler-lagrange}
with $\eta_\pm$, giving
\[
 \int |\nabla v|^{p-2}\nabla v\cdot\nabla\eta_\pm
 =\Lambda\int v^{q-1}\eta_\pm\ge0.
\]
Since
$DF_{a_v(x)}(\nabla v)
=2^{p-3}(p+2)a_v^{p-2}\nabla v$,
the last integral in the convexity
inequality is nonnegative on $\mathcal A_+$.
The same argument with $-\phi$ applies
on $\mathcal A_-$.
Outside $\mathcal A$, we have $w=\psi$.
Dividing by $t^2$ proves the first inequality
in \eqref{eq:sq-clipping}, while the second
follows from \eqref{eq:sq-modular-subadditive}.
Since $|w|,|z|\le|\psi|$,
\eqref{eq:sq-qbound} also gives
$\|w\|_q,\|z\|_q\le Q_0$ whenever
$\mathfrak E_t(\psi)\le1$.

\subsubsection{Weighted estimates and finite spectral approximation}

We first assume that $q\le2$. Necessarily $n\ge3$. Set
\[
b=\frac{n(2-p)}p-2\ge0,\qquad
B=W(1+r^k)^{p-1},\qquad r=|x|.
\]
Here $W=v^{q-2}=(1+r^k)^b$ is the weight from
Subsection~\ref{sec:exact-spectrum}.
For every $x\in\mathbb R^n$,
$\eta>0$, and $a,c\ge0$ satisfying
$ta\le v(x)$, we have
\begin{equation}\label{eq:sq-pointwise-transfer}
B a^{2-p}c^p\le\eta Wa^2+
\frac{C_\eta}{1+r^k}(a_v+tc)^{p-2}c^2,
\qquad C_\eta=K^{2-p}\eta^{-(2-p)/p}+1.
\end{equation}
The term $(a_v(x)+tc)^{p-2}c^2$ is understood as zero
when $a_v(x)=c=0$, as in the definition of $F_0(0)$.
For $a,c>0$, if
\[
 \frac ca\le\eta^{1/p}(1+r^k)^{-1/k},
\]
then the first term on the right-hand side
of \eqref{eq:sq-pointwise-transfer} suffices.
Otherwise, set $\xi=c/a$.
Since $0<2-p<1$,
\[
B\xi^{p-2}(a_v+ta\xi)^{2-p}
\le B a_v^{2-p}\xi^{p-2}+B(ta)^{2-p}.
\]
The two summands are bounded respectively by
\[
K^{2-p}\eta^{-(2-p)/p}
r^{(2-p)/(p-1)}(1+r^k)^{-2/p}
\le\frac{K^{2-p}\eta^{-(2-p)/p}}{1+r^k},
\qquad Bv^{2-p}=\frac1{1+r^k}.
\]
These estimates prove
\eqref{eq:sq-pointwise-transfer} for $a,c>0$.
The cases $a=0$ and $c=0$ follow directly.

For $f\in C_c^\infty(\mathbb R^n)$,
integration by parts, using
$\operatorname{div}(xW)\ge nW$
and H\"older's inequality, gives
\begin{equation}\label{eq:sq-hardy}
\int W|f|^p\le(p/n)^p\int r^pW|\nabla f|^p
\le(p/n)^p\int B|\nabla f|^p.
\end{equation}
The bound $|w|\le v/t$ first gives
$\int Ww^2\le t^{-2}\int v^q<\infty$.
By the chain rule and \eqref{eq:sq-pointwise-transfer},
$f=|w|^{2/p}$ has finite weighted gradient energy:
\[
 \int B|\nabla f|^p
 \le(2/p)^p\left(\eta\int Ww^2+C_\eta\mathfrak E_t(w)\right)<\infty.
\]
We can apply \eqref{eq:sq-hardy} to
$f=|w|^{2/p}$ by approximation.
Indeed, near the origin we use
$|w|\le1/t$ and $n>p$.
For $r\ge1$, we use
$\int Ww^2<\infty$ and $B/r^p\lesssim W$.
These bounds show that the spatial cutoff
errors tend to zero and give
$f\in\dot W^{1,p}$.

Applying \eqref{eq:sq-hardy} and \eqref{eq:sq-pointwise-transfer}
when $\mathfrak E_t(\psi)\le1$ gives
\[
 \int Ww^2\le(2/n)^p\left(\eta\int Ww^2+C_\eta\right).
\]
Define
\[
\eta_0=\frac{n^p}{2^{p+1}},\qquad
C_0=2(2/n)^pC_{\eta_0},\qquad
C_F=(2/p)^p(\eta_0C_0+C_{\eta_0}).
\]
Absorbing the first term yields
\begin{equation}\label{eq:sq-clipped-bounds}
\begin{gathered}
\int Ww^2\le C_0,\qquad
\int B|\nabla(|w|^{2/p})|^p\le C_F,\\
\|w\|_{r_*}\le Q_*:=(S^{-p}C_F)^{1/2},
\quad r_*:=2q/p>2.
\end{gathered}
\end{equation}
The last assertion follows from $B\ge1$ and the Sobolev inequality
\eqref{eq:sobolev} applied to $|w|^{2/p}$.

For fixed $n,p$, the following bounds hold uniformly for $0<t\le1$
and $\mathfrak E_t(\psi)\le1$. Set
\[
a_*=kb,\quad
C_{\rm ext}=2^b\left(\frac2{n+a_*}\right)^p,
\quad \eta_1=(2C_{\rm ext})^{-1},\quad
C_T=2C_{\rm ext}C_{\eta_1},\quad
C_I=2^bQ_*^2\omega_n^{1-2/r_*}.
\]
Then
\begin{equation}\label{eq:sq-weighted-tails}
\int_{|x|>R}Ww^2\le C_TR^{-k}\quad(R\ge1),\qquad
\int_{|x|<\rho}Ww^2\le C_I\rho^p\quad(0<\rho\le1).
\end{equation}
To prove the exterior estimate, we integrate
by parts on $\{|x|>R\}$, using
\[
 \operatorname{div}(xr^{a_*})
 =(n+a_*)r^{a_*}.
\]
The boundary term at $r=R$ is nonpositive,
so
\[
\int_{|x|>R}r^{a_*}|f|^p
\le\left(\frac p{n+a_*}\right)^p
\int_{|x|>R}r^{a_*+p}|\nabla f|^p.
\]
We apply this inequality to $f=|w|^{2/p}$.
For $r\ge1$, the weights $r^{a_*}$ and $W$
are comparable. Using
\eqref{eq:sq-pointwise-transfer} and the
choice $\eta_1=(2C_{\rm ext})^{-1}$,
we absorb half of the term involving
$Ww^2$. Since $\mathfrak E_t(w)\le\mathfrak E_t(\psi)\le1$,
the remaining term is bounded by $C_TR^{-k}$.
The inner estimate follows from
\eqref{eq:sq-clipped-bounds} and H\"older's inequality, because
$n(1-2/r_*)=p$.

For $\psi\in\dot W^{1,p}$ and $0<t\leq1$, put
\begin{equation}\label{eq:scalar-remainder-definition}
 N_t(\psi)=t^{-2}\int B_q(v,t\psi),\qquad c_q=q(q-1)/2.
\end{equation}

\begin{lemma}\label{lem:short-discarded}
Let $1<p<2$, let $n>p$ be an integer, and assume $q\le2$.
There exists $C=C(n,p)>0$ such that for every $0<t\le1$ and
$\psi\in\dot W^{1,p}(\R^n)$ with $\mathfrak E_t(\psi)\le1$,
where $w$ and $z=\psi-w$ are the truncation and remainder
defined above, we have
\begin{equation}\label{eq:short-discarded-power}
 t^{q-2}\|z\|_q^q\le C t^\alpha,
 \qquad \alpha=\frac{4(p-1)}{n-p}>0.
\end{equation}
Moreover,
\begin{equation}\label{eq:short-direct-clipping}
 |N_t(\psi)-N_t(w)|\le C t^\gamma,
 \qquad \gamma=\frac\alpha q=\frac{4(p-1)}{np}>0.
\end{equation}
\end{lemma}
\begin{proof}
Since $|w|=v/t$ on $\mathcal A$,
\begin{equation}\label{eq:short-clipped-mass}
 \int_{\mathcal A}v^q=t^2\int_{\mathcal A}Ww^2\leq C_0t^2.
\end{equation}
Let $\vartheta=(n-2)/n\in(0,1)$. H\"older's inequality gives
\[
 \int_{\mathcal A} a_v^p\leq K_H\left(\int_{\mathcal A}v^q\right)^\vartheta,
 \qquad
 K_H=K^p\left[\frac{n\omega_n}{k}
 \mathrm B\left(\frac nk+\frac n2,\frac n2-\frac nk\right)
                  \right]^{2/n}<\infty.
\]
Indeed, $K_H=\|(a_v^p/v^{q\vartheta})\|_{n/2}$, and $k>2$
makes both Beta parameters positive. Since $\nabla z$ vanishes almost everywhere
outside $\mathcal A$, we split $\mathcal A$
into the regions $t|\nabla z|\le a_v$ and
$t|\nabla z|>a_v$.
Using H\"older's inequality on the first
region and $\mathfrak E_t(z)\le4$
on both regions, we obtain
\[
 \|\nabla z\|_p^p
 \leq(4\,2^{2-p})^{p/2}\left(\int_{\mathcal A} a_v^p\right)^{1-p/2}
                  +4\,2^{2-p}t^{2-p}
 \leq C t^{\vartheta(2-p)}.
\]
Applying \eqref{eq:sobolev}, we obtain
\eqref{eq:short-discarded-power}, since
\[
 q-2+\frac qp\vartheta(2-p)
 =\frac{4(p-1)}{n-p}.
\]

Since $w=\psi$ outside $\mathcal A$,
we only need to estimate the difference
on $\mathcal A$.
There, $tw/v=\operatorname{sgn}(\psi)$,
while $t\psi/v$ has the same sign and
absolute value greater than one.
Thus every intermediate value satisfies
$|s|\ge1$, where
\[
 |\partial_s B_q(1,s)|\le6|s|^{q-1}.
\]
Set $A_q=\int_{\mathcal A}|w|^q$ and
$Z_q=\int|z|^q$.
Using $(1+a)^{q-1}\le1+a^{q-1}$
for $a\ge0$, we obtain
\[
 |N_t(\psi)-N_t(w)|
 \leq6t^{q-2}(A_q^{1-1/q}Z_q^{1/q}+Z_q).
\]
By \eqref{eq:short-clipped-mass},
\[
 t^{q-2}A_q
 =\int_{\mathcal A}Ww^2\le C_0.
\]
Set $X=t^{q-2}Z_q$.
The preceding bound is then at most
$6(C_0^{1-1/q}X^{1/q}+X)$.
By \eqref{eq:short-discarded-power},
we have $X\le Ct^\alpha$.
Since $0<t\le1$, this is bounded by
a constant depending only on $n,p$
times $t^{\alpha/q}$, proving
\eqref{eq:short-direct-clipping}.
\end{proof}

We now return to the full range $n>p$, allowing both $q\le2$ and $q>2$.

\begin{proposition}\label{prop:short-scalar-capture}
For every $\epsilon>0$, there exist a finite set $J$
of spectral indices including the indices corresponding to $\mathcal T$,
and $t_\epsilon\in(0,1]$
such that for every $\psi\in\dot W^{1,p}(\R^n)$ and
$0<t\le t_\epsilon$ satisfying $\mathfrak E_t(\psi)\le1$,
\begin{equation}\label{eq:short-scalar-capture}
 N_t(\psi)\le\frac{q(q-1)}2
       \sum_{j\in J}\left|\int W\psi e_j\right|^2+\epsilon.
\end{equation}
No assumption that $\psi\in L^2(W)$ is required.
\end{proposition}
\begin{proof}
Recall $N_t$ and $\mathfrak E_t$ from
\eqref{eq:scalar-remainder-definition}
and \eqref{eq:modular-energy}.

We first consider $q\le2$.
By \eqref{eq:sq-clipped-bounds},
$w\in L^2(W)$. We will estimate
$\|(I-P_J)w\|_{L^2(W)}$ uniformly
over $\mathfrak E_t(\psi)\le1$.
For $R\geq2$ and
$0<\rho\leq1/2$, let $\chi$ be a radial cutoff
equal to one on $2\rho\leq r\leq R$, supported on
$\rho\leq r\leq2R$, with $|\nabla\chi|\leq2/\rho+2/R$.
Set $u=\chi w$ and $g=\varrho_s*u$, where
\[
\varrho_s(x)=s^{-n}\prod_{i=1}^n(1-|x_i|/s)_+,
\qquad s\sqrt n\leq\rho/2.
\]
Then
\[
 \|\nabla u\|_p\leq B_1+(2/\rho+2/R)Q_0|B_{2R}|^{1/p-1/q},
 \quad \|u\|_{r_*}\leq Q_*,\quad
 \|u\|_p\leq Q_0|B_{2R}|^{1/p-1/q}.
\]
By Young's convolution inequality,
$g\in W^{1,\infty}_c(\mathbb R^n)$
and is supported in
$\{\rho/2\le|x|\le3R\}$.
Since $\mathbb A$ and $W$ are bounded
on a slightly larger annulus, smooth
approximation within that annulus
shows that $g\in D(Q)$.

For $L>0$, choose an integer $M$
satisfying \eqref{eq:stage4-short-cutoff},
and let $J$ contain the indices of all eigenfunctions
with $0\le\ell,j\le M$, including
their angular multiplicities.
Then $P_J=P_M$ and every omitted
eigenvalue exceeds $L$. Since $M\ge1$, $P_J$ acts as
the identity on $\mathcal T$.

Interpolating the translation estimate
in $L^p$ with the $L^{r_*}$ bound,
and applying \eqref{eq:stage4-spectral-tail}
to $g$, we obtain
\begin{equation}\label{eq:short-three-errors}
\begin{aligned}
\|(I-P_J)w\|_{L^2(W)}
&\leq\bigl[C_I(2\rho)^p+C_TR^{-k}\bigr]^{1/2}
       +C_{\rho,R}s^\theta\\
&\quad+C'_{\rho,R}s^{-1-n/p}L^{-1/2},
\qquad \theta=\frac{1/2-1/r_*}{1/p-1/r_*}>0.
\end{aligned}
\end{equation}
To verify the last two terms, use
\[
 \|u-g\|_2\leq(s\sqrt n\|\nabla u\|_p)^\theta
                         (2Q_*)^{1-\theta},\qquad
 \|\nabla\varrho_s\|_k\leq\sqrt n\,2^{n/k}s^{-1-n/p},
\]
and
\[
 \|(I-P_J)g\|_{L^2(W)}^2\leq L^{-1}\int\mathbb A\nabla g\cdot\nabla g
 \leq L^{-1}\max_{\rho/2\leq r\leq3R}a_v^{p-2}
       |B_{3R}|\|\nabla\varrho_s\|_k^2\|u\|_p^2.
\]
Here $W$ is nondecreasing, and since
$p<2$ and $a_v$ has a single maximum,
the maximum of $a_v^{p-2}$ on the
annulus occurs at one of its endpoints.
Thus all the constants in
\eqref{eq:short-three-errors} can be
determined from the preceding bounds.
Contraction of $I-P_J$
proves \eqref{eq:short-three-errors}. The positive powers permit
successive finite choices of $R$, $\rho$, $s$, and $L$, with $J$
then fixed by \eqref{eq:stage4-short-cutoff}. 

We next replace the finite spectral coefficients of $w$ by those
of $\psi$, without assuming $\psi\in L^2(W)$. Namely, $We_j=v^{q-1}b_j$ and $v<t|\psi|$ on $\mathcal A$ give
\[
 \left|\int Wze_j\right|\leq\|b_j\|_\infty Q_0^q t^{q-1}.
\]
Since $|\int Wwe_j|\leq\sqrt{C_0}$, the difference between the
two sums of squared moments is at most $D_Jt^{q-1}$, with
\[
 D_J=\sum_{j\in J}(2\sqrt{C_0}D_j+D_j^2),
 \qquad D_j=\|b_j\|_\infty Q_0^q.
\]

It remains to compare the scalar remainder with $c_q\|w\|_{L^2(W)}^2$.
If $q=2$, the identity $B_2(1,s)=s^2$ gives
$N_t(w)=\|w\|_{L^2(W)}^2$ exactly; the two scalar Taylor-error terms
below can therefore be omitted.

For $q<2$, use the following Taylor bound.
For $0<a\leq1$ and $|s|\leq1$,
\begin{equation}\label{eq:short-scalar-taylor}
 |B_q(1,s)-c_qs^2|\leq64|s|^{2+a}.
\end{equation}
Take $a=\min\{1,(r_*-2)/2\}>0$.
Since $|tw|\leq v$, applying
\eqref{eq:short-scalar-taylor} on $B_R$
and its complement gives, for every $R\geq1$,
\begin{align}
 |N_t(w)-c_q\|w\|_{L^2(W)}^2|
 &\leq64 C_TR^{-k}+C_{R}t^a,
 \label{eq:short-scalar-quadratic}\\
 C_{R}
 &=64(1+R^k)^{b+ma}
       |B_{R}|^{1-(2+a)/r_*}Q_*^{2+a}.\nonumber
\end{align}
Indeed, on $\{|x|>R\}$, the bound
$|tw/v|\le1$ allows us to control
the Taylor error by\linebreak
$64\int_{|x|>R}Ww^2$,
which is bounded by
\eqref{eq:sq-weighted-tails}.
On $B_R$, we use H\"older's inequality
and $\|w\|_{r_*}\le Q_*$, noting that
$2+a\le r_*$.

For $q<2$, combining
\eqref{eq:short-three-errors},
\eqref{eq:short-scalar-quadratic},
Parseval's identity for $w\in L^2(W)$,
and \cref{lem:short-discarded}, we obtain
\[
 N_t(\psi)\leq c_q\sum_{j\in J}\left|\int W\psi e_j\right|^2
       +c_q\|(I-P_J)w\|_{L^2(W)}^2
       +64C_TR^{-k}+C_{R}t^a+Ct^\gamma+c_qD_Jt^{q-1}.
\]
When $q=2$, the same estimate holds
without the two scalar Taylor-error terms
$64C_TR^{-k}$ and $C_Rt^a$.
Choose $R$ first so that the contribution of $C_TR^{-k}$
to the spectral error is sufficiently small, and, if $q<2$,
so that $64C_TR^{-k}<\epsilon/4$ as well.
Next choose $\rho$, then $s$, then
$L$ and its integer cutoff $J$ so that the spectral error is less
than $\epsilon/4$. Finally, choose $t_\epsilon\in(0,1]$
so that the terms
$Ct^\gamma$ and $c_qD_Jt^{q-1}$,
together with $C_Rt^a$ when $q<2$,
have sum at most $\epsilon/2$
for every $0<t\le t_\epsilon$.
This proves \eqref{eq:short-scalar-capture}
when $q\le2$.

We now consider $q>2$.
In this case, H\"older's inequality
gives $\psi\in L^2(W)$.
More precisely, for every measurable
set $\Omega\subset\mathbb R^n$,
\[
 \int_\Omega W\psi^2\,dx
 \leq\|\psi\|_q^2
       \left(\int_\Omega v^q\,dx\right)^{1-2/q}
 \leq Q_0^2\left(\int_\Omega v^q\,dx\right)^{1-2/q}.
\]
For the Sobolev extremal $v$, we have
\[
 \int_{|x|<\rho}v^q\leq\omega_n\rho^n,\qquad
 \int_{|x|>R}v^q\leq\frac{n\omega_n}{np'-n}R^{n-np'}
 \quad(0<\rho\leq1\leq R).
\]
Consequently,
\begin{align*}
 \int_{|x|<\rho}W\psi^2
 &\le Q_0^2(\omega_n\rho^n)^{1-2/q},\\
 \int_{|x|>R}W\psi^2
 &\le Q_0^2\left(\frac{n\omega_n}{np'-n}\right)^{1-2/q}
       R^{(n-np')(1-2/q)}.
\end{align*}
The inner power is positive and the outer power is negative.
We use the same radial cutoff $\chi$
and convolution kernel $\varrho_s$,
now with $u=\chi\psi$, $r_*=q$,
and $Q_0$ in place of $Q_*$.
The translation and convolution
estimates used above remain valid with
\[
 \theta=\frac{1/2-1/q}{1/p-1/q}>0.
\]
The smoothed function belongs to $D(Q)$ by the same
approximation within an annulus.
Choosing $R$, then $\rho$, then $s$,
and finally $L$ and its finite spectral
set $J$, we can make
\[
 \|(I-P_J)\psi\|_{L^2(W)}^2
 =\|\psi\|_{L^2(W)}^2
  -\sum_{j\in J}\left|\int W\psi e_j\right|^2
\]
arbitrarily small, uniformly over
$\mathfrak E_t(\psi)\le1$.
The equality follows from Parseval's identity.

Put $a=\min\{1,q-2\}>0$. Taylor's theorem for $|s|\le1/2$
and a direct estimate for $|s|>1/2$
give
\[
 |B_q(1,s)-c_qs^2|\leq C_q(|s|^{2+a}+|s|^q),
\]
where one may take
\[
 C_q=\max\left\{\frac{q(q-1)|q-2|}{6}
       \max\{2^{3-q},(3/2)^{q-3}\},\,
       3^q+2^q+q2^{q-1}+c_q2^{q-2}\right\}.
\]
Since $q-2-a\geq0$, H\"older and $\|\psi\|_q\leq Q_0$ imply
\[
 |N_t(\psi)-c_q\|\psi\|_{L^2(W)}^2|
 \leq C_qt^a\bigl(\|v\|_q^{q-2-a}Q_0^{2+a}+Q_0^q\bigr).
\]
We first choose the finite spectral
set $J$ so that
\[
 c_q\|(I-P_J)\psi\|_{L^2(W)}^2
 \le\epsilon/2
\]
uniformly over $\mathfrak E_t(\psi)\le1$.
With $J$ fixed, choose
$t_\epsilon\in(0,1]$ so that
\[
 C_qt_\epsilon^a
 \bigl(\|v\|_q^{q-2-a}Q_0^{2+a}+Q_0^q\bigr)
 \le\epsilon/2.
\]
The preceding estimates then prove
\eqref{eq:short-scalar-capture}
for $q>2$.
\end{proof}

\subsubsection{Estimates using finitely many eigenfunctions}

For $a,z\in\mathbb R^n$, set $G_a(z)=B_p(a,z)$.
By \cref{lem:eff-bregman},
\begin{equation}\label{eq:sq-bregman-lower}
G_a(z)\ge\gamma_p(|a|+|z|)^{p-2}|z|^2,
\qquad \gamma_p=\frac{p(p-1)}2.
\end{equation}

For every $\eta>0$, there is an explicit
constant $C(p,\eta)>0$ such that, for $a\ne0$,
the convex conjugate satisfies
\begin{equation}\label{eq:sq-dual-taylor}
G_a^*(Y)\le\frac{1+\eta}{2p}|a|^{2-p}Y^{\mathsf T}M_a^{-1}Y
+C(p,\eta)|Y|^k,
\quad M_a=I+(p-2)\frac{a\otimes a}{|a|^2}.
\end{equation}
The convex conjugate of $|z|^p$ is $c|\xi|^k$, where
$c=(p-1)p^{-k}$. Thus $G_a^*$ is the Bregman remainder of this function
at $\xi_0=p|a|^{p-2}a$. On $|\xi-\xi_0|\le|\xi_0|/2$, direct
differentiation bounds the relative Hessian error by
$H_k|\xi-\xi_0|/|\xi_0|$, where
\[
H_k=(k-1)(k-2)\max\{2^{3-k},(3/2)^{k-3}\}+4(k-2).
\]
Set $\tau=\min\{1/2,\eta/H_k\}$.
When $|Y|\le\tau|\xi_0|$, we integrate
the Hessian estimate along the segment
from $\xi_0$ to $\xi_0+Y$ to obtain
the quadratic bound in
\eqref{eq:sq-dual-taylor}.
When $|Y|>\tau|\xi_0|$, the elementary
power inequality gives
\[
G_a^*(Y)\le c(2^k+1+k)\tau^{-k}|Y|^k.
\]
Thus $C(p,\eta)=c(2^k+1+k)\tau^{-k}$ suffices. For $a=0$, the conjugate is
$G_0^*(Y)=(p-1)p^{-p'}|Y|^{p'}$.

Suppose now that $\mathfrak E_t(\psi)=1$, $h=t\psi$, and put
$m_j=\int W\psi e_j$. For a finite set $J$ of spectral indices, choose
\[
g=(1+\eta)^{-1}\sum_{j\in J}m_je_j,
\qquad Y=pt\mathbb A\nabla g.
\]
By the weak identity \eqref{eq:stage4-weak-pairings},
\[
 t^{-2}\int Y\cdot\nabla h
 =p\int\mathbb A\nabla g\cdot\nabla\psi
 =\frac p{1+\eta}\sum_{j\in J}\mu_jm_j^2.
\]
By \eqref{eq:stage4-weak-pairings},
the quadratic term in the upper bound
\eqref{eq:sq-dual-taylor}, after division
by $t^2$, equals
\[
 \frac{p}{2(1+\eta)}
 \sum_{j\in J}\mu_jm_j^2.
\]
The Fenchel--Young inequality now gives
\begin{equation}\label{eq:sq-finite-dual}
\frac1{t^2}\int G_{\nabla v}(\nabla h)
\ge\frac p{2(1+\eta)}\sum_{j\in J}\mu_jm_j^2
-C_Jt^{k-2}.
\end{equation}
We may take
\[
C_J=C(p,\eta)p^k
\left[\frac{Q_0}{1+\eta}
\sum_{j\in J}\|We_j\|_{q'}\|\mathbb A\nabla e_j\|_k\right]^k.
\]
By \cref{prop:stage4-spectral-toolkit},
the norms in the definition of $C_J$
are finite and have effective upper bounds.
In particular, the integrals
$\int W\psi e_j$ and
$\int\mathbb A\nabla e_j\cdot\nabla\psi$
are well defined for
$\psi\in\dot W^{1,p}$.

\begin{proof}[Proof of \cref{thm:short-subquadratic-local}]
Choose
\[
 \sigma=\frac{\mu_*-\mu_{\mathrm{tan}}}{4\mu_*}\in(0,1/4),\qquad
 \eta=\sigma,\qquad
 \epsilon=\frac{q\sigma\gamma_p}{4p\Lambda}.
\]
Choose $J,t_\epsilon$ by \cref{prop:short-scalar-capture} and
then take the corresponding $C_J$ above. Let $h=t\psi$ be
tangent-orthogonal, with $0<t\le1$ and $\mathfrak E_t(\psi)=1$,
and write $m_j=\int W\psi e_j$.
Since $J$ contains the indices of all tangent eigenfunctions,
tangent orthogonality gives $m_j=0$
for these indices. For every remaining
index in $J$, we have $\mu_j\ge\mu_*$
by \cref{prop:stage4-spectral-toolkit}.
Thus \eqref{eq:sq-finite-dual} gives
the first inequality below.
The second follows from
\eqref{eq:sq-bregman-lower}
and $\mathfrak E_t(\psi)=1$:
\begin{align}
 t^{-2}\int B_p(\nabla v,\nabla h)
 &\geq\frac{p\mu_*}{2(1+\eta)}\sum_{j\in J}m_j^2-C_Jt^{k-2},
 \label{eq:short-finite-dual}\\
 t^{-2}\int B_p(\nabla v,\nabla h)
 &\geq\gamma_p\mathfrak E_t(\psi)=\gamma_p.
 \label{eq:short-reserve}
\end{align}

The bubble equation \eqref{eq:bubble-euler-lagrange} and concavity
of $s^{p/q}$ give
\[
 F(v+h)/t^2\geq t^{-2}\int B_p(\nabla v,\nabla h)
                               -(p/q)\Lambda N_t(\psi).
\]
We use \eqref{eq:short-finite-dual}
for $(1-\sigma)$ of the Bregman term
and \eqref{eq:short-reserve} for the
remaining $\sigma$ part.
By \eqref{eq:short-scalar-capture},
the coefficient of the finite spectral
sum is positive, since
\[
 \frac{1-\sigma}{1+\sigma}\mu_*-\mu_{\mathrm{tan}}
 \geq\mu_*-\mu_{\mathrm{tan}}-2\sigma\mu_*
 =\tfrac12(\mu_*-\mu_{\mathrm{tan}})>0.
\]

We may therefore discard this
nonnegative term. Set
\[
 t_0=\min\left\{1,\ t_\epsilon,
   \left(\frac{\sigma\gamma_p}{4(1+C_J)}\right)^{1/(k-2)}\right\}.
\]
Since $k=p'>2$, this gives $C_Jt_0^{k-2}\le\sigma\gamma_p/4$.
Also $(p/q)\Lambda\epsilon=\sigma\gamma_p/4$. Thus, for
$0<t\le t_0$,
\[
 F(v+h)\geq\tfrac12\sigma\gamma_p t^2.
\]

Finally, let $h\ne0$ be tangent-orthogonal,
and set
\[
 t^2=\int(a_v+|\nabla h|)^{p-2}|\nabla h|^2,\qquad\psi=h/t.
\]
Then $\mathfrak E_t(\psi)=1$ and H\"older's inequality gives
\[
 t\leq\|\nabla h\|_p^{p/2},\qquad
 t^2\geq\frac{\|\nabla h\|_p^2}
 {(\mathcal G+\|\nabla h\|_p)^{2-p}}.
\]
Set
\[
 \delta_0=\min\{1,t_0^{2/p}\},\qquad
 c_0=\frac{\sigma\gamma_p}{2(\mathcal G+1)^{2-p}}.
\]
If $\|\nabla h\|_p\le\delta_0$, then
$t\le t_0$ and
\[
 (\mathcal G+\|\nabla h\|_p)^{2-p}
 \le(\mathcal G+1)^{2-p}.
\]
Combining these bounds with the
preceding estimate gives
\[
 F(v+h)\ge c_0\|\nabla h\|_p^2,
\]
as required.
\end{proof}

\section{Global stability and proof of the main theorem}\label{sec:global-main}

In this section, we pass from local to global stability and complete
the proof of \cref{thm:main}. In Subsection~\ref{sec:entrance}, we use radial transport
 to prove \cref{prop:both-radial}. In
Subsection~\ref{sec:global}, we use polarization to obtain
global stability for nonnegative functions and then extend
the estimate to arbitrary functions in
\eqref{eq:short-real-global}. In
Subsection~\ref{sec:assembly}, we combine the local and global estimates
to obtain the lower bound, with $c_p$ in \eqref{eq:finalcp}.
Subsection~\ref{sec:upper} proves the upper bound using two widely
separated bubbles.

\subsection{Radial approximation}\label{sec:entrance}

\subsubsection{Transport identities}\label{sec:transport}

Fix $p>1$ and an integer $n>p$. We first consider a smooth,
nonnegative, radially nonincreasing function $u$, with
$\norm*{u}_q=1$, support $\overline B_R$, and $u>0$ in $B_R$.
We do not require $u$ to be strictly decreasing.
Set $D=\mathcal D(u)$ and assume $D\leq1/2$.
We follow the radial transport argument in
\cite[Section~2, (2.8)--(2.17)]{CFMP}.
Set
\[
 a_n=\frac{n-p}{p},\qquad b_n=\frac{n-p}{p-1}=p'a_n,
 \qquad q_\#=\frac{p(n-1)}{n-p},\qquad
 \varkappa_D=(1-D)^{-1/p}.
\]
Let
\[
 v_0(r)=A_0(1+r^{p'})^{-a_n},\qquad \norm*{v_0}_q=1,
 \qquad C_v=A_0^{-1/a_n},\qquad \lambda_0=b_nC_v.
\]
With $\tau=(p-1)/p$, the radial variable
$t=r^{p'}/(1+r^{p'})$ under the normalized measure $v_0^q\,dx$ has
distribution $\operatorname{Beta}(n\tau,n/p)$, as in
\eqref{eq:beta-measures}. Consequently,
\begin{equation}\label{eq:radial-K}
 K_n:=\int \abs*{x}^{p'}v_0^q\,dx
   =\frac{n(p-1)}{n-p},\qquad
 V_\#:=\int v_0^{q_\#}\,dx=C_v(1+K_n).
\end{equation}
Indeed, $v_0^{q_\#}=C_v(1+r^{p'})v_0^q$ and
\[
 -v_0'v_0^{-q/p}=\lambda_0r^{p'-1}.
\]
The equality case in \eqref{eq:sobolev} therefore gives
\begin{equation}\label{eq:radial-normalizations}
 E(v_0)=\lambda_0^pK_n=S_{n,p}^p,
 \qquad \frac{q_\#}{n}\lambda_0K_n=C_v(1+K_n)=V_\#.
\end{equation}

Define the increasing radial transport $T:(0,R)\to(0,\infty)$ by
\begin{equation}\label{eq:radial-CDF}
 \int_0^r u(s)^q s^{n-1}\,ds
   =\int_0^{T(r)}v_0(s)^q s^{n-1}\,ds,
 \qquad M=\frac{T}{r},\qquad \sigma=\frac{T'}{M}.
\end{equation}
Both cumulative distributions are strictly increasing on their open
supports, so $T'(r)>0$ and $\sigma>0$. Since both measures have
unit mass and $v_0>0$ on $\mathbb R^n$, $T(r)\to\infty$ as
$r\uparrow R$. Set
$\mathcal T(x)=T(\abs*{x})x/\abs*{x}$. The transport identities
(2.9)--(2.11) in the cited proof give
\begin{equation}\label{eq:radial-Jacobian}
 u(r)^q=v_0(T(r))^qM(r)^n\sigma(r),
\end{equation}
and
\begin{equation}\label{eq:radial-det}
 V_\#=\int u^{q_\#}M\sigma^{1/n}\,dx.
\end{equation}

The rescaled function $\widetilde u(s)=R^{a_n}u(Rs)$ has
support $[0,1]$ and preserves the $L^q$ normalization.
Applying the integration-by-parts argument of
\cite[(2.12)--(2.15), pp.~1110--1111]{CFMP} to
$\widetilde u$ with $\widetilde T(s)=T(Rs)$, and using
$a_nq_\#+1-n=0$, we obtain
\begin{equation}\label{eq:radial-IBP}
 \int u^{q_\#}\operatorname{div}\mathcal T\,dx
      =q_\#\int T(-u')u^{q_\#-1}\,dx<\infty.
\end{equation}
The integration-by-parts argument in
\cite[(2.13)--(2.15)]{CFMP} does not require $u$
to be strictly decreasing. The integrability estimate there
ensures that the boundary term vanishes.

Define $\zeta(s)=s+n-1-ns^{1/n}$ and
\begin{equation}\label{eq:transport-defects}
 h_n(s)=s^{1-1/n}+(n-1)s^{-1/n}-n,
 \qquad \mathcal H_p(z)=z^p-pz+p-1.
\end{equation}
The functions $h_n$ and $\mathcal H_p$ are nonnegative and vanish only at $1$.
By \cref{eq:radial-normalizations},
$q_\# S_{n,p}K_n^{1/p'}=nV_\#$. The deficit in the cited proof is
$\delta(u)=E(u)^{1/p}/S_{n,p}-1=\varkappa_D-1$.
Thus its remainder estimate (2.17) reads
\[
 \int u^{q_\#}M\zeta(\sigma)\,dx
 \leq nV_\#(\varkappa_D-1).
\]
By \cref{eq:radial-Jacobian},
\[
 u^{q_\#-q}=u^{-1/a_n}
       =C_v(1+T^{p'})M^{-1}\sigma^{-1/n},
\]
so $h_n(s)=s^{-1/n}\zeta(s)$ gives
\begin{equation}\label{eq:radial-AMGM-exact}
 \frac1n\int u^{q_\#}M\zeta(\sigma)\,dx
 =\frac{C_v}{n}\int_{\R^n}(1+T^{p'})h_n(\sigma)u^q\,dx.
\end{equation}
Combining \eqref{eq:radial-AMGM-exact}
with the preceding remainder estimate,
we obtain
\begin{equation}\label{eq:transport-AMGM}
 \int_{\R^n}(1+T^{p'})h_n(\sigma)u^q\,dx
 \leq n(1+K_n)(\varkappa_D-1).
\end{equation}

We next estimate the remainder in H\"older's inequality. Set
\[
 A=(-u')u^{-q/p},\qquad
 \Lambda=\frac{E(u)^{1/p}}{K_n^{1/p}}=\varkappa_D\lambda_0,
 \qquad \rho=\frac{A}{\Lambda T^{p'-1}}.
\]
These quantities are defined for $0<r<R$.
On this interval, $\rho=0$ wherever $u'=0$.
Transport and energy normalization give
\[
 \int_{\R^n}\frac{T^{p'}}{K_n}u^q\,dx=1,
 \qquad
 \int_{\R^n}\frac{T^{p'}}{K_n}\rho^pu^q\,dx=1.
\]
By the arithmetic--geometric mean inequality and
\eqref{eq:radial-det},
\[
 \operatorname{div}\mathcal T=M(\sigma+n-1)\ge nM\sigma^{1/n},
 \qquad
 \int u^{q_\#}\operatorname{div}\mathcal T\,dx\ge nV_\#.
\]
Together with \eqref{eq:radial-IBP}, this gives
\[
 \int_{\R^n}\frac{T^{p'}}{K_n}\rho u^q\,dx
 =\frac1{\Lambda K_n}\int T(-u')u^{q_\#-1}\,dx
 \geq\frac{nV_\#}{q_\#\Lambda K_n}=\varkappa_D^{-1}.
\]
Hence
\begin{equation}\label{eq:radial-Holder}
 \int_{\R^n}\frac{T^{p'}}{K_n}\mathcal H_p(\rho)u^q\,dx
 \leq p\bigl[1-(1-D)^{1/p}\bigr]\leq2D.
\end{equation}

Regard $\sigma,M$, and $\rho$ as functions of $x=\log T$
through the inverse transport. Set $t=T^{p'}/(1+T^{p'})$,
$Y=u^{-1/a_n}$ and $P=p't$. From
\cref{eq:radial-Jacobian},
\[
 Y=C_v(1+T^{p'})\frac rT\sigma^{-1/n},\qquad
 \frac{dr}{dT}=\frac r{T\sigma},\qquad
 \frac{dY}{dr}=p'C_v\varkappa_DT^{p'-1}\rho.
\]
Differentiating the first expression with respect to $r$ gives
\[
 \frac{dY}{dr}=\frac{C_v(1+T^{p'})}{T}\sigma^{-1/n}
      \left[1+(P-1)\sigma-\frac Tn\frac{d\sigma}{dT}\right].
\]
Comparing these expressions, and using $dT/dx=T$, gives
\begin{equation}\label{eq:transport-ODE}
 \sigma_x=n\bigl[1+(P-1)\sigma-P\varkappa_D\rho\sigma^{1/n}\bigr],
 \qquad (\log M)_x=1-1/\sigma.
\end{equation}

To control $\rho-1$ in
\eqref{eq:transport-ODE}, we use the
following bound for $\mathcal H_p$ defined in
\eqref{eq:transport-defects}, valid for every $z\geq0$:
\begin{equation}\label{eq:radial-scalar-H}
 \abs*{z-1}\leq A_p^H\sqrt{\mathcal H_p(z)}+B_p^H\mathcal H_p(z),
\end{equation}
where one may take
\[
 A_p^H=\left(\frac{2}{p(p-1)2^{-\abs*{p-2}}}\right)^{1/2},\qquad
 B_p^H=\frac1{\min\{2\mathcal H_p(1/2),\mathcal H_p(2)\}}.
\]
On $[1/2,2]$, we use $\mathcal H_p(1)=\mathcal H_p'(1)=0$ and integrate
the bound $\mathcal H_p''\geq p(p-1)2^{-\abs*{p-2}}$.
On $[0,1/2]$ and $[2,\infty)$, we use
$\mathcal H_p(z)/(1-z)\geq2\mathcal H_p(1/2)$ and
$\mathcal H_p(z)/(z-1)\geq \mathcal H_p(2)$, respectively.
This proves \eqref{eq:radial-scalar-H} for every $p>1$.

For any $k>0$, set $w_k(r)=k^{a_n}v_0(kr)$ and
$\theta(r)=\log(k/M(r))$. Since $kr=e^\theta T$, the Jacobian identity gives
\begin{equation}\label{eq:radial-reconstruction}
 \log\frac{w_k}{u}
 =a_n\left[\theta-
       \log\frac{1+e^{p'\theta}T^{p'}}{1+T^{p'}}\right]
                       -\frac1q\log\sigma.
\end{equation}
The preceding identities hold for every integer $n>p>1$.
We impose further restrictions on $n$ only when we need
estimates uniform in the dimension.

\subsubsection{A radial window estimate}

We follow \cite[proof of Theorem~2, Steps~II--IV]{CFMP}:
we choose an interval carrying most of the $L^q$ mass
and show that the transport map is close to a dilation there.
The bounds \eqref{eq:transport-AMGM} and
\eqref{eq:radial-Holder} allow us to keep track of
constants in both fixed and large dimensions.

\begin{lemma}\label{lem:normalized-tails}
Let $q>1$ and $m,\delta\geq0$. Let
$u,w\in L^q(\R^n)$ be nonnegative functions with
$\|u\|_q=\|w\|_q=1$.
Let $C\subseteq\R^n$ be measurable and suppose that
$u>0$ almost everywhere on $C$. If
\[
 \int_{C^c}u^q\leq2m,\qquad
 |\log(w/u)|\leq\delta\quad\text{almost everywhere on }C,
\]
then
\begin{equation}\label{eq:normalized-tail-bound}
 \int_{C^c}w^q\leq2m+q\delta,\qquad
 \|u-w\|_q^q\leq(e^\delta-1)^q
                  +2^{q-1}(4m+q\delta).
\end{equation}
\end{lemma}
\begin{proof}
Since $w^q\geq e^{-q\delta}u^q$ on $C$, normalization gives
$\int_{C^c}w^q\leq1-e^{-q\delta}(1-2m)\leq2m+q\delta$.
Integrating the $q$th power of $|u-w|\leq(e^\delta-1)u$ on $C$
and the bound $|u-w|^q\leq2^{q-1}(u^q+w^q)$ on $C^c$ gives
the second assertion.
\end{proof}

Let $\pi$ be the pushforward of $u^q(\xi)\,d\xi$
under $\xi\mapsto\log T(|\xi|)$.
By \eqref{eq:radial-CDF}, it is also the distribution
of $\log|\xi|$ under $v_0^q(\xi)\,d\xi$.
It is distinct from the gradient measure $\nu$
in \cref{sec:preliminaries}.
Set $K=K_n$ as in \eqref{eq:radial-K},
and let $D=\mathcal D(u)$.
Define
\[
 A_{\mathrm{win}}=n(1+K)\frac{2^{1+1/p}}p.
\]
Since $1+T^{p'}\geq1$, \eqref{eq:transport-AMGM}
and \eqref{eq:radial-Holder} give
\begin{equation}\label{eq:window-defects}
 \int h_n(\sigma)\,d\pi\leq A_{\mathrm{win}}D,\qquad
 \int \frac{T^{p'}}K \mathcal H_p(\rho)\,d\pi\leq2D
 \quad(D\leq1/2).
\end{equation}

\begin{lemma}\label{lem:one-window}
Let $u$ and its radial transport satisfy the assumptions of
\cref{sec:transport}, and let $\pi,K,D,A_{\mathrm{win}}$
be as defined above. Fix $\lambda\geq1$ and $x_*\in\R$, and set
$y=\lambda(x-x_*)$. Writing $d\pi=g_x(x)\,dx$, define
\[
 \pi_y=(x\mapsto\lambda(x-x_*))_\#\pi,\qquad
 d\pi_y=g(y)\,dy,\qquad
 g(y)=\lambda^{-1}g_x(x_*+y/\lambda).
\]
We regard $T,\sigma,M,\rho,P,t$ as functions of $y$ through
$x=x_*+y/\lambda$. Let $m\geq0$, $\ell>0$,
$g_*,w_*>0$ and $B\geq0$. Suppose that an interval
$I$ in the $y$ coordinate satisfies
\[
 \pi_y(I)\geq1-2m>\tfrac12,\quad |I|\leq\ell,\quad
 g\geq g_*>0,\quad T^{p'}/K\geq w_*>0,\quad |P-1|\leq B
 \quad\text{on }I.
\]
Let $A_p^H,B_p^H$ be the constants in
\eqref{eq:radial-scalar-H}. Set $a=\lambda^2D$ and define
\begin{align*}
 Q&=2/(g_*w_*),\qquad R_\lambda=n/\lambda^2,\\
 J&=A_p^H\sqrt{\ell Q}+B_p^HQ,\qquad U=nB/\lambda+2p'/\lambda,\\
 S&=e^{U\ell}\left(\sqrt{32A_{\mathrm{win}}}/\lambda
                          +4p'R_\lambda J+8p'R_\lambda\ell\right),\\
 W&=2\ell S,\\
 L&=a_n\left(BW/\lambda+p'^2W^2/(8\lambda^2)\right)+2S/q.
\end{align*}
If
\begin{equation}\label{eq:one-window-smallness}
 a\leq\min\{1,\lambda^2/2,\lambda^2/(512A_{\mathrm{win}}),(8S)^{-2}\},
\end{equation}
then there exists $k>0$ such that the normalized centered
bubble $w_k$ defined above satisfies
\begin{equation}\label{eq:one-window-estimate}
 \|u-w_k\|_q^q\leq(e^{L\sqrt a}-1)^q
                     +2^{q-1}(4m+qL\sqrt a).
\end{equation}
\end{lemma}
\begin{proof}
Using \eqref{eq:window-defects} and the lower bounds
$g\ge g_*$ and $T^{p'}/K\ge w_*$ on $I$, we obtain
$\int_I \mathcal H_p(\rho)\,dy\le QD$.
Hence, by \eqref{eq:radial-scalar-H} and the
Cauchy--Schwarz inequality,
$\int_I|\rho-1|\,dy\le J\sqrt D$ for $D\le1$.
By \eqref{eq:window-defects} and $\pi_y(I)>1/2$, we may choose
$y_0\in I$ with $h_n(\sigma(y_0))\leq2A_{\mathrm{win}}D$.
Since
\[
 h_n'(s)=\frac{n-1}{n}s^{-1-1/n}(s-1),
\]
we have
\[
 h_n(s)\geq
 \begin{cases}
  (s-1)^2/16,&1/2\leq s\leq2,\\
  1/64,&s\notin[1/2,2],
 \end{cases}
 \qquad s>0.
\]
Thus
$|\sigma(y_0)-1|\leq\sqrt{32A_{\mathrm{win}}}/\lambda\sqrt a\leq1/4$.

As long as $\sigma\in[1/2,2]$, \eqref{eq:transport-ODE} gives
\[
 |\sigma_y|\leq U|\sigma-1|
       +\frac{4np'}\lambda|\rho-1|+\frac{8np'}\lambda D.
\]
Here we used $|1-\sigma^{1/n}|\leq2|1-\sigma|/n$,
$\varkappa_D\leq2$, and $\varkappa_D-1\leq4D$.
We apply Gronwall's inequality to the left and right
of $y_0$, as long as $\sigma\in[1/2,2]$.
This gives $|\sigma-1|\le S\sqrt a$.
Since $S\sqrt a\le1/8$, continuity shows that
$\sigma$ stays in $[1/2,2]$ throughout $I$.

Choose $k=M(y_0)$. The second identity in \eqref{eq:transport-ODE} gives
\[
 |\theta|\leq W\sqrt a/\lambda\quad\hbox{on }I.
\]
For fixed $T$, the expression in brackets in
\eqref{eq:radial-reconstruction} has derivative
$1-P$ at $\theta=0$ and second derivative bounded
in absolute value by $p'^2/4$.
Using $|P-1|\le B$ and the bounds for $\theta$
and $\sigma$, we obtain
$|\log(w_k/u)|\le L\sqrt a$ on the annulus
corresponding to $I$.
This annulus carries at least $1-2m$ of the
$u^q$ mass, so \cref{lem:normalized-tails}
gives \eqref{eq:one-window-estimate}.
\end{proof}

\begin{proposition}\label{prop:both-radial}
Let $p>1$, let $n>p$ be an integer, and let $0<\beta<1$.
There exists $a_{\rm rad}(n,p,\beta)>0$ such that every
nonnegative, radially nonincreasing function
$u\in\dot W^{1,p}(\R^n)$ with $\|u\|_q=1$ satisfies
\[
 \mathcal D(u)\leq a_{\rm rad}(n,p,\beta)
 \quad\Longrightarrow\quad
 \inf_{\substack{V>0\text{ centered bubble}\\\|V\|_q=1}}
 \|u-V\|_q<\beta.
\]
Moreover, for every fixed $p>1$ and $0<\beta<1$,
there exist $N_R(p,\beta)>p$ and $a_R(p,\beta)>0$
such that the same conclusion holds for every integer
$n\geq N_R(p,\beta)$ and every such $u$ whenever
$n\mathcal D(u)\leq a_R(p,\beta)$.
\end{proposition}
\begin{proof}
\emph{Fixed dimension.}
Let $\tau=(p-1)/p$, $\alpha=n\tau$, $b=n/p$, and
$B_0=\mathrm B(\alpha,b)$. Given $0<\beta<1$, set
\[
 m=\frac{\beta^q}{2^{q+4}},\qquad
 h=\min\left\{\frac14,(m\alpha B_0)^{1/\alpha},
  (mbB_0/2^{(1-\alpha)_+})^{1/b}\right\}.
\]
Take $\lambda=1$ and $x_*=0$, so that $y=x$.
Let $I$ be the interval in the $y$ coordinate
corresponding to $h\leq t\leq1-h$. Since $b>1$, the two Beta tails satisfy
\[
 \pi_y(t<h)\le\frac{h^\alpha}{\alpha B_0},\qquad
 \pi_y(t>1-h)\le\frac{2^{(1-\alpha)_+}h^b}{bB_0}.
\]
Both tails are at most $m$, also when $\alpha<1$.
In particular, $\pi_y(I)\geq1-2m$.

The density is
\[
 g(x)=\frac{t^\alpha(1-t)^b}{\tau B_0}.
\]
The hypotheses of \cref{lem:one-window} hold with
$\ell=2\tau\log((1-h)/h)$, $g_*=h^n/(\tau B_0)$,
$w_*=h/K$, and $B=1+p'$.

All the constants in \cref{lem:one-window} are now fixed.
Choose a positive threshold for $D$ such that,
whenever $a=D$ is below this threshold,
\eqref{eq:one-window-smallness} holds and
\[
 L\sqrt a\leq
 \min\left\{1,\frac\beta{4e},\frac{\beta^q}{q2^{q+2}}\right\}.
\]
The right-hand side of \eqref{eq:one-window-estimate}
is then smaller than $\beta^q/2$.

\emph{Large dimensions.}
Write $d_* =\tau(1-\tau)$, $k_* =\min\{\tau,1-\tau\}$,
and $m=\beta^{2p}/2^{2p+4}$. Set
\[
 \lambda=\sqrt n,\quad
 y=\sqrt n\bigl(x-\tau\log(p-1)\bigr),\qquad
 |t-\tau|\leq\sqrt{d_*/(mn)}.
\]
The last inequality defines $I$.
Assume $n\geq2p$ and $n\geq4d_*/(mk_*^2)$.
Under $\pi_y$, the variable $t$ has mean $\tau$
and variance $d_* /(n+1)$.
Chebyshev's inequality therefore gives $\pi_y(I)\geq1-m$.
Since $t(1-t)\geq d_*/4$ on $I$,
\[
 |I|\leq\ell:=8p\sqrt{d_*/m},\qquad
 |P-1|\leq\frac{p'\sqrt{d_*/m}}{\sqrt n},\qquad
 T^{p'}/K\geq\frac1{4p}.
\]

The density has its maximum at $t=\tau$, and
\[
 \mathsf K(t)=\tau\log\frac\tau t
 +(1-\tau)\log\frac{1-\tau}{1-t},\qquad
 \mathsf K'(t)=\frac{t-\tau}{t(1-t)},\qquad
 \frac{g(y)}{g(0)}=e^{-n\mathsf K(t)}.
\]
By our choice of $n$, $|t-\tau|\le k_*/2$ on $I$.
Thus $t\ge\tau/2$ and $1-t\ge(1-\tau)/2$, so
\[
 \mathsf K''(t)=\frac{\tau}{t^2}
  +\frac{1-\tau}{(1-t)^2}\le\frac4{d_*}
 \qquad(t\in I).
\]
Since $\mathsf K(\tau)=\mathsf K'(\tau)=0$,
Taylor's theorem gives $n\mathsf K(t)\le2/m$ on $I$.
Moreover, $g(y)\le g(0)$, and hence
$g(0)|I|\ge\pi_y(I)\ge1-m$. Therefore,
\[
 g(y)\geq g_*:=\frac{1-m}{\ell}e^{-2/m}.
\]

In Lemma~\ref{lem:one-window}, $A_{\mathrm{win}}/n\leq8$,
$R_\lambda=1$,
$nB/\lambda\leq p'\sqrt{d_*/m}$,
$a_nB/\lambda\leq p'\sqrt{d_*/m}/p$, and
$a_n/\lambda^2\leq1/p$.
Thus the constants $S$ and $L$ in
\cref{lem:one-window} have explicit upper bounds
depending only on $p$ and $m$.
Since $p<q\le2p$, we can choose a positive
threshold for $a=nD$, depending only on
$p$ and $\beta$, such that whenever $a$ is
below this threshold,
\eqref{eq:one-window-smallness} holds and
\[
 L\sqrt a\leq
 \min\left\{1,\frac\beta{4e},
               \frac{\beta^{2p}}{(2p)2^{2p+2}}\right\}.
\]
The right-hand side of \eqref{eq:one-window-estimate}
is again smaller than $\beta^q/2$.

\medskip
\phantomsection
\noindent\emph{Approximation by smooth radial functions.}\label{sec:radialapproximation}
For a nonnegative radial nonincreasing function
$u$ with $\|u\|_q=1$, set
\[
 f_{M,\varepsilon}=\min\{(u-\varepsilon)_+,M\}.
\]
Since $u\in L^q$, its $\varepsilon$ superlevel ball has finite radius;
thus $f_{M,\varepsilon}$ is bounded and compactly supported. Dominated
convergence gives $f_{M,\varepsilon}\to u$ in $L^q$ as
$\varepsilon\downarrow0$, $M\uparrow\infty$, while the chain rule gives
\[
 \nabla f_{M,\varepsilon}
    =\mathbf1_{\{\varepsilon<u<\varepsilon+M\}}\nabla u
      \quad\hbox{a.e.}
\]
Since $\nabla u=0$ almost everywhere on $\{u=0\}$, the gradients
also converge strongly in $L^p$.

We convolve $f_{M,\varepsilon}$ with a smooth,
nonnegative, radially nonincreasing kernel that is
positive in the interior of its support ball.
By the layer-cake formula, the convolution is radial
and nonincreasing, since the volume of the intersection
of two balls decreases with the distance between their
centers. The convolution is also positive in the interior
of its support: there the positivity sets of the two
factors overlap in positive measure. For fixed $n$,
mollification converges in $L^q$ and $\dot W^{1,p}$.
After normalization in $L^q$, we obtain
smooth radial nonincreasing functions $u_j$,
positive in the interiors of their supports,
such that
\[
 u_j\longrightarrow u\quad\hbox{in }L^q,
 \qquad \nabla u_j\longrightarrow\nabla u\quad\hbox{in }L^p,
 \qquad \mathcal D(u_j)\longrightarrow\mathcal D(u).
\]

The preceding estimates apply to smooth functions
as in Subsection~\ref{sec:transport}.
We now halve the positive deficit thresholds chosen
above for smooth functions. If the deficit of $u$
is at most the resulting threshold, the convergence
$\mathcal D(u_j)\to\mathcal D(u)$ allows us to apply
the smooth estimate to $u_j$ for all sufficiently
large $j$.
Since the distance to the set of normalized centered
bubbles is $1$-Lipschitz in $L^q$, passing to the
limit gives a distance at most $2^{-1/q}\beta<\beta$.
The halved thresholds have the stated dependence
on $n,p,\beta$ and on $p,\beta$, respectively.
\end{proof}

\subsection{From local to global stability}\label{sec:global}

Following the local-to-global argument of
\cite[Section~3.1]{DEFLL}, we use finitely many
polarizations and continuous paths joining each
function to its polarization. These paths preserve
the $L^q$ norm and gradient $p$-energy.
Together with \cref{prop:both-radial}, this allows us
to pass from local to global stability.

For $e\in\mathbb S^{n-1}$ and $a\geq0$, let
$H_{a,e}=\{x:x\cdot e\leq a\}$ and
$\sigma_{a,e}x=x+2(a-x\cdot e)e$. For a function $g$, define
\begin{equation}\label{eq:polarization}
 P_{a,e}g(x)=
 \begin{cases}
 \max\{g(x),g(\sigma_{a,e}x)\},&x\in H_{a,e},\\
 \min\{g(x),g(\sigma_{a,e}x)\},&x\notin H_{a,e}.
 \end{cases}
\end{equation}
\begin{lemma}\label{lem:continuous-polarization}
Let $e\in\mathbb S^{n-1}$ and
$g\in W^{1,p}(\R^n)\cap L^q(\R^n)$.
Then the map $a\mapsto P_{a,e}g$ is continuous from
$[0,\infty)$ into both $W^{1,p}(\R^n)$ and $L^q(\R^n)$.
If, in addition, $g\geq0$ and $\operatorname{supp}g\subseteq B_R$
for some $R>0$, then the path $0\leq a\leq R$ preserves
all $L^s$ norms for which they are finite, the gradient
$p$-energy, and support in $B_R$. Moreover, $P_{R,e}g=g$.
\end{lemma}
\begin{proof}
For a fixed reflection, the maximum and minimum
of two Sobolev functions depend continuously on
the functions in $W^{1,p}$. Indeed, by the Sobolev
chain rule, the gradients of the two functions agree
almost everywhere on their coincidence set.
Under strong $W^{1,p}$ convergence, the gradients
of the maxima and minima converge in measure.
Their $p$-th powers are uniformly integrable,
and their integrals outside sufficiently
large balls are uniformly small.
Together with convergence in measure, this
gives strong convergence of the gradients
in $L^p$. The pointwise inequalities
for maxima and minima also give continuity in
$L^p$ and $L^q$. Applying these facts on the two
half-spaces proves continuity of $P_{0,e}$; the
two branches have the same trace on the plane.

For $y\in\mathbb R^n$, write $\mathsf T_y g(x)=g(x-y)$. Then
\[
 P_{a,e}=\mathsf T_{ae}P_{0,e}\mathsf T_{-ae}.
\]
Translations are strongly continuous in $W^{1,p}$ and $L^q$.
Together with the continuity of $P_{0,e}$ proved above, this
shows that $a\mapsto P_{a,e}g$ is continuous in both spaces.

For each pair $x,\sigma_{a,e}x$, polarization either keeps the
two values or exchanges them. When values are exchanged, their gradients are exchanged
up to the orthogonal reflection. Thus the
$L^s$ norms and the gradient $p$-energy are unchanged.

Now suppose that $g\geq0$ and $\operatorname{supp}g\subseteq B_R$.
For $x\in H_{a,e}$ and $a\geq0$, we have
$|\sigma_{a,e}x|\geq|x|$, while on the other half-space
$0\leq P_{a,e}g\leq g$. Hence $P_{a,e}g$ is also supported
in $B_R$. At $a=R$, the reflecting plane is tangent to $B_R$,
so $P_{R,e}g=g$.
\end{proof}

\begin{proposition}\label{prop:exact-crossing}
Let $p>1$, let $n>p$ be an integer, and fix
$\varepsilon,b_0>0$. Suppose that every nonzero nonnegative
$u\in\dot W^{1,p}(\R^n)$ satisfies
\[
 d(u)\leq\varepsilon\quad\Longrightarrow\quad
 \mathcal D(u)\geq b_0d(u)^{r_p}.
\]
Set
\[
 \alpha_p=\max\{p,p'\},\quad
 K_p=2[\alpha_p\max\{1/p,1/2\}]^{1/\alpha_p},\quad
 \eta=\tfrac12\min\{\varepsilon,1\}.
\]
Let $\beta=\tfrac13(\eta/(2K_p))^{\alpha_p}\in(0,1)$,
and let $A_R=A_R(n,p,\beta)>0$ be a fixed-dimensional
radial approximation threshold from \cref{prop:both-radial}.
Choose
\begin{equation}\label{eq:short-global-choice}
 0<a_0\leq\min\{A_R,\beta,\tfrac12b_0\eta^{r_p}\},\qquad
 c_+=\min\{a_0,b_0\}.
\end{equation}
Then every nonzero nonnegative $u\in\dot W^{1,p}(\R^n)$
satisfies $\mathcal D(u)\geq c_+d(u)^{r_p}$.
For every nonzero $u\in\dot W^{1,p}(\R^n)$, one has
\begin{equation}\label{eq:short-real-global}
 \mathcal D(u)\geq
 2^{-r_p}\min\{c_+,2(1-2^{-p/n})\}\,d(u)^{r_p}.
\end{equation}
\end{proposition}
\begin{proof}
Let $f\geq0$ be smooth and compactly supported, with $\|f\|_q=1$
and $D=\mathcal D(f)<a_0$. Let $f^*$ denote the symmetric decreasing rearrangement of $f$.
By the P\'olya--Szeg\H{o} inequality,
$\|\nabla f^*\|_p\leq\|\nabla f\|_p$ and
$\|f^*\|_q=\|f\|_q$. Thus $\mathcal D(f^*)\leq D$,
so $f^*$ is within $\beta$ in $L^q$ of a normalized centered
bubble $V$. 

Choose $R_f,M_f>0$ such that
$\operatorname{supp}f\subset B_{R_f}$ and
$0\le f\le M_f$. We apply
\cite[Proposition~5.2 and (5.1)]{BF}.
In that paper, the reflection parameter is the
distance $2a=|\sigma_{a,e}(0)|$, so a uniform
choice from $(0,2R_f)$ corresponds to $a$
uniformly distributed on $(0,R_f)$ in
\eqref{eq:polarization}. Both conventions place
the larger value on the side containing the origin.
Let $P_Nf$ be the result of $N$ successive
independent polarizations with parameters chosen
uniformly from $(0,R_f)\times\mathbb S^{n-1}$.
The cited estimate gives
\[
 \mathbb E\|P_Nf-f^*\|_1
        \leq\frac{2n|B_{2R_f}|M_f}{N}.
\]
Since $0\le P_Nf,f^*\le M_f$, we have
$\|P_Nf-f^*\|_q^q\le M_f^{q-1}\|P_Nf-f^*\|_1$.
For sufficiently large finite $N$, the expectation
of the left-hand side is less than $\beta^q$.
Thus some choice of the polarizations has an
endpoint $g$ satisfying $\|g-f^*\|_q<\beta$.

Each polarization in the construction preserves support in
$B_{R_f}$, the $L^q$ norm, and the gradient $p$-energy.
For each step, let $h$ be the current function
and $(a,e)$ the chosen polarization parameters.
Let the parameter in $P_{t,e}h$ decrease
continuously from $R_f$ to $a$. By \cref{lem:continuous-polarization}, this gives a path
from $h=P_{R_f,e}h$ to $P_{a,e}h$, continuous in $W^{1,p}$
and $L^q$ and preserving the $L^q$ norm and gradient
$p$-energy. Concatenating the finitely many paths, we obtain a
path from $f$ to $g$ along which these quantities are constant.

We next estimate the gradient distance between
$g$ and $V$.
Let $\|f_1\|_q=\|f_2\|_q=1$, suppose that
$E(f_1),E(f_2)\le E$, and put
$z=\|f_1-f_2\|_q\le1$.
Since $\|(f_1+f_2)/2\|_q\ge1-z/2$,
the Sobolev inequality \eqref{eq:sobolev} and
Clarkson's inequalities for vector-valued functions
\cite[(2.1)--(2.2)]{Neumayer} give
\begin{equation}\label{eq:polarization-gradient-distance}
\frac{\|\nabla(f_1-f_2)\|_p}{E^{1/p}}
\leq2\left[
1-(1-D)^{\alpha_p/p}(1-z/2)^{\alpha_p}
\right]^{1/\alpha_p}
\leq K_p(D+z)^{1/\alpha_p},
\end{equation}
where $D=1-S_{n,p}^p/E$. The last inequality
follows from Bernoulli's inequality.
For the endpoint $g$ and the centered bubble $V$
chosen above, we have $\|g-V\|_q<2\beta<1$.
Since $D=\mathcal D(f)<a_0\leq\beta$,
applying \eqref{eq:polarization-gradient-distance}
with $f_1=g$, $f_2=V$, and $E=E(f)$ gives
$d(g)<K_p(3\beta)^{1/\alpha_p}=\eta/2$.

The path preserves both $E$ and the $L^q$ norm,
so $\mathcal D$ equals $D$ throughout the path.
The path is continuous in gradient $L^p$ and has
constant positive energy, so $d$ is continuous
along it. If $d(f)>\eta$, then $d(g)<\eta/2$
implies that some $f_s$ on the path satisfies
$d(f_s)=\eta$. The assumed local estimate gives
$D=\mathcal D(f_s)\ge b_0\eta^{r_p}>a_0$,
a contradiction. Thus $d(f)\le\eta$.

\medskip
\phantomsection
\noindent\emph{Approximation by smooth nonnegative functions.}\label{sec:approximation}
Fix $n,p$, and let $u\ge0$ belong to $\dot W^{1,p}(\R^n)$.
Choose $\chi\in C_c^\infty(\R^n)$ with $0\le\chi\le1$,
$\chi=1$ on $B_1$, and $\chi=0$ outside $B_2$, and set
$\chi_R(x)=\chi(x/R)$. Then $u\chi_R\to u$ in
$L^q$ and $(\chi_R-1)\nabla u\to0$ in $L^p$. Since
$1/p=1/q+1/n$,
\[
 \norm*{u\nabla\chi_R}_p
 \leq\norm*{u}_{L^q(B_{2R}\setminus B_R)}\norm*{\nabla\chi_R}_n
 \leq C_n\norm*{u}_{L^q(\R^n\setminus B_R)}\longrightarrow0.
\]
For each $R$, convolve $u\chi_R$ with a nonnegative smooth,
compactly supported mollifier. Choosing the mollification scale
sufficiently small as $R\to\infty$ gives nonnegative functions
in $C_c^\infty(\R^n)$ converging to $u$ in $L^q$ and
$\dot W^{1,p}$.
For $u\ne0$, normalize the approximants in $L^q$ and denote them
by $u_j$. Their normalization factors tend to $\norm*{u}_q^{-1}$,
and both deficits and gradient distances converge
along the normalized approximating sequence.

If $\mathcal D(u)<a_0$, the convergence
$\mathcal D(u_j)\to\mathcal D(u)$ gives
$\mathcal D(u_j)<a_0$ for all sufficiently large
$j$. Thus $d(u_j)\le\eta$ by the smooth case,
and continuity of $d$ yields $d(u)\le\eta$.
We can now apply the assumed local estimate.
If $\mathcal D(u)\ge a_0$, we use $d(u)\le1$.
In either case,
$\mathcal D(u)\ge c_+d(u)^{r_p}$.

If $u$ has one sign, the result follows from
the nonnegative case, since
$\mathcal D(-u)=\mathcal D(u)$ and $d(-u)=d(u)$.
Otherwise, both $u_+$ and $u_-$ are nonzero.
We normalize $E(u)=1$ and, replacing $u$ by
$-u$ if necessary, assume that
$w=E(u_+)\ge1/2$ and put $z=E(u_-)=1-w$.
The deficit and distance satisfy
\[
 \mathcal D(u)\geq\tfrac12w\mathcal D(u_+)
                 +(1-2^{-p/n})z,\qquad
 d(u)^{r_p}\leq2^{r_p-1}\{w d(u_+)^{r_p}+z\}.
\]
The first inequality follows by averaging
\[
\mathcal D(u)\ge
w\mathcal D(u_+)+z\mathcal D(u_-)
\]
and
\[
\mathcal D(u)\ge
2(1-2^{-p/n})z.
\]
The latter bound follows from the concavity of
\[
s\mapsto
1-\bigl((1-s)^{q/p}+s^{q/p}\bigr)^{p/q}
\]
on $[0,1/2]$ and its endpoint values
$0$ and $1-2^{-p/n}$.

For the distance estimate, the triangle
inequality gives
\[
d(u)\le w^{1/p}d(u_+)+z^{1/p}.
\]
Raising this inequality to the power $r_p$
and using $r_p\ge p$ gives the second
inequality above.
Combining the two estimates above with
$\mathcal D(u_+)\ge c_+d(u_+)^{r_p}$
gives \eqref{eq:short-real-global}.
\end{proof}

\subsection{Proof of the lower bound}\label{sec:assembly}

\begin{proof}
Let $\varepsilon_p,b_p,N_p^L$ be the constants from
\cref{prop:local}. Let $K_p,\alpha_p$ be as in
\cref{prop:exact-crossing}. Set
$\eta_p^{\rm gl}=\tfrac12\min\{\varepsilon_p,1\}$ and
$\beta_p=\tfrac13(\eta_p^{\rm gl}/(2K_p))^{\alpha_p}$.
With $a_R(p,\beta_p)$ from \cref{prop:both-radial}, set
\[
 a_* =\min\{a_R(p,\beta_p),\beta_p,\tfrac12b_p(\eta_p^{\rm gl})^{r_p}\},
 \qquad m_p^{\rm gl}=\min\{a_*,b_p\}.
\]
Choose an integer
$N(p)>\max\{N_p^L,N_R(p,\beta_p),2p\}$.
For $n\ge N(p)$, apply \cref{prop:exact-crossing} with
$b_0=b_p/n$, $A_R=a_R(p,\beta_p)/n$, and
$a_0=a_*/n$.
By the definition of $a_*$, these choices satisfy
\eqref{eq:short-global-choice}, and the coefficient
for nonnegative functions is
$c_+=m_p^{\rm gl}/n$.
Since $n\ge2p$, the inequality
$1-e^{-x}\ge x/2$ for $0\le x\le (\log2)/2$ gives
$2(1-2^{-p/n})\ge p\log2/n$.
The estimate for arbitrary functions in
\cref{prop:exact-crossing} therefore gives
\begin{equation}\label{eq:large-dimensional-global-bound}
\mathcal D(u)\geq\frac{c_\infty(p)}n d(u)^{r_p},\qquad
c_\infty(p)=2^{-r_p}\min\{m_p^{\rm gl},p\log2\}>0
\quad(n\geq N(p)).
\end{equation}

For each integer $n$ with $p<n<N(p)$, use
\cref{thm:super-local} if $p\ge2$ and
\cref{thm:short-subquadratic-local} if $1<p<2$.
In \cref{prop:iv-invariant}, take
$a=cE_v^{r_p/p-1}$ and $\rho=t_*/E_v^{1/p}$,
where $c$ and $t_*$ are the coefficient and gradient
radius in the corresponding orthogonal estimate
($t_*=\delta_0$ when $1<p<2$).
This gives a local radius $\varepsilon_{n,p}>0$
and coefficient $b_{n,p}>0$.
The fixed-dimensional part of \cref{prop:both-radial},
together with \cref{prop:exact-crossing}, then yields
$\mathcal D(u)\ge c(n,p)d(u)^{r_p}$ for every
nonzero $u\in\dot W^{1,p}(\R^n)$, with $c(n,p)>0$.

Combining \eqref{eq:large-dimensional-global-bound}
with the bounds obtained for the finitely many
remaining dimensions, we choose
\begin{equation}\label{eq:finalcp}
 c_p=\min\{c_\infty(p),\ n c(n,p):
                   n\in\mathbb N,\ p<n<N(p)\}>0,
\end{equation}
omitting the finite sublist when it is empty.
This proves the lower bound in
\cref{eq:mainconstants} for every integer $n>p$.

All constants entering this finite minimum are determined by the
preceding estimates.
\end{proof}

\subsection{Proof of the upper bound}\label{sec:upper}

Set $\widehat v=E_v^{-1/p}v$, so that
$\|\nabla\widehat v\|_p=1$.
We prove the upper bound by considering the sum
of two widely separated translates of $\widehat v$.
We need to estimate its distance to the full
manifold \eqref{eq:manifold}, allowing arbitrary
amplitudes, scales, and centers.

\begin{proposition}\label{prop:upper}
Let $p>1$ and let $n>p$ be an integer.
Let $e_1$ be the first coordinate unit vector, and, for $R>0$, set
$u_R=\widehat v+\widehat v(\cdot-Re_1)$.
Then
\begin{equation}\label{eq:twobubbledistance}
 \operatorname{dist}_{\nabla,p}(u_R)\longrightarrow1,\qquad
 d(u_R)^p\longrightarrow\tfrac12
 \quad\text{as }R\longrightarrow\infty.
\end{equation}
Consequently,
\[
 \mathsf C_{n,p}\leq2^{r_p/p}(1-2^{-p/n})
                \leq\frac{2^{r_p/p}p\log2}{n}.
\]
\end{proposition}

\begin{proof}
For any $L>0$, the restrictions of $\widehat v$
and $\nabla\widehat v$ to $B_L$ have disjoint
translates when $R>2L$. The errors in $L^q$
and $L^p$ are bounded by twice the corresponding
tail norms. Letting first $R\to\infty$ and then
$L\to\infty$, we obtain
\begin{equation}\label{eq:twobubblemass}
 E(u_R)\longrightarrow2,\qquad
 S_{n,p}^p\norm*{u_R}_q^p\longrightarrow2^{p/q},\qquad
 \mathcal D(u_R)\longrightarrow1-2^{-p/n}.
\end{equation}

Taking $V=\widehat v$ in \eqref{eq:distance} gives
$\operatorname{dist}_{\nabla,p}(u_R)\leq1$ for every $R$.
Choose $R_j\to\infty$ such that
$\operatorname{dist}_{\nabla,p}(u_{R_j})$
tends to its lower limit. For each $j$, choose
$W_j\in\mathcal M_{n,p}$ satisfying
\[
 W_j=c_j\lambda_j^{(n-p)/p}\widehat v(\lambda_j(x-y_j)),
 \qquad
 \norm*{\nabla(u_{R_j}-W_j)}_p
 \leq\operatorname{dist}_{\nabla,p}(u_{R_j})+j^{-1},
\]
where $c_j\in\R\setminus\{0\}$, $\lambda_j>0$,
and $y_j\in\R^n$ are unrestricted.
Since $\norm*{\nabla u_{R_j}}_p\leq2$ and
$\norm*{\nabla W_j}_p=\abs*{c_j}$,
the triangle inequality gives
$\abs*{c_j}\leq3+j^{-1}$.
After passing to a subsequence, we may assume
that $\lambda_j\to0$, $\lambda_j\to\lambda\in(0,\infty)$,
or $\lambda_j\to\infty$.

If $\lambda_j\to0$, take $x_j=0$. Boundedness of
$\nabla \widehat v$ gives
\[
 \norm*{\nabla W_j}_\infty
 \leq C_{n,p}\abs*{c_j}\lambda_j^{n/p}\longrightarrow0
\]
regardless of the centers. Moreover,
$\nabla u_{R_j}(x)\to\nabla \widehat v(x)$ for every fixed $x$.

If $\lambda_j\to\lambda\in(0,\infty)$, at least
one of the two centers $0$ and $R_je_1$ is at
distance at least $R_j/2$ from $y_j$.
After passing to a subsequence, choose that
same component for every $j$, and denote its
center by $x_j\in\{0,R_je_1\}$.
Then $\abs*{y_j-x_j}\to\infty$, so the decay
of $\nabla\widehat v$ gives
$\nabla W_j(x_j+x)\to0$ for every fixed $x$,
while $\nabla u_{R_j}(x_j+x)\to\nabla\widehat v(x)$.

Finally, suppose $\lambda_j\to\infty$ and take
$x_j=0$. After passing to a subsequence, either
$y_j\to y\in\R^n$ or $\abs*{y_j}\to\infty$.
The explicit bubble satisfies
\[
 \abs*{\nabla\widehat v(x)}
 \leq C_{n,p}\abs*{x}^{-\gamma}\quad(\abs*{x}\geq1),
 \qquad
 \gamma=\frac{n-1}{p-1}>\frac np.
\]
For every fixed $x\ne y$ in the first case,
and every fixed $x$ in the second,
$\abs*{x-y_j}$ is bounded below by a positive
number for all sufficiently large $j$.
Consequently, $\lambda_j\abs*{x-y_j}\geq1$
eventually, and
\[
 \abs*{\nabla W_j(x)}
 \leq C_{n,p}\abs*{c_j}\lambda_j^{n/p-\gamma}
       \abs*{x-y_j}^{-\gamma}
 \longrightarrow0\quad\text{for almost every }x.
\]
Also $\nabla u_{R_j}(x)\to\nabla\widehat v(x)$
for every fixed $x$.

In each case, we have chosen
$x_j\in\{0,R_je_1\}$ such that
$\nabla(u_{R_j}-W_j)(x_j+\cdot)$ converges
almost everywhere to $\nabla\widehat v$.
By translation invariance and Fatou's lemma,
\[
 \liminf_j\norm*{\nabla(u_{R_j}-W_j)}_p^p
 \geq\int_{\R^n}\abs*{\nabla\widehat v}^p\,dx=1.
\]
Since $W_j$ are approximate minimizers and
$\operatorname{dist}_{\nabla,p}(u_R)\leq1$,
we conclude that
$\operatorname{dist}_{\nabla,p}(u_R)\to1$.
Combining this convergence with
\eqref{eq:twobubblemass} proves
\eqref{eq:twobubbledistance}.
Consequently, by the definition \eqref{eq:bestconstant},
\[
 \mathsf C_{n,p}
 \leq\lim_{R\to\infty}\frac{\mathcal D(u_R)}{d(u_R)^{r_p}}
 =2^{r_p/p}(1-2^{-p/n})
 \leq\frac{2^{r_p/p}p\log2}{n}.
\]
Together with the lower bound in \cref{sec:assembly}, this
proves \cref{thm:main}.
\end{proof}

\appendix
\crefalias{section}{appendix}
\crefalias{subsection}{appendix}
\crefname{appendix}{Appendix}{Appendices}
\Crefname{appendix}{Appendix}{Appendices}

\section{Gradient and small-support estimates}\label{app:coordinate}

We prove gradient bounds for the inverses of the coordinate
operators introduced in Subsection~\ref{sec:basic-tools},
together with weak integration-by-parts identities and
small-support estimates used in \cref{lem:core,lem:capacity}.

\begin{lemma}\label{lem:coordinate}
Fix $p>1$ and an integer $n>p$ such that
$n\geq\max\{4,\kappa^{-1}\}$. Let $A,B_S$ and their closed forms
be as in \cref{sec:basic-tools}. Then the following hold.
\begin{enumerate}
\item[(i)] Suppose $g\in L^\infty(\mu_t)$ satisfies
$\int g\,d\mu_t=0$, and set $\psi=A^{-1}g$. Then
\begin{equation}\label{eq:coordinate-flux}
 \psi'(t)=-\frac{\int_0^t m(s)g(s)\,ds}{a(t)m(t)}.
\end{equation}
In particular,
\begin{equation}\label{eq:coordinate-flux-bound}
 \abs*{a(t)m(t)\psi'(t)}
 \leq\norm*{g}_\infty\min\{\mu_t((0,t)),\mu_t((t,1))\},
\end{equation}
and
\begin{equation}\label{eq:coordinate-radial-gradient}
 \norm*{\psi'}_\infty\leq\frac{R_p}{\sqrt n}\norm*{g}_\infty,
 \qquad R_p=\frac2\kappa+
                       \frac\pi{\sqrt{\kappa d_*}}.
\end{equation}
If $g\in L^\infty(\sigma)$ satisfies $\int g\,d\sigma=0$, then
\begin{equation}\label{eq:coordinate-sphere-gradient}
 \norm*{\nabla_S B_S^{-1}g}_\infty
       \leq\frac\pi{2\sqrt{n-2}}\norm*{g}_\infty.
\end{equation}
The bounds \eqref{eq:coordinate-radial-gradient} and
\eqref{eq:coordinate-sphere-gradient} also apply to measurable
families satisfying the same hypotheses, uniformly whenever
$\|g\|_\infty$ is uniformly bounded.

We also have $a(t)m(t)\psi'(t)\to0$ as
$t\to0^+$ and $t\to1^-$. For every bounded,
locally absolutely continuous function $h$
with $a\psi'h'\in L^1(\mu_t)$, we have
\begin{equation}\label{eq:coordinate-weak-pairings}
 \int_0^1 a\psi'h'\,d\mu_t=\int_0^1 gh\,d\mu_t.
\end{equation}
The analogous identity on the sphere holds for
$h\in W^{1,1}(\mathbb S^{n-1})$.

\item[(ii)] Let $\mathcal H$ be a Hilbert space,
regarded as a real Hilbert space if necessary.
Suppose that $G:(0,1)\to\mathcal H$ is bounded and locally
absolutely continuous, with
\[
 \int_0^1 a(t)\|G'(t)\|_{\mathcal H}^2\,d\mu_t(t)<\infty.
\]
Then
\begin{equation}\label{eq:coordinate-stein}
 n\int_0^1(t-\tau)\abs*{G}^2\,d\mu_t
 =2\int_0^1 a\langle G,G'\rangle\,d\mu_t.
\end{equation}
Every bounded member of the $\mathcal H$-valued closed radial
form domain has a representative satisfying these hypotheses.

\item[(iii)] Let $D=A$ with $\pi=\mu_t$, or $D=B_S$ with
$\pi=\sigma$. Suppose that $z\in L^2(\pi)$ satisfies
$\pi(z\ne0)\leq m$ for some $m\in(0,\mathrm e^{-1}]$,
and set $L=\log(1/m)$. Then
\begin{equation}\label{eq:coordinate-support-inverse}
 \left\langle z-\int z\,d\pi,
       D^{-1}\left(z-\int z\,d\pi\right)\right\rangle_{L^2(\pi)}
       \leq\frac{20}{nL}\norm*{z}_2^2.
\end{equation}
If $F$ belongs to the closed form domain associated with the
product Dirichlet energy in \eqref{eq:coordinate-product-lsi}
and $\mu(F\ne0)\leq m$ for some $m\in(0,1)$, then
\begin{equation}\label{eq:coordinate-support-energy}
 \int\big(a\abs*{F_t}^2+\abs*{\nabla_SF}^2\big)\,d\mu
       \geq\frac n8\log\frac1m\int F^2\,d\mu.
\end{equation}
\end{enumerate}
\end{lemma}

\begin{proof}
Let $D=A$ or $B_S$, set $P_s=e^{-sD}$, and write
$\Gamma=\Gamma_D$. By the curvature estimates derived from
\eqref{eq:coordinate-curvature} and
the sphere Bochner formula, we can apply
\cite[Theorem~4.7.2(ii),(iv) and (4.7.6),
pp.~209--210]{BGL}
to these closed operators. We obtain
\begin{equation}\label{eq:coordinate-gradient-contraction}
 \Gamma(P_s f)\leq e^{-2Ks}P_s\Gamma(f),
\end{equation}
and
\begin{equation}\label{eq:coordinate-reverse-poincare}
 \Gamma(P_s f)\leq
 \frac K{e^{2Ks}-1}\big(P_sf^2-(P_sf)^2\big)
 \leq\frac K{e^{2Ks}-1}\norm*{f}_\infty^2.
\end{equation}
For bounded mean-zero $g$, the spectral gaps in
\eqref{eq:coordinate-spectra} give
$D^{-1}g=\int_0^\infty P_sg\,ds$ in $L^2$.
The gradients are summable in $L^\infty$, since
\begin{equation}\label{eq:coordinate-gradient-integral}
 \int_0^\infty\sqrt{\frac K{e^{2Ks}-1}}\,ds
 =\frac\pi{2\sqrt K}.
\end{equation}
Combining \eqref{eq:coordinate-reverse-poincare}
with \eqref{eq:coordinate-gradient-integral}
and taking weak derivatives in the
resolvent representation gives
\eqref{eq:coordinate-sphere-gradient} and
\[
 \norm*{\sqrt a\,(A^{-1}g)'}_\infty
 \leq\frac\pi{2\sqrt{n\kappa}}\norm*{g}_\infty.
\]
Approximating in $L^2$ by uniformly bounded functions extends
both gradient estimates to bounded $g$.

For $\psi=A^{-1}g$, the weak equation gives
$(am\psi')'=-mg$ on compact subintervals.
Since $mg\in L^1(0,1)$, $am\psi'$ has
finite limits at both endpoints. A nonzero limit would make
\[
 \mathcal E_A(\psi,\psi)=\int_0^1\frac{\abs*{am\psi'}^2}{am}\,dt
\]
infinite, since $\alpha,\beta\geq1$. Both limits are zero.
Integration gives \cref{eq:coordinate-flux} and
\cref{eq:coordinate-flux-bound}. If $t\leq\tau/2$, then
\cref{eq:coordinate-density-derivative} gives
\[
 am(t)\geq\frac{n\tau}{2}\int_0^t m(s)\,ds,\qquad
 \abs*{\psi'(t)}\leq\frac2{n\tau}\norm*{g}_\infty.
\]
For $t\geq(1+\tau)/2$, the corresponding bound is
$\abs*{\psi'}\leq2\norm*{g}_\infty/[n(1-\tau)]$.
In the intervening interval, $a\geq d_*/2$, so
\[
 \abs*{\psi'}\leq\frac\pi{\sqrt{2n\kappa d_*}}\norm*{g}_\infty.
\]
Together, these bounds give \eqref{eq:coordinate-radial-gradient}
with the stated $R_p$.

For \cref{eq:coordinate-weak-pairings}, integrate on
$[\varepsilon,1-\varepsilon]$: the boundary term $(am\psi')h$ vanishes
by \cref{eq:coordinate-flux-bound} and boundedness of $h$, and the
integrals converge by the assumed absolute integrability.
On the sphere, the bounded gradient and distributional Laplacian $-g$
of $B_S^{-1}g$ give the pairing by smooth density in $W^{1,1}$.
For measurable families satisfying these hypotheses, Fubini's theorem
gives the same identities whenever the corresponding product integrands
are absolutely integrable. In \cref{sec:local},
absolute integrability follows from the bounded dual field and
$\mathsf Zh\in L^p(\nu)$.

A bounded member of the $\mathcal H$-valued closed radial form
domain has a locally absolutely continuous representative: on compact
subintervals of $(0,1)$, the form norm controls the $W^{1,2}$ norm.
Since $G$ is bounded and has finite weighted derivative energy,
$(\abs*{G}^2)'=2\langle G,G'\rangle$ is integrable against
$a\,d\mu_t$ by the Cauchy--Schwarz inequality. Boundedness also gives
$am\abs*{G}^2\to0$ at both endpoints. Integrating
\eqref{eq:coordinate-density-derivative} against $\abs*{G}^2$ on
compact subintervals and passing to the limit proves
\eqref{eq:coordinate-stein}.

We now prove (iii), starting with
\eqref{eq:coordinate-support-inverse}.
Assume that $\pi(z\ne0)\leq m\leq e^{-1}$
and set $L=\log(1/m)\geq1$.
H\"older's inequality gives
$\norm*{z}_r\leq m^{1/r-1/2}\norm*{z}_2$ for $1\leq r\leq2$.
Thus \cref{eq:coordinate-hypercontractivity} implies
\begin{equation}\label{eq:coordinate-supported-heat}
 \langle z,P_sz\rangle=\norm*{P_{s/2}z}_2^2
 \leq e^{-L\tanh(ns/8)}\norm*{z}_2^2.
\end{equation}
Write $z_0=z-\int z\,d\pi$. The spectral theorem gives
\[
 \langle z_0,D^{-1}z_0\rangle
 =\int_0^\infty\langle z_0,P_sz_0\rangle\,ds.
\]
For $0\leq s\leq4/n$, use
$\tanh(ns/8)\geq ns/16$ in
\cref{eq:coordinate-supported-heat}, together with
$\langle z_0,P_sz_0\rangle\leq\langle z,P_sz\rangle$.
This part of the integral is at most $16\norm*{z}_2^2/(nL)$.
Both coordinate gaps are at least $n/2$, so for $s\geq4/n$,
\[
 \langle z_0,P_sz_0\rangle
 \leq e^{-n(s-4/n)/2}\langle z_0,P_{4/n}z_0\rangle
 \leq e^{-n(s-4/n)/2}e^{-L/3}\norm*{z}_2^2,
\]
where $\tanh(1/2)\geq1/3$ was used.
Integrating this bound over $s\geq4/n$ gives
\[
 \frac{2e^{-L/3}}n\|z\|_2^2\leq\frac4{nL}\|z\|_2^2.
\]
Adding the two parts proves the constant
$20$ in \cref{eq:coordinate-support-inverse}.

Since $\mu(F\ne0)\leq m$,
Jensen's inequality gives
\[
 \operatorname{Ent}_\mu(F^2)\geq\log(1/m)\int F^2\,d\mu.
\]
Combining this with \eqref{eq:coordinate-product-lsi} gives
\eqref{eq:coordinate-support-energy}.
\end{proof}

\section{Derivative and moment estimates for extremals}\label{app:bubble}

We establish derivative and moment estimates for the
Sobolev extremals $w_\theta$ used in \cref{lem:iv-modulation}.

\begin{lemma}\label{lem:bubble}
Let $p>1$ and let $n>p$ be an integer. There exist constants
$L,M,N\geq1$ such that, for every parameter $\theta$ with
$|\log\lambda|\leq1$ and every $c,d\in\mathbb R^{n+1}$,
\begin{align}
 E_v^{-1/p}\|\nabla Dw_\theta[c]\|_p&\leq L|c|,
 \label{eq:iv-chart-gradient}\\
 V_q^{-1/q'}\|D(w_\theta^{q-1})[c]\|_{q'}&\leq M|c|,
 \label{eq:iv-chart-first}\\
 V_q^{-1/q'}\|D^2(w_\theta^{q-1})[c,d]\|_{q'}&\leq N|c||d|.
 \label{eq:iv-chart-second}
\end{align}
The maps $\theta\mapsto w_\theta$,
$\theta\mapsto w_\theta^{q-1}$, and
$\theta\mapsto w_\theta^q$ are, respectively,
$C^1$ into $\dot W^{1,p}$, $C^2$ into $L^{q'}$,
and $C^2$ into $L^1$.

The constants $L,M,N$ can be obtained using finitely many
algebraic operations, powers, and Beta integrals with positive
parameters. For each fixed $p>1$, they can be chosen independently
of $n$ for all sufficiently large integers $n$, with a dimension
threshold given by explicit inequalities in $p$.
\end{lemma}

\begin{proof}
Write $c=(c_0,c_z)$ and $d=(d_0,d_z)$.
In the variables $x'=\lambda(x-\eta z)$, set
\[
 e_d=(t-\tau)/\sigma_t,\qquad
 e_a=\frac{\sqrt n}{K}k(t)(a\cdot\omega),\qquad
 \ell_\theta[c]=D\log w_\theta[c]=c_0e_d+\lambda e_{c_z}.
\]
Since $\int t\,d\mu_t=\tau$,
$\operatorname{Var}_{\mu_t}t=\sigma_t^2$, and
$\int\omega_i\omega_j\,d\sigma=\delta_{ij}/n$, we have
\begin{equation}\label{eq:iv-chart-gram}
 \int\ell_\theta[c]\,d\mu_\theta=0,\qquad
 \int\ell_\theta[c]\ell_\theta[d]\,d\mu_\theta
 =c_0d_0+\lambda^2c_z\cdot d_z,
 \quad d\mu_\theta=V_q^{-1}w_\theta^qdx.
\end{equation}
Define $d\nu_\theta=E_v^{-1}|\nabla w_\theta|^p\,dx$. The map
$x\mapsto x'=\lambda(x-\eta z)$ sends $\mu_\theta$ and $\nu_\theta$
to $\mu$ and $\nu$, respectively, as defined in \eqref{eq:local-measures}.

\Needspace{8\baselineskip}
We first choose $L,M,N$ for fixed $n,p$. Set
\begin{align*}
 C_{\log}&=\sigma_t^{-1}+\mathrm e\sqrt n/K,&
 U&=\left[\frac{B(\alpha+2-p,\beta-1)}
                    {B(\alpha+1,\beta-1)}\right]^{1/p},\\
 \xi&=\frac{(p-2)_+}{p},&
 Z&=\left[\frac{B(\alpha-\xi q',\beta)}{B(\alpha,\beta)}
                                              \right]^{1/q'}.
\end{align*}
All the Beta integrals above are finite. Indeed,
$\alpha+2-p>1$ because $n>p$;
for $p>2$, $\xi q'<\xi p'=(p-2)/(p-1)<1<\alpha$, while
$\xi=0$ for $p\leq2$.  Also
$K^2=B(\alpha+2/p,\beta+2\tau)/B(\alpha,\beta)>0$.
By \eqref{eq:ground-gradient}, we obtain
\[
 |\mathsf Ze_d|\leq\frac{p}{(n-p)\sigma_t},\qquad
 |\mathsf Ze_a|\leq\eta(p'+1)t^{-\tau}|a|.
\]
Since both $\|\ell_\theta[c]\|_{L^s(\mu_\theta)}$ and
$\|\ell_\theta[c]\|_{L^s(\nu_\theta)}$ are at most $C_{\log}|c|$
for every $s\geq1$, one can take
\[
 L=\max\left\{1,C_{\log}+\frac{p}{(n-p)\sigma_t}
                           +\mathrm e\eta(p'+1)U\right\}.
\]

For $H_\theta=D_\theta^2\log w_\theta$, the Hessian blocks are
\begin{align*}
 H_{00}&=-\frac{p't(1-t)}{b\sigma_t^2},\\
 H_{0z}[a]&=-\frac{\eta\lambda p'}{\sigma_t}
                       \frac{t(1-t)}{|x'|}(a\cdot\omega),\\
 H_{zz}[a,a']&=-b\eta^2\lambda^2\frac{k}{|x'|}
 \left[a\cdot a'+(p'-2-p't)(a\cdot\omega)(a'\cdot\omega)\right].
\end{align*}
The elementary bounds $t(1-t)/|x'|\leq1$ and
$k/|x'|\leq t^{-\xi}$ give
$\bigl\|\|H_\theta\|_{\mathrm{op}}\bigr\|_{L^{q'}(\mu_\theta)}\le H_*$, where
\[
 H_*:=\frac{p'}{4b\sigma_t^2}+\frac{2\mathrm e\eta p'}{\sigma_t}
              +\mathrm e^2b\eta^2(1+|p'-2|+p')Z.
\]
Thus $M=\max\{1,(q-1)C_{\log}\}$ and
$N=\max\{1,(q-1)H_*+(q-1)^2C_{\log}^2\}$ suffice, using
\[
 D(w_\theta^{q-1})=(q-1)w_\theta^{q-1}\ell_\theta,\qquad
 D^2(w_\theta^{q-1})=w_\theta^{q-1}
           [(q-1)H_\theta+(q-1)^2\ell_\theta\otimes\ell_\theta].
\]

We next show that $L,M,N$ can be chosen independently of
$n$ when $p$ is fixed and $n$ is sufficiently large.
We use moments of the Beta distributions and of the
uniform measure on the sphere. Choose integers $j,m$ with
$2j\geq\max\{p,4p'\}$ and $m\geq p-1$. The laws
$\mu_t=\operatorname{Beta}(\alpha,\beta)$ and
$\nu_t=\operatorname{Beta}(\alpha+1,\beta-1)$ both have total shape
parameter $n$; their means differ by $1/n$. For the law $\mu_t$ with
mean $\tau$, integration by parts gives, for smooth $h$,
\[
\int_0^1(t-\tau)h(t)\,d\mu_t
 =\frac1n\int_0^1t(1-t)h'(t)\,d\mu_t.
\]
Writing $D_j=((2j-1)!!)^{1/(2j)}$ and $k_0=k(\tau)$,
this identity and the sphere moments give, for $1\leq s\leq2j$,
\[
 \|e_d\|_{L^s(\mu)},\ \|e_d\|_{L^s(\nu)}
 \leq\frac{D_j+2}{\sqrt{2d_*}},\qquad
 \|e_a\|_{L^s(\mu)},\ \|e_a\|_{L^s(\nu)}
 \leq\frac{\sqrt8D_j}{k_0}|a|.
\]
Indeed $\mathbb E|t-\mathbb Et|^{2j}\leq(2j-1)!!/(4n)^j$;
$K\geq k_0/\sqrt8$ follows from Chebyshev's inequality for
$n\geq2/\min(\tau,1-\tau)^2$; and
$\mathbb E(\sqrt n\,\omega\cdot a)^{2j}\leq(2j-1)!!|a|^{2j}$.

Furthermore, for $n\tau\geq2m$ and $n\tau\geq2$,
\[
 U\leq(2/\tau)^\tau,\qquad Z\leq(2/\tau)^\xi.
\]
The Beta moment formulas and H\"older's inequality give these bounds;
the second also uses $\int t^{-1}\,d\mu_t\leq2/\tau$.
For $n\geq\max\{2p,2\}$,
\[
 \frac{p}{(n-p)\sigma_t}\leq\frac{2\sqrt2p}{\sqrt{d_*}},\qquad
 \eta\leq\frac{2(p-1)}{K\sqrt n},\qquad
 \frac1{b\sigma_t^2}\leq3p^2,\quad
 \frac\eta\sigma_t\leq\frac{2\sqrt2(p-1)}{K\sqrt{d_*}},\quad
 b\eta^2\leq\frac{2(p-1)}{K^2}.
\]
Combining these bounds with the derivative formulas and using
$q-1\leq2p-1$, $2q'\leq2p'$, and the moment estimates in place
of the pointwise bound $C_{\log}$, we obtain constants $L,M,N$
depending only on $p$, through finitely many explicit operations.
Also, $\gamma\leq1$ when $n\geq8(p-1)^2/d_*$. Taking $n$ above
all the preceding thresholds completes the uniform estimates.

We finally justify the differentiability statements.
When $p>2$, second spatial derivatives may be singular at the
center, but they are $O(r^{p'-2})$ there. The boundary term
in integration by parts is $O(\varepsilon^{n+p'-2})$ and hence
tends to zero as $\varepsilon\downarrow0$.
The positive Beta parameters checked above give the required
$L^{q'}$ and $L^p$ integrability. At infinity, each translation
derivative improves the decay of the bubble, while a dilation
derivative preserves it. Translations and critical dilations
are strongly continuous in the relevant $L^s$ spaces.
Writing their difference quotients as integrals of the
corresponding weak derivatives gives the stated differentiability
in $\dot W^{1,p}$, $L^{q'}$, and $L^1$.
\end{proof}

\section{Nonlinear estimates and weighted quadratic forms}\label{app:nonlinear}

In this appendix, we prove several estimates used in the
local stability arguments. In
Subsection~\ref{subsec:vector-transfer}, we establish
vector and scalar inequalities, including
\eqref{eq:bregman-radial-growth}--\eqref{eq:appendix-scalar-transfer}.
In Subsection~\ref{subsec:form-limits}, we justify weighted
integration by parts and approximation in the closed form domain.
Finally, in Subsection~\ref{sec:shifted}, we prove the gradient
estimates \eqref{eq:shifted-high} and \eqref{eq:shifted-low}
for shifted perturbations.

\subsection{Vector and scalar inequalities}\label{subsec:vector-transfer}

\begin{lemma}\label{lem:transfer}
Fix $p>1$ and a unit vector $e\in\R^n$, and set
$\mathfrak r_p=2\max\{p,p'\}$ and
$\zeta=p/(2\mathfrak r_p)=\frac14\min\{1,p-1\}$.
The following estimates hold.
\begin{enumerate}
\item[(i)] For every $z\in\R^n$ and $\lambda\geq1$,
\begin{equation}\label{eq:bregman-radial-growth}
 z\cdot DG_1(z)\leq\mathfrak r_pG_1(z),\qquad
 G_1(\lambda z)\leq\lambda^{\mathfrak r_p}G_1(z).
\end{equation}

\item[(ii)] For $c>0$, define
\[
 \Phi_c(s)=(s^\zeta-c^\zeta)_+^{1/\zeta},\qquad
 \Psi_c(s)=(s-c)_+-\Phi_c(s),\qquad s\geq0.
\]
Then, for every $s>c$ and $z\in\R^n$,
\begin{equation}\label{eq:appendix-nonlinear-split}
 G_s(z)\geq G_{\Phi_c(s)}(\Phi_c'(s)z)
 +\min\{1,p/2\}G_{c+\Psi_c(s)}(\Psi_c'(s)z).
\end{equation}

\item[(iii)] If $c=M+1$ with $M\geq\mathrm e$,
there exist $R_0,C_{p,M}>0$, depending only on
$p,M$, such that, for every $p<q\leq2p$ and $r>M$,
with $\epsilon=q-p$,
\begin{equation}\label{eq:appendix-scalar-transfer}
 \left(\frac pq+\frac\epsilon{8p}\right)r^q-r^p
 \leq\left(\frac pq+\frac\epsilon{4p}\right)
 \Phi_c(1+r)^q-\Phi_c(1+r)^p
 +C_{p,M}\epsilon1_{\{M<r<R_0\}}.
\end{equation}
\end{enumerate}
\end{lemma}

\begin{proof}
We first prove (i). Write $w=e+z$, $r=\abs*{w}$, and
$u=e\cdot w\in[-r,r]$. Direct calculation gives
\[
 G_1(z)=r^p-pu+p-1,
 \qquad z\cdot DG_1(z)
 =p[r^p+1-(r^{p-2}+1)u].
\]
The difference $\mathfrak r_pG_1-z\cdot DG_1$ is affine in $u$,
so it suffices to consider $u=\pm r$.
For $w=te$, $t\geq0$, set
$B=t^p-pt+p-1$ and $D=p(t-1)(t^{p-1}-1)$.
Convexity of $t^{p-1}$ for $p\geq2$ and concavity of $t^{2-p}$
for $1<p<2$ give, respectively,
\begin{align*}
 pB-D&=p[t^{p-1}-(p-1)t+p-2]\geq0,\\
 p'B-D&=\frac p{p-1}[(2-p)t^p+(p-1)t^{p-1}-t]\geq0.
\end{align*}
For $w=-re$, the corresponding expressions are
$B=r^p+pr+p-1$ and $D=p(r+1)(r^{p-1}+1)$.
The inequalities $r^{p-1}\leq r^p+1$ and
$B\geq\min\{1,p-1\}(r^p+r+1)$ imply $D\leq\mathfrak r_pB$.
The case $w=0$ follows by continuity. Thus the first inequality
in (i) holds.

For $z\ne0$, strict convexity gives $G_1(az)>0$ and
$\frac{d}{da}\log G_1(az)\leq\mathfrak r_p/a$ for $a>0$.
Integrating over $[1,\lambda]$ proves the second inequality
in (i); the case $z=0$ is immediate.

To prove (ii), fix $c>0$ and $s>c$, and set
$\theta=\Phi_c(s)/s\in(0,1)$.
Since $0<\zeta<1$, subadditivity of $x^\zeta$ gives
$0\leq\Phi_c,\Psi_c\leq(s-c)_+$, and differentiation yields
\[
 \Phi_c'(s)=\theta^{1-\zeta},\qquad
 \Psi_c'(s)=1-\theta^{1-\zeta}.
\]
Both functions are $1$-Lipschitz. Homogeneity and
\eqref{eq:bregman-radial-growth} give
\[
 G_{\Phi_c(s)}(\Phi_c'(s)z)
 =s^p\theta^pG_1(\theta^{-\zeta}z/s)
 \leq\theta^{p/2}G_s(z).
\]
Also $c+\Psi_c(s)=s(1-\theta)$ and
$(1-\theta^{1-\zeta})/(1-\theta)\leq1$, so convexity gives
\[
 G_{c+\Psi_c(s)}(\Psi_c'(s)z)
 \leq(1-\theta)^pG_s(z).
\]
Finally,
$1-\theta^{p/2}\geq\min\{1,p/2\}(1-\theta)
\geq\min\{1,p/2\}(1-\theta)^p$.
Combining these estimates proves
\eqref{eq:appendix-nonlinear-split},
including the case $se+z=0$.

For (iii), fix $M\geq\mathrm e$ and $c=M+1$, and let
$p<q\leq2p$, with $\epsilon=q-p$.
Set $A=\Phi_c(1+r)$ and $J_q(x)=(p/q)x^q-x^p$.
Since $0\leq A\leq r-M<r$,
\[
 J_q(r)-J_q(A)
 =p\int_A^r x^{p-1}(x^\epsilon-1)\,dx
 \leq p\epsilon r^{q-1}\log r\,(r-A).
\]
The integrand is nonpositive for $x\leq1$. For $x>1$, we use
$1-e^{-\epsilon\log x}\leq\epsilon\log x$.
Bernoulli's inequality gives $r-A\leq C_{\mathrm{tr}}r^{1-\zeta}$, where
$C_{\mathrm{tr}}=2(M+1)^\zeta/\zeta$.
Choose $R_0=R_0(p,M)\geq\max\{M+1,\mathrm e\}$ so that
\[
 C_{\mathrm{tr}}R_0^{-\zeta}\leq\frac1{8p},\qquad
 \frac{2C_{\mathrm{tr}}}{\mathrm e\zeta}R_0^{-\zeta/2}\leq\frac1{16p^2}.
\]
For $r\geq R_0$, these conditions and
$\log r\leq2r^{\zeta/2}/(\mathrm e\zeta)$ imply
$A^q\geq r^q(1-qC_{\mathrm{tr}}r^{-\zeta})\geq3r^q/4$.
For $r\geq R_0$, the left-hand side of
\eqref{eq:appendix-scalar-transfer} minus its right-hand side
is at most
\[
 \epsilon r^q\left[pC_{\mathrm{tr}}r^{-\zeta}\log r
 +\frac1{8p}-\frac3{16p}\right]
 \leq0\qquad(r\geq R_0).
\]

For $M<r<R_0$, the preceding integral formula bounds the error by
\[
 \epsilon R_0^{2p}\bigl(p\log R_0+1/(8p)\bigr).
\]
Thus we may take
$C_{p,M}=R_0^{2p}(p\log R_0+1/(8p))$.
Both $R_0$ and $C_{p,M}$ are independent of $q$, proving (iii).
\end{proof}

\subsection{Weighted integration by parts and approximation}\label{subsec:form-limits}

\begin{lemma}\label{lem:form-domain}
Let $1<p<n$, where $n$ is an integer.
Every bounded admissible function $b$ has locally
absolutely continuous radial slices and spherical
$W^{1,1}$ slices almost everywhere. Under the
hypotheses of \cref{lem:coordinate}(i), the radial weak
pairing \eqref{eq:coordinate-weak-pairings} and its
spherical counterpart remain valid after integration
over the other coordinate, provided the corresponding
product integrands are absolutely integrable.

For the quadratic-form estimates, suppose in addition
that $\alpha=n\tau>1$ and $\beta=n/p>1$.
Every bounded function with locally weak first derivatives
on $(0,1)\times\mathbb S^{n-1}$ and finite quadratic
energy in \eqref{eq:exact-Z-form} belongs to the associated
closed form domain on $L^2(\mu)$. Consequently, the
product spectral-gap and logarithmic Sobolev inequalities,
as well as the support estimate under its stated support
condition, extend to these functions; see
\eqref{eq:coordinate-spectra},
\eqref{eq:coordinate-product-lsi}, and
\eqref{eq:coordinate-support-energy}.

Finally, suppose $p\geq2$, let $c>0$, and let $w\geq0$
be admissible. For $j=1,2,\ldots$, set
$X_j=(c+\min\{w,j\})^{p/2}-c^{p/2}$. Then
\begin{equation}\label{eq:power-form-bound}
 \int\abs*{\mathsf ZX_j}^2\,d\nu
 \leq\frac{p^2}{4}\int(c+w)^{p-2}
 \abs*{\mathsf Zw}^2\,d\nu.
\end{equation}
If the right-hand side is finite, set
$X=(c+w)^{p/2}-c^{p/2}$. Then
$X_j\to X$ strongly in $L^2(\mu)$ and
$\mathsf ZX_j\to\mathsf ZX$ strongly in $L^2(\nu)$.
In particular, $X$ belongs to the same closed form domain
and satisfies the corresponding form bound.
No support restriction is needed for this limit.
\end{lemma}

\begin{proof}
For bounded admissible functions, slicing on compact annuli gives
locally absolutely continuous radial slices and spherical $W^{1,p}$
slices, hence $W^{1,1}$ slices, almost everywhere. The two Beta
densities give
\[
 \frac{d\mu}{d\nu}=\frac{\alpha}{\beta-1}\frac{1-t}{t}.
\]
By the gradient bounds in
\cref{lem:coordinate}(i) and
\eqref{eq:ground-gradient}, the integrands
in the radial and spherical weak pairings,
when expressed with respect to $d\nu$,
are bounded in absolute value by
$C(n,p)|\mathsf Zb|$.
Since $\nu$ is a probability measure and
$\mathsf Zb\in L^p(\nu)$ with $p>1$,
both integrands are absolutely integrable.
We can therefore apply the radial and spherical
weak integration-by-parts identities, followed
by Fubini's theorem.

In \cref{eq:exact-Z-form}, $\alpha=n\tau>1$ and $\beta=n/p>1$
for the dimensions used in \cref{sec:local}. For a bounded finite-form
function, endpoint cutoffs changing on $(\varepsilon,2\varepsilon)$
add $O(\varepsilon^{\alpha-1})$ radial energy at zero and
$O(\varepsilon^{\beta-1})$ at one. The radial and angular energy on the discarded
strips tends to zero because the total energy
is finite. On the remaining compact cylinder,
mollification in $t$ and smooth approximation
on the sphere give strong convergence in the
form norm.

Thus the function belongs to the closed form
domain, and the spectral-gap and logarithmic
Sobolev inequalities extend to it. Under the
support condition in
\eqref{eq:coordinate-support-energy}, the
support estimate also applies.

For the truncated powers, the chain rule gives
\eqref{eq:power-form-bound}. To prove the convergence
statements, assume that the right-hand side of this
inequality is finite. This holds in particular when
$\mathcal R(c+w)<\infty$, by
\eqref{eq:full-bregman-bound}.

Since $vw\in\dot W^{1,p}(\R^n)$, the Sobolev
inequality gives $w\in L^q(\mu)$.
As $q>p\geq2$ and $\mu$ is a probability measure,
$w\in L^p(\mu)\cap L^2(\mu)$.
The bound $X_j^2\leq C_{p,c}(w^p+w^2)$ and
dominated convergence therefore give
$X_j\to X$ strongly in $L^2(\mu)$.

The chain rule gives
\[
 \mathsf ZX_j=\frac p2(c+w)^{p/2-1}\mathbf1_{\{w<j\}}\mathsf Zw
 \quad\text{a.e.}
\]
The squared difference of the derivatives is dominated by a constant
multiple of $(c+w)^{p-2}|\mathsf Zw|^2$, so
$\mathsf ZX_j\to\mathsf ZX$ strongly in $L^2(\nu)$.
Each bounded $X_j$ belongs to the closed form
domain by the approximation argument above.
The strong convergence in the form norm then
shows that $X$ also belongs to this domain, and
\eqref{eq:power-form-bound} passes to the limit.
\end{proof}

\subsection{Gradient estimates for shifted perturbations}\label{sec:shifted}

\begin{lemma}\label{lem:shifted}
For every $p>1$, there exists $c_p^{\rm sh}>0$,
depending only on $p$, such that the following hold
for every integer $n\geq N_p^{\rm supp}(2)$, every
$c\geq1$, and every nonnegative admissible function
$w$ satisfying $\mu(w\ne0)\leq\mathrm e^{-2}$.

If $p\geq2$, then
\begin{equation}\label{eq:shifted-high}
\mathcal R(c+w)\geq\frac{c_p^{\rm sh}}{n}
 \frac{\norm*{\nabla(vw)}_p^p}{E_v}.
\end{equation}

If $1<p<2$, assume in addition that
$\|w\|_{L^p(\mu)},\|w\|_{L^p(\nu)}$ and
$\|\mathsf Zw\|_{L^p(\nu)}$ are all at most $c$.
Then
\begin{equation}\label{eq:shifted-low}
\mathcal R(c+w)\geq\frac{c_p^{\rm sh}c^{p-2}}n
 \left(\frac{\norm*{\nabla(vw)}_p}{E_v^{1/p}}\right)^2.
\end{equation}
An explicit coefficient for both cases is given in
\eqref{eq:shifted-coefficient}.
\end{lemma}

\begin{proof}
If $\mathcal R(c+w)=\infty$, both inequalities are immediate.
We may therefore assume that $\mathcal R(c+w)<\infty$.
Set $X=\|w\|_{L^p(\nu)}$, $U=\|w\|_{L^p(\mu)}$,
$Y=\|\mathsf Zw\|_{L^p(\nu)}$, and $J_p=p/(1-2^{-p})$.

We first consider $p\geq2$.
By \eqref{eq:full-bregman-bound}, $\mathcal R(c+w)$
controls both integrals
\[
 \int|\mathsf Zw|^p\,d\nu,\qquad
 \int(c+w)^{p-2}|\mathsf Zw|^2\,d\nu.
\]
It follows from \cref{lem:form-domain} that
$X_*=(c+w)^{p/2}-c^{p/2}$ belongs to the
closed form domain. Moreover,
$\mu(X_*\ne0)\leq\mu(w\ne0)\leq e^{-2}$.
Applying \eqref{eq:capacity-high} to $X_*$ with
$L_0=2$, and using \eqref{eq:power-form-bound}
and $X_*^2\geq w^p$, we obtain
\[
\int w^p\,d\mu\leq C_pn\mathcal R(c+w),\qquad
\int\abs*{\mathsf Zw}^p\,d\nu\leq C_p\mathcal R(c+w).
\]
By \eqref{eq:local-picone-signed} and Young's inequality,
$X^p\leq C_p(U^p+Y^p)$.
Combining this with \eqref{eq:normalized-gradient-product}
and the estimates above proves \eqref{eq:shifted-high},
with a coefficient independent of $c$.

We now consider $1<p<2$ and assume the additional
$L^p$ bounds in \cref{lem:shifted}.
Put $a_j=2^jc$, $\chi_j=(c+w-a_j)_+\wedge a_j$, and
$D_j=\{a_j<c+w<2a_j\}$ for $j=0,1,2,\ldots$.
The Lipschitz truncation $\chi_j$ is admissible,
$0\le\chi_j\le a_j$, and $\{\chi_j\ne0\}\subseteq\{w\ne0\}$.
Its gradient support $D_{\chi_j}=\{\mathsf Z\chi_j\ne0\}$ is
contained in $D_j$ up to null sets, and $\mathsf Z\chi_j=\mathsf Zw$
almost everywhere on $D_j$. On $D_{\chi_j}$, the function $f=c+w$
satisfies $0\le f\le2a_j$ and $\mathsf Z\chi_j=\mathsf Zf$.
Thus \eqref{eq:capacity-low} applies with $L_0=2$; the same
coefficient $\ell_p$ works at every level.
The sets $D_{\chi_j}$ are pairwise disjoint up to null sets,
so we may sum \eqref{eq:capacity-low} over $j$.

If $w\leq c$, then $\chi_0=w$ and the zeroth term satisfies
$c^{p-2}w^2\geq(c+w)^{p-2}w^2$. If $w>c$, choose $j$ with
$2a_j\leq c+w<4a_j$.
Then $\chi_j=a_j$, and
\[
 a_j^p\geq4^{-p}(c+w)^{p-2}w^2.
\]
Tonelli therefore gives
\[
J:=\int(c+w)^{p-2}w^2\,d\mu\leq C_pn\mathcal R(c+w).
\]

H\"older's inequality with exponents $2/p$ and $2/(2-p)$ gives
\begin{align*}
\norm*{w}_{L^p(\mu)}^2
 &\leq J(c+\norm*{w}_{L^p(\mu)})^{2-p},\\
\norm*{\mathsf Zw}_{L^p(\nu)}^2
 &\leq C_p\mathcal R(c+w)
   (c+\norm*{w}_{L^p(\nu)}+\norm*{\mathsf Zw}_{L^p(\nu)})^{2-p}.
\end{align*}
The second line follows from \eqref{eq:full-bregman-bound}, with weight
$(c+w+\abs*{\mathsf Zw})^{p-2}$. It does not require
$\mathsf Zw\in L^2(\nu)$.
By \eqref{eq:local-picone-signed},
$X^p\leq U^p+pX^{p-1}Y$.
If $X\geq2U$, then
$(1-2^{-p})X\leq pY$; otherwise, $X<2U$.
Thus in either case $X\leq2U+J_pY$.
By \eqref{eq:normalized-gradient-product},
\[
 \frac{\|\nabla(vw)\|_p}{E_v^{1/p}}
 =\|we+\mathsf Zw\|_{L^p(\nu)}
 \le X+Y\le2U+(J_p+1)Y.
\]
Together with the preceding estimates, this proves
\eqref{eq:shifted-low}.

We now give an explicit choice of the constant.
Let $\ell_p$ and $\gamma_p$ be the constants in
\cref{lem:capacity} and \cref{lem:eff-bregman}, respectively.
The estimates above give
$U^{r_p}\leq C_U n\mathcal R(c+w)$ and
$Y^{r_p}\leq C_Y\mathcal R(c+w)$, where
\[
 (C_U,C_Y)=
 \begin{cases}
 (p^2/(8\gamma_p\ell_p),\gamma_p^{-1}),&p\geq2,\\
 (4^p(2c)^{2-p}/(2\ell_p),(3c)^{2-p}/\gamma_p),&p<2.
 \end{cases}
\]
The preceding gradient bound gives a common sufficient coefficient:
\begin{equation}\label{eq:shifted-coefficient}
 k_{\rm sh}=\left(2^{r_p-1}
       [2^{r_p}C_U+(J_p+1)^{r_p}C_Y]\right)^{-1}>0,
 \qquad
 \mathcal R(c+w)\geq\frac{k_{\rm sh}}n
       \left(\frac{\|\nabla(vw)\|_p}{E_v^{1/p}}\right)^{r_p}.
\end{equation}
For $p\geq2$, the coefficient $k_{\rm sh}$ depends only on $p$.
For $1<p<2$, under the stated norm assumptions, the formulas for
$C_U$ and $C_Y$ show that $k_{\rm sh}$ is $c^{p-2}$ times a positive
explicit constant depending only on $p$.
Thus the stated coefficient $c_p^{\rm sh}$ can be chosen to depend
only on $p$ in both cases.
\end{proof}

\Needspace{12\baselineskip}
\section*{Declarations}
\addcontentsline{toc}{section}{Declarations}

\noindent\textbf{Conflict of interest.}
The authors declare that they have no conflicts of interest.

\medskip

\noindent\textbf{Data availability.}
Data sharing is not applicable to this article, as no datasets were generated or analyzed during the current study.

\medskip

\noindent\textbf{Generative AI disclosure.}
The authors initiated this research and led it throughout its development. ChatGPT was used to assist in discussing certain mathematical details, checking mathematical arguments and references, and preparing the manuscript. The authors have carefully reviewed the AI-assisted material, further refined the exposition, and take full responsibility for the mathematical results and the final manuscript.

\end{document}